\pdfoutput=1
\documentclass{amsart}
\usepackage{calc}
\usepackage[nomessages]{fp}
\usepackage{tikz,amsmath,amssymb}
\usepackage{pgfplots}
\usepackage{quiver}
\usetikzlibrary{matrix,arrows,backgrounds,shapes.misc,shapes.geometric,patterns,calc,positioning,intersections}
\usetikzlibrary{decorations.pathmorphing, angles, quotes}
\pgfplotsset{compat=1.15}
\usepackage{mathrsfs}
\usepackage{rotating}
\usepackage{tkz-euclide} 
\usetikzlibrary{decorations.markings}
\usetikzlibrary{shapes,fit}
\usepackage[english]{babel}
\usepackage[hidelinks]{hyperref}
\usepackage{amsmath,amsfonts,amssymb}
\usepackage{amsthm}
\usepackage[all]{xy}
\usepackage{ stmaryrd }
\usepackage{ mathrsfs }
\theoremstyle{plain} 

\usepackage{epsfig}
\usepackage[T1]{fontenc}
\usepackage{lmodern}
\usepackage{eurosym}
\usepackage{graphicx}
\usepackage{enumitem}
\usepackage{amsbsy}
\usepackage{rotating}
\usepackage[noabbrev,capitalise]{cleveref}
\usepackage{silence}
\usepgfplotslibrary{fillbetween}
\makeatletter
\tikzset{use path/.code=\tikz@addmode{\pgfsyssoftpath@setcurrentpath#1}}
\makeatother
\pgfdeclarelayer{bg}
\pgfsetlayers{bg,main}

\usepackage{stackengine,scalerel}
\stackMath
\def\Resleq{\ThisStyle{\mathrel{%
  \stackinset{r}{.75pt+.15\LMpt}{t}{.1\LMpt}{\rule{.3pt}{1.1\LMex+.2ex}}{\SavedStyle\leqslant}%
}}}

\definecolor{midnightblue}{rgb}{0.1, 0.1, 0.44}
\definecolor{plum}{rgb}{0.56, 0.27, 0.52}
\definecolor{Plum}{rgb}{0.56, 0.27, 0.52}
\definecolor{patriarch}{rgb}{0.5, 0.0, 0.5}
\definecolor{darkgreen}{rgb}{0.0, 0.2, 0.13}
\definecolor{darkcerulean}{rgb}{0.03, 0.27, 0.49}
\definecolor{jade}{rgb}{0.0, 0.66, 0.42}

\makeatletter
\newcommand{\addbar}[3]{{\vphantom{#3}\mathpalette\add@bar{{#1}{#2}{#3}}}}
\newcommand{\add@bar}[2]{\add@@bar{#1}#2}
\newcommand{\add@@bar}[4]{%
  \begingroup
  \sbox\z@{$\m@th#1#4$}%
  \ooalign{%
    \hidewidth\kern#2\wd\z@\add@@@bar{#1}{#3}\hidewidth\cr
    \box\z@\cr
  }%
  \endgroup
}
\newcommand{\add@@@bar}[2]{%
  \sbox\tw@{$\m@th#1\newmcodes@\if\relax#2\relax-\else\bm{-}\fi$}%
  \raisebox{\dimexpr(\ht\z@-\ht\tw@)/2}{\usebox\tw@}%
}
\makeatother

\hypersetup{
	colorlinks=true,
	linkcolor=darkcerulean,
	citecolor=jade,
	filecolor=magenta,      
	urlcolor=cyan,
	pdfpagemode=FullScreen,
}

\renewcommand{\mod}{\operatorname{mod}}

\newcommand{\rep}{\operatorname{rep}}
\newcommand{\add}{\operatorname{add}}
\newcommand{\Ker}{\operatorname{Ker}}

\newcommand{\Hom}{\operatorname{Hom}}

\newcommand{\Ext}{\operatorname{Ext}}

\newcommand{\ind}{\operatorname{ind}}
\newcommand{\proj}{\operatorname{proj}}

\newcommand{\Res}{\operatorname{Res}}
\newcommand{\ResOrd}{\operatorname{\pmb{Res}}}

\newcommand{\Homleq}{\rightarrow}

\newcommand{\THomleq}{\pmb{\pmb{\rightarrow}}}

\newcommand{\Supp}{\operatorname{Supp}}
\newcommand{\opQ}{\operatorname{\mathbf{Q}}}
\newcommand{\opR}{\operatorname{\mathbf{R}}}

\newcommand{\opResAc}{\operatorname{\pmb{\mathscr{R}}}}
\newcommand{\Prj}{\operatorname{Prj}}
\newcommand{\MM}{\operatorname{M}}
\newcommand{\NP}{\operatorname{N}_{\proj}}
\newcommand{\Surf}{\operatorname{\pmb{\mathcal{S}}}}
\newcommand{\PP}{\operatorname{P}}
\newcommand{\tc}{t_{\operatorname{cell}}}
\renewcommand{\sc}{s_{\operatorname{cell}}}
\renewcommand{\Ker}{\operatorname{Ker}}

\newcommand{\ArExt}{\operatorname{ArExt}}
\newcommand{\OvExt}{\operatorname{OvExt}}
\newcommand{\Accord}{\operatorname{\pmb{\mathscr{A}}}}
\newcommand{\gldim}{\operatorname{gldim}}
\newcommand{\pdim}{\operatorname{pdim}}
\newcommand{\AR}{\operatorname{AR}}

\date{\today}

\usepackage{thmtools}
\declaretheorem[numberwithin=section,name=Theorem,
refname={Theorem,Theorems},
Refname={Theorem,Theorems}]{theorem}
\declaretheorem[numberlike=theorem,name=Lemma,
refname={Lemma,Lemmas},
Refname={Lemma,Lemmas}]{lemma}
\declaretheorem[numberlike=theorem,name=Proposition,
refname={Proposition,Propositions},
Refname={Proposition,Propositions}]{prop}
\declaretheorem[numberlike=theorem,name=Corollary,
refname={Corollary,Corollaries},
Refname={Corollary,Corollaries}]{cor}

\declaretheorem[style=definition,numberlike=theorem,name=Definition,
refname={Definition,Definitions},
Refname={Definition,Definitions}]{definition}
\declaretheorem[style=definition,numberlike=theorem,name=Convention,
refname={Convention,Conventions},
Refname={Convention,Conventions}]{conv}
\declaretheorem[style=definition,numberlike=theorem,name=Algorithm,
refname={Algorithm,Algorithms},
Refname={Algorithm,Algorithms}]{algo}

\declaretheorem[style=definition,numberlike=theorem,name=Example,
refname={Example,Examples},
Refname={Example,Examples}]{ex}

\declaretheorem[style=remark,numberlike=theorem,name=Remark,
refname={Remark,Remarks},
Refname={Remark,Remarks}]{remark}
\declaretheorem[style=remark,numberlike=theorem,name=Notation,
refname={Notation,Notations},
Refname={Notation,Notations}]{Notation}

\usepackage{todonotes}

\definecolor{aquamarine}{rgb}{0.5, 1.0, 0.83}

\newcommand{\new}[1]{\textit{\textbf{\color{patriarch}{#1}}}}

\definecolor{dark-green}{RGB}{14,150,2}
\newcommand{\gpoint}{{\color{dark-green}{\circ}}}
\newcommand{\rpoint}{\textcolor{red}{\bullet}}
\definecolor{darkgray}{rgb}{0.66, 0.66, 0.66}
\definecolor{darkpink}{rgb}{0.91, 0.33, 0.5}
\newcommand{\gsquare}{\color{dark-green}{\pmb{\square}}}
\newcommand{\rsquare}{\color{red}{{\blacksquare}}}
\newcommand{\osquare}{\color{orange}{\pmb{\boxtimes}}}
\newcommand{\psquare}{\color{darkpink}{\pmb{\otimes}}}

\definecolor{mypurple}{rgb}{0.63, 0.36, 0.94}

\begin{document}
\title{Resolving subcategories for gentle algebras I: Monogeneous resolving subcategories for gentle trees}

\author[B.~Dequêne]{Benjamin Dequêne}
\address[B.~Dequêne]{School of Mathematics, University of Leeds, Woodhouse, LS2 9JT, Leeds, United Kingdom}
\email{b.d.Dequene@leeds.ac.uk}

\author[M.~Schoonheere]{Michaël Schoonheere
}
\address[M.~Schoonheere]{LAMFA, Université de Picardie Jules Verne, 33 rue Saint-Leu, 80000, Amiens, France.}
\email{michael.schoonheere@u-picardie.fr}

\begin{abstract}
		This paper is the first part of a series that intends to study resolving subcategories for gentle algebras over an algebraically closed field $\mathbb{K}$. 
        
        In a general setting, we improve the precision of Takahashi's algorithm for resolving closure calculations in well-behaved abelian categories. Then, we modify the geometric model of Baur--Coelho-Simões and Opper--Plamondon--Schroll to compute such subcategories for gentle quivers that have a finite global dimension. 
        
        Finally, we focus on gentle quivers $(Q,R)$ such that $Q$ is a directed tree, and we study monogeneous resolving subcategories, which are the ones generated by a single non-projective indecomposable $\mathbb{K}Q/\langle R \rangle$-module. Moreover, we prove that these subcategories are the join-irreducible elements of the poset of all resolving subcategories ordered by inclusion. 
	\end{abstract}

\maketitle

\tableofcontents

	\section{Introduction}
	\label{sec:Intro}
	\pagestyle{plain}
 The notion of a resolving subcategory of a module category dates at least to \cite{Auslander1969}. This notion refers to subcategories containing the projective modules,  closed under direct sums, summands, extensions, and kernels of epimorphisms. This definition shows that the class of resolving subcategories generalizes the subcategory of objects of Gorenstein dimension 0. The dual notion of a coresolving subcategory generalizes the subcategory of Gorenstein injective objects. Biresolving subcategories, that is subcategories that are both resolving and coresolving, are a generalization of the subcategory of projective-injective objects. Inspired by the case of Frobenius categories, \cite{Rump2021} states that, given a biresolving subcategory $\mathscr{R}$ of an exact category $\mathscr{E}$, the localization of the additive quotient $\mathscr{E}/(\mathscr{R})$ at its morphisms that are both monomorphisms and epimorphisms is triangulated.

Another natural occurrence of (contravariantly finite) resolving subcategories arises from hereditary cotorsion pairs \cite{Auslander1991, Adachi2022}. Hereditary cotorsion pairs are in bijection with contravariantly finite resolving subcategories, and this bijection preserves the natural order given by inclusion. This one-to-one correspondence motivates the computation of all contravariantly finite resolving subcategories. Takahashi \cite{T09} addressed this question, and gave an algorithm to generate the smallest resolving subcategory $\mathscr{R}$ containing a family of objects $\mathcal{X}$. The subcategory $\mathscr{R}$ is called the \new{resolving closure} of $\mathcal{X}$ and denoted by $\Res(\mathcal{X})$. This algorithm is designed to work over every algebra of finite representation type and thus does not exploit any insights into the algebra's combinatorics.

Gentle algebras are algebras arising from quivers with relations and appeared for the first time in \cite{Assem1981,Assem1982}. They have a good combinatorial behavior: almost-split sequences between indecomposable objects have at most two indecomposable direct summands \cite{BS80,WW85,GP68,SW83,DS87}, homspaces and extensions are endowed with a basis parametrized by string and band combinatorics \cite{CB89,K88,CPS21}. Moreover, the module category over the path algebra, the derived category, and the two-term homotopy category of projectives are in bijection with dissected surfaces in several geometric models \cite{BCS21,OPS18,APS19,PPP18,PPP182}, which provide a combinatorial way to compute in these categories.

We aim to compute the lattice of the resolving subcategories of $\mod{A}$ for a gentle algebra $A$ of finite representation type. We first focus on the case of gentle algebras arising from an unpunctured disc, or equivalently, those gentle algebras whose associated
quiver is a tree. We call \new{monogeneous} a resolving subcategory of the form $\Res(X)$ for $X$ a non-projective indecomposable object. Our first result concerns the structure of the lattice of resolving subcategories.

\begin{theorem}[\ref{thm:Monoareallthejoinirred}]
The join-irreducible elements in the poset of resolving subcategories in $\rep(Q,R)$ ordered by inclusion are exactly given by the monogeneous resolving subcategories.
\end{theorem}

To introduce the main result of this work, we need to look more closely at the combinatorics of surfaces, which helps compute extensions and kernels. The dissected disc is endowed with a set of red points on the boundary: one in each cell of the dissected surface. These points are the endpoints of the \emph{accordions} that are curves in bijection with indecomposable modules (\cref{thm:GeomandRep}). Given such an accordion $\delta$, we denote by $\MM(\delta)$ its corresponding indecomposable module. To each accordion $\delta$ we associate the set $\NP(\delta)$ (\cref{def:ExtprojCover}) of accordions in bijection with projective modules appearing in the projective resolution of $\MM(\delta)$ or having a non-trivial extension with $\MM(\delta)$. The set of endpoints of curves in $\NP(\delta)$ is colored according to a combinatorial rule explained in \cref{def:colourendpoints} as red, orange, or green. We define $\opResAc(\delta)$ as the set of accordions whose sources are coloured either red or orange and targets are coloured either orange or green. We can now state the main result:

\begin{theorem}
Let $(Q,R)$ be a gentle tree, and let $\Surf(Q,R) = (\pmb{\Sigma}, \mathcal{M}, \Delta^{\gpoint})$ be the dissected surface associated with $(Q,R)$. For any accordion $\delta$ the monogeneous resolving closure  $\Res(\MM(\delta))$ is the additive closure of the union of indecomposable projectives and  of $\{\MM(\varsigma) \mid \varsigma \in \opResAc(\delta)\}$.
\end{theorem}

The lattice of resolving subcategories is not semidistributive in general. However, in the follow-up work \cite{DS252} we introduce a specific decomposition of all resolving subcategories as a union of monogeneous resolving subcategories. This decomposition allows us to combinatorially describe all resolving closures and, thus, all resolving subcategories. 
 We also consider the bijection between resolving subcategories and tilting modules from \cite{Auslander1991} and give a combinatorial construction of these tilting modules. 

The paper is organized as follows. \cref{sec:Res} contains generalities about resolving subcategories.  \cref{sec:Gentle} recalls what is needed about gentle algebras arising from trees for the study of resolving subcategories.  With these tools, we can begin studying resolving closures and the resolving poset in \cref{sec:TreePart1}. \cref{sec:Surface} gives a translation of the previous results in the surface model of gentle algebras. Then we compute monogeneous resolving subcategories using the geometric model in \cref{sec:TreePart12}.

	\section{Resolving subcategories}
	\label{sec:Res}
	\pagestyle{plain}
This section is organized as follows. In \cref{ss:ImproAlgo} we improve the description of the steps of Takahashi's algorithm used to compute the resolving closure using the Krull-Schmidt property of $\mathscr{C}$. We will show that it is possible to refine the description of the different steps of this algorithm by adding a few hypotheses, and we do so in \cref{ss:UsabAlgo}.
\subsection{General context and definitions}
\label{ss:Gencontext}
Let $\mathscr{C}$ be a Krull--Schmidt abelian category. 
\begin{conv} \label{conv:reducedepi}
A subcategory is said to be \emph{additive}, or \emph{additively closed}, if it is closed under direct sums and summands. Consider $\mathscr{C}$ to be an additive category and $\mathcal{X}$ a class of its objects; we denote by $\add(\mathcal{X})$ the smallest additive subcategory of $\mathscr{C}$ containing $\mathcal{X}$.
\end{conv}

Given $A,B \in \mathscr{C}$ and an epimorphism $\begin{tikzcd}
	f: A & B
	\arrow[two heads, from=1-1, to=1-2]
\end{tikzcd}$, we recall that $f$ is \emph{minimal} whenever there are no summand in $\Ker(f)$ that is isomorphic to one of $A$. Assume that $\mathscr{C}$ \new{has enough projective objects}: for any object $M \in \mathscr{C}$, there exist a projective object $P \in \mathscr{C}$ and an epimorphism $\begin{tikzcd}
	P & M
	\arrow[two heads, from=1-1, to=1-2]
\end{tikzcd}.$ 
If such an epimorphism is minimal, we say that $P$ is a \new{projective cover} of $M$. We recall that a projective cover of a given object of $\mathscr{C}$ is unique up to isomorphism. A \new{minimal projective resolution} of $M$ is a complex \[\begin{tikzcd}
	(\PP_M^\ast): & \cdots & \PP_M^2 & \PP_M^1 & \PP_M^0 & M
	\arrow["h_3", from=1-2, to=1-3]
	\arrow["h_2", from=1-3, to=1-4]
	\arrow["h_1", from=1-4, to=1-5]
	\arrow["h_0",two heads, from=1-5, to=1-6]
\end{tikzcd} \] such that:
\begin{enumerate}[label=$\bullet$,itemsep=1mm]
   \item $\PP_M^0$ is the projective cover of $M$ and $h_0$ is an epimorphism; and,
    \item for all $i \in \mathbb{N}$, $\PP_M^{i+1}$ is the projective cover of $\Ker(h_i)$, and $h_{i+1}$ is the composition of the inclusion map $\begin{tikzcd}
	\Ker(h_i) & \PP_M^i
	\arrow[tail,from=1-1, to=1-2]
\end{tikzcd}$ followed by an epimorphism defining $\PP_M^{i+1}$.
\end{enumerate}
From $(\PP_M^\ast)$ as defined as above, for any $i \geqslant 1$, we can define the \new{$i$th syzygy object} of $M$ to be the image of $h_i$, or, equivalently, the kernel of $h_{i-1}$. 

For any subcategory $\mathscr{D} \subseteq \mathscr{C}$, we write $\ind(\mathscr{D})$ for the set of isomorphism classes of indecomposable objects in $\mathscr{D}$. 
We set $\pmb{\ind \setminus \proj}(\mathscr{D}) = \ind(\mathscr{D}) \setminus \proj(\mathscr{C})$.

Let $\mathscr{G} \subseteq \mathscr{C}$ be a subcategory of $\mathscr{C}$. We say that an object $M \in \mathscr{C}$ \new{admits a cover in} $\mathscr{G}$ if there exists an epimorphism $\begin{tikzcd}
	G & M
	\arrow[two heads, from=1-1, to=1-2]
\end{tikzcd}$ with $G \in \mathscr{G}$. We say that the subcategory $\mathscr{G} \subseteq \mathscr{C}$  \new{generates} $\mathscr{C}$ whenever every object in $\mathscr{C}$ admits a cover in $\mathscr{G}$. For instance, $\proj(\mathscr{C})$ generates $\mathscr{C}$ by hypothesis. 

\begin{definition} \label{def:resolv}
A full subcategory $\mathscr{R} \subseteq \mathscr{C}$ is called \new{resolving} if it satisfies the following conditions: 
\begin{enumerate}[label=$(\mathsf{R \arabic*})$, itemsep=1mm]
\setcounter{enumi}{-1}
\item \label{R0} $\mathscr{R}$ is additively closed,
\item \label{R1} $\mathscr{R}$ is generating $\mathscr{C}$,
\item \label{R2} $\mathscr{R}$ is closed under extensions, and,
\item \label{R3} $\mathscr{R}$ is closed under kernels of epimorphisms.
\end{enumerate}
\end{definition}

\begin{remark}
\label{rem:iso_clo}
Note that resolving subcategories are closed under isomorphism.
\end{remark}

We can simplify some conditions to check that a subcategory is resolving under the hypotheses we enforce on $\mathscr{C}$.

\begin{lemma} \label{lem:othercharactresolv}
Let $\mathscr{R} \subseteq \mathscr{C}$ be a full additive subcategory. We have the following equivalences:
\begin{enumerate}[label=$(\roman*)$, itemsep=1mm]
\item \label{ires} $\mathscr{R}$ satisfies \ref{R1} if, and only if, $\proj(\mathscr{C}) \subseteq \mathscr{R}$;
\item \cite[Lemma 3.2]{Y05} $\mathscr{R}$ is resolving if, and only if,  $\mathscr{R}$ satisfies \ref{R0}, \ref{R1}, \ref{R2}, and $\mathscr{R}$ is closed under syzygies.
\end{enumerate}
\end{lemma}

\begin{ex} \label{ex:1stResCats} $ $
\begin{enumerate}[label=$\bullet$, itemsep=1mm]
    \item The largest resolving subcategory of $\mathscr{C}$ is $\mathscr{C}$ itself.
    \item The smallest resolving subcategory of $\mathscr{C}$ is $\proj(\mathscr{C})$ since $\mathscr{C}$ is supposed to have enough projectives.
    \item Given $\{\mathscr{R}_i\}_{i \in I}$ an arbitrary family of resolving subcategories of $\mathscr{C}$, using \cref{lem:othercharactresolv} \ref{ires}, the intersection $\bigcap_{i \in I} \mathscr{R}_i$ contains all the projectives. Thus $\bigcap_{i \in I} \mathscr{R}_i$ is a resolving subcategory. \qedhere
\end{enumerate}
\end{ex}

\subsection{Resolving closure}
\label{ss:ResClos}
The last example allows us to introduce the following notion.

\begin{definition} \label{def:resclos}
Let $\mathcal{X}$ be a set of objects of $\mathscr{C}$. The \new{resolving closure} of $\mathcal{X}$, denoted by $\Res(\mathcal{X})$, is the smallest resolving subcategory of $\mathscr{C}$ containing $\mathcal{X}$. If $\mathcal{X} = \{X\}$ for some object $X\in\mathscr{C}$, we will write $\Res(X)$ for $\Res(\mathcal{X})$.

A resolving subcategory $\mathscr{R} \subseteq \mathscr{C}$ is said to be \new{monogeneous} if there exists a non-projective indecomposable object $X \in \mathscr{C}$ such that $\mathscr{R}$ is equal to $\Res(X)$.
\end{definition}

The following result is a direct consequence of \cref{def:resolv} \ref{R0} and the Krull--Schmidt property we imposed on $\mathscr{C}$. 

\begin{cor} \label{cor:resclos}  Every resolving subcategory of $\mathscr{C}$ is the resolving closure of its non-projective indecomposable objects. 
\end{cor}

Therefore, from now on, we focus exclusively on resolving closures of subsets of $\pmb{\ind \setminus \proj}(\mathscr{C})$.

\begin{lemma} \label{lem:cloresop}
Let $\mathcal{X}$ and $\mathcal{Y}$ be two collections of objects in $\mathscr{C}$. We have the following properties:
\begin{enumerate}[label=$(\roman*)$, itemsep=1mm]
\item $\mathcal{X} \subseteq \Res(\mathcal{X})$,
\item $\Res(\Res(\mathcal{X})) = \Res(\mathcal{X})$, and,
\item if $\mathcal{X} \subseteq \mathcal{Y}$, then $\Res(\mathcal{X}) \subseteq \Res(\mathcal{Y})$.
\end{enumerate}
\end{lemma}

The previous lemma shows that $\Res$ behave like a \emph{closure operator}, whose abstract definition is given below.

\begin{definition} \label{def:closureoperator}
Let $\mathfrak{A}$ be a set. A \new{closure operator} on $\mathfrak{A}$ is a map \[\operatorname{cl}: \mathcal{P}(\mathfrak{A}) \rightarrow \mathcal{P}(\mathfrak{A})\] which satisfies the following properties:
\begin{enumerate}[label=$\bullet$, itemsep=1mm]
    \item for all $\mathcal{X} \in  \mathcal{P}(\mathfrak{A})$, $\mathcal{X} \subseteq \operatorname{cl}(\mathcal{X})$,
    \item for all $\mathcal{X} \in  \mathcal{P}(\mathfrak{A})$, $\operatorname{cl}(\mathcal{X}) = \operatorname{cl}(\operatorname{cl}(\mathcal{X}))$, and,
    \item for all $\mathcal{X}, \mathcal{Y} \in \mathcal{P}(\mathfrak{A})$, if $\mathcal{X} \subseteq \mathcal{Y}$, then $\operatorname{cl}(\mathcal{X}) \subseteq \operatorname{cl}(\mathcal{Y})$.
\end{enumerate}
\end{definition}

We focus on the operator $(-)^{\Res}$ on subsets of $\pmb{\ind \setminus \proj} (\mathscr{C})$ which to every $\mathcal{X} \subseteq \pmb{\ind \setminus \proj }(\mathscr{C})$ associates the class of non-projective indecomposable objects of the subcategory $\Res(\mathcal{X})$. The result below follows from \cref{lem:cloresop} and the Krull--Schmidt property. 

\begin{lemma} \label{lem:resclosop}
The operator $(-)^{\Res}$ is a closure operator on $\pmb{\ind \setminus \proj} (\mathscr{C})$ .
\end{lemma}

\subsection{Algorithm for explicit calculations of resolving closures}
\label{ss:ImproAlgo}
Given a commutative noetherian ring $R$, Takahashi introduces an inductive construction of the resolving closure of any collection of objects in $\mathscr{C} = \mod{R}$ \cite{T09}. We realized that this construction can be used in the broader setting of Krull-Schmidt abelian categories with enough projectives. Moreover, by \cref{cor:resclos}, we improve this algorithm by simplifying each step of the computation. 

\begin{theorem} \label{thm:rescondonindec} An additive subcategory $\mathscr{R} \subseteq \mathscr{C}$ is resolving if and only if all the following properties are satisfied:
\begin{enumerate}[label=$(\mathsf{R' \arabic*})$,itemsep=1mm]
\item \label{R'1} $\proj(\mathscr{C}) \subseteq \mathscr{R}$,
\item \label{R'2} For any short exact sequence \[\begin{tikzcd}
	X & E & Y
	\arrow["k",tail, from=1-1, to=1-2]
	\arrow["f",two heads, from=1-2, to=1-3]
\end{tikzcd},\] if $X,Y \in \ind(\mathscr{R})$, then $E \in \mathscr{R}$, and,
\item \label{R'3} For any short exact sequence \[\begin{tikzcd}
	K & M & Y
	\arrow["k",tail, from=1-1, to=1-2]
	\arrow["f",two heads, from=1-2, to=1-3]
\end{tikzcd},\] if $M \in \mathscr{R}$ and $Y \in \ind(\mathscr{R})$, then $K \in \mathscr{R}$.
\end{enumerate}
\end{theorem}

\begin{prop}\label{prop:red_ker}
Let  $\mathscr{R}$  be an additive subcategory of $\mathscr{C}$. Then $\mathscr{R}$ satisfies \ref{R3} if, and only if, $\mathscr{R}$ satisfies \ref{R'3}. 
\end{prop}
\begin{proof} If $\mathscr{R}$ satisfies \ref{R3}, then it is obvious that $\mathscr{R}$ satisfies \ref{R'3}. 
Assume that $\mathscr{R}$ satisfies \ref{R'3}. We will prove that if there is a short exact sequence 
\[\begin{tikzcd}
	K & M & N
	\arrow["k",tail, from=1-1, to=1-2]
	\arrow["f",two heads, from=1-2, to=1-3]
\end{tikzcd}\]
where $M,N \in \mathscr{R}$ and $K \in  \mathscr{C}$, then $K \in \mathscr{R}$. Using the Krull--Schmidt property of $\mathscr{C}$, we can, proceed by induction on the number of indecomposable summands of $N$. The initialization is obvious by \ref{R'3}.

Now assume that $N$ is not indecomposable. In this case, we write $N=N_1\oplus N'$ with $N_1 \in \ind(\mathscr{C})$. As $\mathscr{R}$ is additive, we have that $N_1, N' \in \mathscr{R}$. Therefore, we have the following short exact sequences. \[\begin{tikzcd}
	K & M & N
	\arrow["k",tail, from=1-1, to=1-2]
	\arrow["f",two heads, from=1-2, to=1-3]
\end{tikzcd} \text{ and } \begin{tikzcd}
	N_1 & N & N'
	\arrow["i",tail, from=1-1, to=1-2]
	\arrow["p",two heads, from=1-2, to=1-3]
\end{tikzcd}\]
By taking the pullback of $f$ along $i$, we obtain the following diagram of short exact sequences.
\[\begin{tikzcd}
	K & {M'} & {N_1} \\
	K & M & N \\
	& {N'} & {N'}
	\arrow[tail, from=1-1, to=1-2]
	\arrow[Rightarrow, no head, from=1-1, to=2-1]
	\arrow["{f_1}", two heads, from=1-2, to=1-3]
	\arrow[tail, from=1-2, to=2-2]
	\arrow["{i}", tail, from=1-3, to=2-3]
	\arrow["k"', tail, from=2-1, to=2-2]
	\arrow["f"', two heads, from=2-2, to=2-3]
	\arrow["{f'}"', two heads, from=2-2, to=3-2]
	\arrow["{p}",two heads, from=2-3, to=3-3]
	\arrow[Rightarrow, no head, from=3-2, to=3-3]
\end{tikzcd}\]
 We already know that $N_1 \in \ind(\mathscr{R})$. Moreover, \[\begin{tikzcd}
	M' & M & N'
	\arrow[tail, from=1-1, to=1-2]
	\arrow["f'",two heads, from=1-2, to=1-3]
\end{tikzcd}\] is a short exact sequence, where $M,N' \in \mathscr{R}$, and the number of indecomposable summands of $N'$ is strictly smaller than the number of indecomposable summands of $N$. By induction hypothesis, we have that $M' \in \mathscr{R}$. Therefore, by \ref{R'3}, we have that $K \in \mathscr{R}$ on the exact sequence  \[\begin{tikzcd}
	K & M' & N_1.
	\arrow[tail, from=1-1, to=1-2]
	\arrow["f_1",two heads, from=1-2, to=1-3]
\end{tikzcd}\] This ends the induction.
\end{proof}
\begin{prop} \label{prop:red_ext}
Let  $\mathscr{R}$ be an additive subcategory of $\mathscr{C}$ satisfying \ref{R3}. Then $\mathscr{R}$ satisfies \ref{R2} if, and only if, $\mathscr{R}$ satisfies \ref{R'2}.
\end{prop}
\begin{proof}
It is clear that if $\mathscr{R}$ satisfies \ref{R2}, then $\mathscr{R}$ satisfies \ref{R'2}.
Assume that $\mathscr{R}$ satisfies \ref{R'2}. We will show that if we have a short exact sequence 
\[\begin{tikzcd}
	M & E & N
	\arrow["k",tail, from=1-1, to=1-2]
	\arrow["f",two heads, from=1-2, to=1-3]
\end{tikzcd}\]
where $M,N \in \mathscr{R}$ and $E \in  \mathscr{C}$, then $E \in \mathscr{R}$. We will proceed in two steps: first, we study the case where $N \in \ind(\mathscr{R})$, and then the general case.

First assume that $N \in \ind(\mathscr{R})$. Consider the decomposition into indecomposable summands: $M=\bigoplus^n_{i=0} M_i$. By pushing out along the $i$th projection $M \rightarrow M_i$ we obtain the family of short exact sequences: \[\begin{tikzcd}
	{M_i} & {E_i} & N
	\arrow[tail, from=1-1, to=1-2]
	\arrow[two heads, from=1-2, to=1-3]
\end{tikzcd}\]
As $M_i, N \in \ind(\mathscr{R})$, we have $E_i \in \mathscr{R}$ by \ref{R'2}. Consider the family of the following diagrams:
\[\begin{tikzcd}
	M & E & N \\
	{M_i} & {E_i} & N
	\arrow[tail, from=1-1, to=1-2]
	\arrow["{\pi_i}"', from=1-1, to=2-1]
	\arrow[two heads, from=1-2, to=1-3]
	\arrow[from=1-2, to=2-2]
	\arrow[Rightarrow, no head, from=1-3, to=2-3]
	\arrow[tail, from=2-1, to=2-2]
	\arrow[two heads, from=2-2, to=2-3]
\end{tikzcd}\]
Summing the bottom rows, we obtain:
\[\begin{tikzcd}
	M & {E} & N \\
	{\displaystyle \bigoplus_{i=0}^n M_i} & {\displaystyle \bigoplus_{i=0}^nE_i} & {\displaystyle \bigoplus_{i=0}^n N}
	\arrow[tail, from=1-1, to=1-2]
	\arrow[Rightarrow, no head, from=1-1, to=2-1]
	\arrow[two heads, from=1-2, to=1-3]
	\arrow[from=1-2, to=2-2]
	\arrow["\Delta", from=1-3, to=2-3]
	\arrow[tail, from=2-1, to=2-2]
	\arrow[two heads, from=2-2, to=2-3]
\end{tikzcd}\]
 By taking the mapping cone (see the short exact sequences as chain complexes), we obtain the following short exact sequence. \[\begin{tikzcd}
	E & {\displaystyle \bigoplus_{i=0}^n (E_i) \oplus N} & {\displaystyle \bigoplus_{i=0}^nN}
	\arrow[tail, from=1-1, to=1-2]
	\arrow[two heads, from=1-2, to=1-3]
\end{tikzcd}\] Thus $E$ is in $\mathscr{R}$ by \ref{R3}.

We now prove the general result by induction on the number of indecomposable summands of $N$, as before. The initialization is done above. Assume that $N$ is not indecomposable. We can write $N=N_1\oplus N'$ with $N_1 \in \ind(\mathscr{R})$ and $N' \in \mathscr{R}$ by additivity of $\mathscr{R}$. As in the proof of \cref{prop:red_ker}, we obtain the following diagram.
\[\begin{tikzcd}
	M & {E'} & {N'} \\
	M & E & N \\
	& {N_1} & {N_1}
	\arrow[tail, from=1-1, to=1-2]
	\arrow[Rightarrow, no head, from=1-1, to=2-1]
	\arrow[two heads, from=1-2, to=1-3]
	\arrow[tail, from=1-2, to=2-2]
	\arrow[tail, from=1-3, to=2-3]
	\arrow[tail, from=2-1, to=2-2]
	\arrow[two heads, from=2-2, to=2-3]
	\arrow[two heads, from=2-2, to=3-2]
	\arrow[two heads, from=2-3, to=3-3]
	\arrow[Rightarrow, no head, from=3-2, to=3-3]
\end{tikzcd}\]
As we have the following short exact sequence \[\begin{tikzcd}
	M & E'& N'
	\arrow[tail, from=1-1, to=1-2]
	\arrow[two heads, from=1-2, to=1-3]
\end{tikzcd},\] by induction hypothesis, we have that $E' \in \mathscr{R}$. Moreover, the following short exact sequence \[\begin{tikzcd}
	E' & E & N_1
	\arrow[tail, from=1-1, to=1-2]
	\arrow[two heads, from=1-2, to=1-3]
\end{tikzcd}\] allows us to conclude that $E \in \mathscr{R}$ by the first step of the proof. 
\end{proof}

\begin{proof}[Proof of \cref{thm:rescondonindec}]
The result follows from \cref{lem:othercharactresolv} \ref{ires}, \cref{prop:red_ker,prop:red_ext}.
\end{proof}

Thanks to our previous result, we update Takahashi's algorithm to focus only on non-projective indecomposable objects in $\mathscr{C}$.

\begin{algo} \label{algo:inductivRes} Let $\mathscr{C}$ be an abelian Krull--Schmidt category which has enough projective objects. Consider $\mathcal{X} \subseteq \pmb{\ind \setminus \proj}(\mathscr{C})$.
\begin{enumerate}[label=$(\arabic*)$,itemsep=1mm]
\item Set $\mathcal{X}^0 = \mathcal{X} \cup (\ind(\mathscr{C}) \cap \proj(\mathscr{C}))$.

\item \label{2algores1} At the $i$th iteration of the algorithm, we define $\mathcal{X}^{i+1}$ as the union of $\mathcal{X}^i$ and indecomposable summands of $Y \in \mathscr{C}$ appearing either:
\begin{enumerate}[label=$(2 \alph*)$,itemsep=1mm]
    \item in a short exact sequence 
    \[\begin{tikzcd}
	A & Y & B
	\arrow[tail, from=1-1, to=1-2]
	\arrow[two heads, from=1-2, to=1-3]
\end{tikzcd}\] with $A,B \in \mathcal{X}^i$, or,
    \item \label{2balgores1} in a short exact sequence \[\begin{tikzcd}
	Y & \bigoplus_{k=0}^p A_k & B
	\arrow[tail, from=1-1, to=1-2]
	\arrow[two heads, from=1-2, to=1-3]
\end{tikzcd}\] for some $p \in \mathbb{N}^*$ with $B,A_0,\ldots,A_p \in \mathcal{X}^i$;
\end{enumerate}
\item If $\mathcal{X}^{i+1} \neq \mathcal{X}^i$, then go back to Step $(2)$;

\item Otherwise, return $\mathcal{X}^{i+1}$.
\end{enumerate}
\end{algo}

\begin{theorem} \label{thm:inductiveRes}
Let $\mathscr{C}$ be an abelian Krull--Schmidt category with enough projective objects. For any set $\mathcal{X} \subseteq \pmb{\ind \setminus \proj}(\mathscr{C})$, we have \[\Res(\mathcal{X}) = \add \left( \bigcup_{i \geqslant 0} \mathcal{X}^i \right), \] and $\mathcal{X}^{\Res} = \left( \bigcup_{i \geqslant 0} \mathcal{X}^i \right) \setminus \proj(\mathscr{C})$
\end{theorem}

\begin{proof}
It follows from \cite[Remark 3.2]{T09} and  \cref{prop:red_ker,prop:red_ext}.
\end{proof}

\subsection{Discussion on the usability of the algorithm}
\label{ss:UsabAlgo}

First, note that \cref{algo:inductivRes} is taking  more steps than the one given by Takahashi. Nevertheless, it enables us to construct step by step $\mathcal{X}^{\Res}$ for all collections $\mathcal{X} \subseteq \pmb{\ind \setminus \proj}(\mathcal{C})$, easily selecting the objects to add. Recall that, by \cref{cor:resclos}, it suffices to describe the resolving closures of subsets of non-projective indecomposable objects to describe any resolving subcategory of $\mathscr{C}$.

Secondly, an important remark is that, in general, \cref{algo:inductivRes} does not terminate. For instance, if $\ind(\mathscr{C})$ is infinite. However, even if we assume that $\ind(\mathscr{C})$ is finite, another issue could arise from the step \ref{2balgores1} as we have to consider all the objects of $\add(\mathcal{X}^i)$. To overcome this problem, we can restrict ourselves to categories $\mathscr{C}$ which satisfy the following condition:
\begin{enumerate}[label = $(\mathsf{Mid})_{p_0}$,itemsep=1mm]
    \item \label{MidHyp}  There exists $p_0 \in \mathbb{N}^*$ such that, for any short exact sequence \[\begin{tikzcd}
	X & Y & Z
	\arrow[tail, from=1-1, to=1-2]
	\arrow[two heads, from=1-2, to=1-3]
\end{tikzcd}\] in $\mathscr{C}$ where $X,Z \in \ind(\mathscr{C})$, the object $Y$ admits at most $p_0$ indecomposable summands, up to isomorphism.
\end{enumerate}
Under those hypotheses, we can ensure that the algorithm terminates.

\begin{lemma} \label{lem:red_ker_mid} Let $\mathscr{C}$ be an abelian Krull--Schmidt category with enough projective objects. Assume that $\mathscr{C}$ satisfies \ref{MidHyp} for some $p_0 \in \mathbb{N}^*$. Let $\mathscr{R} \subseteq \mathscr{C}$ be an additive subcategory satisfying \ref{R2}. Then $\mathscr{R}$ satisfies \ref{R'3} if and only if $\mathscr{R}$ satisfies the following:
\begin{enumerate}[label=$(\mathsf{R''3})_{p_0}$, itemsep=1mm]
\item \label{R''3} for any short exact sequence \[\begin{tikzcd}
	K & \widetilde{M} & Y
	\arrow[tail, from=1-1, to=1-2]
	\arrow[two heads, from=1-2, to=1-3]
\end{tikzcd}\] with $Y \in \ind(\mathscr{R})$ and $\widetilde{M} \in \mathscr{R}$ admitting at most $p_0$ indecomposable summands, then $K \in \mathscr{R}$.
\end{enumerate}
\end{lemma}

\begin{proof}
It is clear that if $\mathscr{R}$ satisfies \ref{R'3} then $\mathscr{R}$ satisfies \ref{R''3}. 

Assume that $\mathscr{R}$ satisfies \ref{R''3}. Consider a short exact sequence 
\[\begin{tikzcd}
	K & M & Y
	\arrow[tail, from=1-1, to=1-2]
	\arrow[two heads, from=1-2, to=1-3]
\end{tikzcd}\] with $Y \in \ind(\mathscr{R})$ and $M \in \mathscr{R}$.
Consider the decomposition of $K$ into indecomposable summands: $K = \bigoplus_{i=0}^n K_i$. By pushing out along the $i$th projection, we get the family of short exact sequences. 
\[\begin{tikzcd}
	K_i & M_i & Y
	\arrow[tail, from=1-1, to=1-2]
	\arrow[two heads, from=1-2, to=1-3]
\end{tikzcd}\] By \ref{MidHyp}, for all $i \in \{0,\ldots,n\}$, we have that $M_i$ admits at most $p_0$ indecomposable summands. By \ref{R2} and additivity of $\mathscr{R}$, we have that $M_i \in \mathscr{R}$ for all $i \in \{0,\ldots, n\}$. By \ref{R''3}, we get that $K_i \in \mathscr{R}$ for all $i \in \{0,\ldots,n\}$ and, by additivity of $\mathscr{R}$, we have that $K \in \mathscr{R}$.  Moreover, from these short exact sequences, by proceeding the same way that we did in the proof of \cref{prop:red_ext}, namely summing the family of short exact sequences and taking some mapping cone, we get the following exact sequence.
\[\begin{tikzcd}
	M & {\displaystyle \left(\bigoplus_{i=0}^n M_i \right)\oplus Y} & {\displaystyle \bigoplus_{i=0}^n Y}
	\arrow[tail, from=1-1, to=1-2]
	\arrow[two heads, from=1-2, to=1-3]
\end{tikzcd}\] 
This sequence completes the proof.
\end{proof}

Let $\mathbb{K}$ be a field. For the following statement, we will restrict ourselves to $\mathbb{K}$-linear categories, meaning, in particular, that the homspaces between objects of $\mathscr{C}$ are $\mathbb{K}$-vector spaces. This restriction is natural when we look at module categories over a $\mathbb{K}$-algebra. We will also impose that $\ind(\mathscr{C})$ is finite, and that homspaces between indecomposable objects of $\mathscr{C}$ are at most one-dimensional. This is the kind of restriction we will address in our setting later (see \cref{sec:TreePart1}).

\begin{prop} \label{prop:algoend}
Let $\mathbb{K}$ be a field and $\mathscr{C}$ be a $\mathbb{K}$-linear Krull--Schmidt category with enough projective objects. Assume that $\ind(\mathscr{C})$ is finite, that $\mathscr{C}$ satisfies \ref{MidHyp} for some $p_0 \in \mathbb{N}^*$, and that forall $X,Y \in \ind(\mathscr{C}),\ \dim_\mathbb{K}(\Hom(X,Y)) \leqslant 1$. Then, for any $\mathcal{X} \subseteq \pmb{\ind \setminus \proj}(\mathscr{C})$, \cref{algo:inductivRes} terminates.
\end{prop}

\begin{proof} 
Since, for any object $U,V \in \mathscr{C}$, we can describe $\Hom(U,V)$ via homspaces from indecomposable summands of $U$ and the ones of $V$, we will restrict ourselves to the indecomposable case. By our assumptions, to describe a morphism of $\Hom(U,V)$, up to a change of basis, we only need to know if the induced morphism from each summand of $U$ to each one of $V$ is a zero map or not. Thus under our assumptions, there is a finite $\mathbb{K}$-basis of the homspace.

Since $\ind(\mathscr{C})$ is finite, thanks to \cref{lem:red_ker_mid} and by previous arguments, we know that the step \ref{2algores1} of \cref{algo:inductivRes} terminates in a finite number of calculations. Moreover, $(\mathcal{X}^i)_{i \geqslant 0}$ is an increasing sequence of subsets of $\ind(\mathscr{C})$. Thus, as $\ind(\mathscr{C})$ is finite, there must be an integer $j \geqslant 0$ such that $\mathcal{X}^j = \mathcal{X}^{j+1}$. This is equivalent to saying that the algorithm terminates.
\end{proof}

\begin{remark} \label{rem:Artincase} Let $\mathbb{K}$ be a field, and $A$ be an artinian $\mathbb{K}$-algebra. Butler and Ringel \cite{BR87} showed that $\mathscr{C} = 
\mod{A}$ satisfies \ref{MidHyp} for a certain $p_0$ which depends on the maximal lengths of both indecomposable projective and injective $A$-modules. 
\end{remark}

In general, it is difficult to give an explicit recipe to get the objects appearing in step \ref{2algores1} to construct $\mathcal{X}^{\Res}$ from any $\mathcal{X} \subseteq \pmb{\ind \setminus \proj}(\mathscr{C})$. 

In what follows, we give explicit computations in the setup of module categories of some gentle algebras, and we refine \cref{algo:inductivRes}. Consequently, it will allow us to describe all resolving subcategories of those categories explicitly.

	\section{Gentle algebras}
	\label{sec:Gentle}
	\pagestyle{plain}

The section starts with the relevant basic definitions on gentle quivers in \cref{ss:GentleQuivAlg}, before recalling the classification of their indecomposable representations in \cref{ss:stringreps}, and the complete description of morphisms and extensions between them in \cref{ss:morphindec,ss:extindec}, kernels of epimorphisms in \cref{ss:GentleKernelEpi} and projective resolutions in \cref{ss:ProjRes}.In particular, we recall that modules over such an algebra satisfy the extra hypotheses required by the improved algorithm. From now on, $\mathbb{K}$ is an algebraically closed field.

\subsection{Gentle quivers and gentle algebras}
\label{ss:GentleQuivAlg} 
A \new{quiver} is given by a quadruplet $Q=(Q_0,Q_1,s,t)$ where:
\begin{enumerate}[label = $\bullet$, itemsep=0.1em]
    \item $Q_0$ is the vertex set of $Q$;
    \item $Q_1$ is the arrow set of $Q$;
    \item $s,t; Q_1 \longrightarrow Q_0$ are the source and target functions, respectively.
\end{enumerate}
In what follows, we assume that $Q$ is \new{finite}, meaning that $Q_0$ and $Q_1$ are finite sets.

A \new{path} of $Q$ is either a formal element $e_q$ for some $q \in Q_0$ called the \emph{lazy path at $q$}, or a finite sequence of arrows $(\alpha_1, \ldots, \alpha_k)$, for some $k \geqslant 1$, such that $s(\alpha_{i+1}) = t(\alpha_i)$ for all $1 \leqslant i < k$. From now on, we write $c = \alpha_k \cdots \alpha_1$ for the non-lazy path $(\alpha_1, \ldots, \alpha_k)$, and we set $s(c) = s(\alpha_1)$ and $t(c) = t(\alpha_k)$. We also set $s(e_q) = t(e_q) = q$ for all $q \in Q_0$. The \new{length} of a path $c$, denoted by $\ell(c)$, is the number of arrows appearing in $c$. A \new{cycle} is a non-lazy path whose source and target coincide.

For $\alpha \in Q_1$, write $\alpha^{-1}$ for the formal inverse of $\alpha$. Set $s(\alpha^{-1}) = t(\alpha)$ and $t(\alpha^{-1}) = s(\alpha) $. A \new{walk} of $Q$ is either a lazy path or a finite sequence of arrows and inverse arrows $\rho = (\beta_1,\ldots, \beta_k)$, for some $k \geqslant 1$, such that $s(\beta_{i+1}) = t(\beta_i)$ for all $1 \leqslant i < k$. From now on, we write $\rho = \alpha_k^{\varepsilon_k} \ldots \alpha_1^{\varepsilon_1}$ for the non-lazy walk $(\alpha_1^{\varepsilon_1},\ldots,\alpha_k^{\varepsilon_k})$, where $\alpha_i \in Q_1$ and $\varepsilon_i \in \{\pm 1\}$. As previously for paths, we set $s(\rho) = s(\alpha_1^{\varepsilon_1})$ and $t(\rho) = t(\alpha_k^{\varepsilon_k})$. 

In what follows, we assume that $Q$ is \new{connected}; meaning that for any pair $(v_1,v_2) \in Q_0$, there exists a walk $\rho$ of $Q$ such that $s(\rho) = v_1$ and $t(\rho) = v_2$.

 The \new{path algebra} of $Q$ (over $\mathbb{K}$), denoted by $\mathbb{K}Q$, is the $\mathbb{K}$-vector space with a basis given by all paths of $Q$ endowed with a multiplication defined on paths as follows: for $c_1$ and $c_2$ two paths of $Q$, we set \[c_2 \cdot c_1 = \begin{cases}
c_2 c_1 & \text{if } s(c_2) = t(c_1), \\ 0 & \text{otherwise.}
\end{cases}\] A \new{relation} on $Q$ is a $\mathbb{K}$-linear combination involving a finite set of paths of length greater than or equal to $2$ that have the same source and the same target. Such a relation is said to be \new{quadratic monomial} if it involves only one path and that path is of length $2$. For any set of relations $R$, write $\langle R \rangle$ for the (bi-sided) ideal in $\mathbb{K}Q$ generated by the elements of $R$. 

For $\ell \in \mathbb{N}$, write $\mathbb{K}Q_{\geqslant \ell}$ for the ideal generated by all paths $c$ such that $\ell(c) \geqslant \ell$. An ideal $I \subseteq \mathbb{K}Q$ is said to be \new{admissible} if $\mathbb{K}Q_{\geqslant N} \subseteq I \subseteq \mathbb{K}Q_{\geqslant 2}$ for some $N \geqslant 2$. Recall that if $I$ is an admissible ideal, the quotient algebra $\mathbb{K}Q/I$ is finite-dimensional.
\begin{definition} \label{def:gentle}
	A \new{gentle quiver} is a pair $(Q,R)$ such that:
	\begin{itemize}
		\item $Q$ is a quiver such that there are at most two incoming arrows and at most two outgoing arrows at each vertex of $Q$;
		\item  $R$ is a set of quadratic monomial relations such that: 
		\begin{itemize}
			\item for any arrow $\alpha \in Q_1$:
			\begin{itemize}
				\item there is at most one $\beta \in Q_1$ such that $s(\beta) = t(\alpha)$ and $\beta \alpha \in R$;
				\item there is at most one $\gamma \in Q_1$ such that $s(\gamma) = t(\alpha)$ and $\gamma \alpha \notin R$;
				\item there is at most one $\beta' \in Q_1$ such that $t(\beta') = s(\alpha)$ and $\alpha \beta' \in R$;
				\item there is at most one $\gamma' \in Q_1$ such that $t(\gamma') = s(\alpha)$ and $\alpha \gamma' \notin R$;
			\end{itemize}
		\item for any cycle $c = \alpha_k \cdots \alpha_1$, either there exists $i \in \{1, \ldots, k-1\}$ such that $\alpha_{i+1} \alpha_i \in R$, or $\alpha_1 \alpha_k \in R$.
		\end{itemize}
	\end{itemize}
	A $\mathbb{K}$-algebra $\Lambda$ is said to be \new{gentle} if there exists a gentle quiver $(Q,R)$ such that $\Lambda \cong \mathbb{K}Q/\langle R \rangle$.
\end{definition}

Remark that if $(Q,R)$ is a gentle quiver, then $\langle R \rangle$ is an admissible ideal in $\mathbb{K}Q$, and thus any gentle algebra is finite-dimensional. 

\smallskip
A \new{representation} of $(Q,R)$ is a pair $E = ((E_q)_{q \in Q_0}, (E_\alpha)_{\alpha \in Q_1})$ such that:
\begin{enumerate}[label = $\bullet$, itemsep=1mm]
	\item for any $q \in Q_0$, $E_q$ is a $\mathbb{K}$-vector space;
	\item for any $\alpha \in Q_1$, $E_\alpha: E_{s(\alpha)} \longrightarrow E_{t(\alpha)}$ is a $\mathbb{K}$-linear transformation;
	\item for any $(\alpha, \beta) \in (Q_1)^2$ such that $\beta \alpha \in R$, then $E_\beta E_\alpha = 0$.
\end{enumerate}
A representation $E$ of $(Q,R)$ is \new{finite dimensional} whenever $E_q$ is a finite-dimensional $\mathbb{K}$-vector space for all $q \in Q_0$.


Denote by $\rep(Q,R)$ the Krull-Schmidt category of finite-dimensional representations of $(Q,R)$. 


Write $\proj(Q,R)$ for the subcategory of projective representations in $\rep(Q,R)$


and $\ind(Q,R)$ for the set of isomorphism classes of indecomposable representations of $(Q,R)$. We say that $(Q,R)$ is \new{representation-finite} whenever $\#\ind(Q,R) < \infty$. In what follows, and for our purpose, we focus on representation-finite gentle quivers.

\subsection{String representations} 
\label{ss:stringreps} 
In this subsection, we recall the string combinatorics that give a complete classification of indecomposable representations of any representation-finite gentle quiver $(Q,R)$. Let us first begin with the following definition.

\begin{definition}
	\label{def:string}
	A \new{string} of $(Q,R)$ is either a lazy path, or a non-lazy walk $\rho = \alpha_k^{\varepsilon_k} \cdots \alpha_1^{\varepsilon_1}$ of $Q$ such that, for all $1 \leqslant i < k$:
	\begin{enumerate}[label = $\bullet$,itemsep=1mm]
		\item $\alpha_{i+1} \neq \alpha_i$ whenever $\varepsilon_{i+1} = - \varepsilon_i$: in this case we say that $\rho$ is \new{reduced};
		\item $\alpha_{i+1} \alpha_i \notin R$ whenever $\varepsilon_{i+1} = \varepsilon_i = 1$; and,
		\item $\alpha_i \alpha_{i+1} \notin R$ whenever $\varepsilon_{i+1} = \varepsilon_i = -1$.
	\end{enumerate}
    We call the first arrow (respectively, last arrow) of a string the arrow $\alpha_1$ (respectively, $\alpha_k$).
\end{definition}

These strings allow us to construct a particular family of isomorphism classes of indecomposable representations of $(Q,R)$.

\begin{definition}
	\label{def:stringrep}
	Let $(Q,R)$ be a gentle quiver. Consider $\rho = \alpha_k^{\varepsilon_k} \cdots \alpha_1^{\varepsilon_1}$ to be a string of $(Q,R)$, the \new{standard string representation} $M(\rho)$ is a representation of $Q$ obtained from $\rho$ as follows:
	\begin{enumerate}[label = $\bullet$,itemsep=1mm]
		\item Let $v_0 = s(\alpha_1^{\varepsilon_1})$ and for all $ i \in \{ 1, \ldots, k\}$, $v_i = t(\alpha_i^{\varepsilon_i})$;
		
		\item For $q \in Q_0$, $M(\rho)_q$ is the vector space having as basis $\{ x_i \mid v_i = q \}$, where $x_i$ are formal elements;
		
		\item For $\beta \in Q_1$, $M(\rho)_\beta: M(\rho)_{s(\beta)} \longrightarrow M(\rho)_{t(\beta)}$ is the linear transformation such that:  \[ M(\rho)_\beta(x_i) =  \begin{cases}
			x_{i-1} & \text{if } \alpha_i = \beta \text{ and } \varepsilon_i = -1 \hfill\\
			x_{i+1} & \text{if } \alpha_{i+1} = \beta \text{ and } \varepsilon_{i+1} = 1 \\
			0 & \text{otherwise.}
		\end{cases} \]
	\end{enumerate}
	A representation $E$ of $(Q,R)$ is a \new{string representation} whenever $E$ is isomorphic to some standard string representation.
\end{definition}

The following well-known result shows that string representations form a complete set of representatives for $\ind(Q,R)$.

\begin{theorem}[\cite{BR87}] \label{thm:BR} Let $\mathbb{K}$ be an algebraically closed field. Let $(Q,R)$ be a representation-finite gentle quiver. Any $E \in \ind(Q,R)$ is isomorphic to a string representation. Moreover, given two strings of $(Q,R)$ $\rho$ and $\rho'$, $M(\rho)$ is isomorphic to $M(\rho')$ if and only if $\rho' = \rho^{\pm 1}$.
\end{theorem}

We identify $\rho$ with its inverse $\rho^{-1}$ as they give isomorphic representations.

\begin{remark} \label{rem:BRgeneral} 
We can state a generalization of the previous result for any gentle quiver by considering, in addition, \emph{band representations}. They are parameterized by primitive cyclic strings of $(Q,R)$, nonzero elements of $\mathbb{K}$, and positive integers. However, as we are restricting ourselves to representation-finite gentle quivers, we will not need to consider band representations in this paper.
\end{remark}

\subsection{Morphisms between indecomposable representations}
\label{ss:morphindec} Now we recall a combinatorial way to describe the homomorphism space between two fixed indecomposable representations of $(Q,R)$.

\begin{definition} \label{def:substringtopbot}
	Let $\rho = \alpha_k^{\varepsilon_k} \cdots \alpha_1^{\varepsilon_1}$ be a string of $(Q,R)$. A \new{substring} of $\rho$ is either a lazy path $e_q$ for some $q \in \Supp_0(\rho)$ or a string $\alpha_j^{\varepsilon_j} \cdots \alpha_i^{\varepsilon_i}$ for some $1 \leqslant i \leqslant j \leqslant k$ together with its position in the string. A substring $\sigma = \alpha_j^{\varepsilon_j} \ldots \alpha_i^{\varepsilon_i}$ of $\rho$ is said to be \new{on top} of $\rho$ if $i =1$ or $\varepsilon_{i-1} = -1$, and $j=k$ or $\varepsilon_{j+1} = 1$. In the same way $\sigma$ is said to be \new{at the bottom} of $\rho$ if $i=1$ or $\varepsilon_{i-1} = 1$, and $j=k$ or $\varepsilon_{j+1} = -1$. Note that any string is both on top of and at the bottom of itself.
\end{definition}
\begin{conv}
Note that in what follows, arrows are pictured pointing downwards, such that on top and at the bottom, substrings will be pictured above or below the string, respectively. The arrow $\alpha_i$ is also pointing left or right depending on whether $\varepsilon_i$ in $\alpha_i^{\varepsilon_i}$ is either -1 or +1.
\end{conv}
\begin{figure}[!ht]
\centering
	\scalebox{0.9}{\begin{tikzpicture}[yshift = 1cm, xshift = 5cm, ->,line width=0.2mm,>= angle 60,color=black, scale=0.8]
			\node (rho) at (-1,.5){$\rho =$};
			\node (1) at (0,0){$\bullet$};
			\node (2) at (2,0){$\bullet$};
			\node (3) at (3,1){$\bullet$};
			\node (4) at (6,1){$\bullet$};
			\node (5) at (7,0){$\bullet$};
			\node (6) at (9,0){$\bullet$};
			\draw[-,decorate, decoration={snake,amplitude=.4mm}] (1) -- (2);
			\draw (3) -- node[above left]{$\alpha_{i-1}$} (2);
			\draw[-,line width=0.7mm,decorate, decoration={snake,amplitude=.4mm}] (3) -- node[above]{$\pmb \sigma$} (4);
			\draw (4) -- node[above right]{$\alpha_{j+1}$} (5);
			\draw[-,decorate, decoration={snake,amplitude=.4mm}]  (5) -- (6);
			\node at (-.3,0.5){$\left(\vphantom{\begin{matrix}
						\\
						\\
						\\
				\end{matrix}}\right.$};
			\node at (2.8,0.5){$\left.\vphantom{\begin{matrix}
						\\
						\\
						\\
				\end{matrix}}\right)$};
			\node at (6.2,0.5){$\left(\vphantom{\begin{matrix}
						\\
						\\
						\\
				\end{matrix}}\right.$};
			\node at (9.3,0.5){$\left.\vphantom{\begin{matrix}
						\\
						\\
						\\
				\end{matrix}}\right)$};
			
			\begin{scope}[yshift=-2cm]
				\node (rho) at (-1,-.5){$\rho =$};
				\node (1) at (0,0){$\bullet$};
				\node (2) at (2,0){$\bullet$};
				\node (3) at (3,-1){$\bullet$};
				\node (4) at (6,-1){$\bullet$};
				\node (5) at (7,0){$\bullet$};
				\node (6) at (9,0){$\bullet$};
				\draw[-,decorate, decoration={snake,amplitude=.4mm}] (1) --  (2);
				\draw (2) -- node[below left]{$\alpha_{i-1}$} (3);
				\draw[-,line width=0.7mm,decorate, decoration={snake,amplitude=.4mm}] (3) -- node[below]{$\pmb \sigma$} (4);
				\draw (5) -- node[below right]{$\alpha_{j+1}$} (4);
				\draw[-,decorate, decoration={snake,amplitude=.4mm}]  (5) -- (6);
				\node at (-.3,-0.5){$\left(\vphantom{\begin{matrix}
							\\
							\\
							\\
					\end{matrix}}\right.$};
				\node at (2.8,-0.5){$\left.\vphantom{\begin{matrix}
							\\
							\\
							\\
					\end{matrix}}\right)$};
				\node at (6.2,-0.5){$\left(\vphantom{\begin{matrix}
							\\
							\\
							\\
					\end{matrix}}\right.$};
				\node at (9.3,-0.5){$\left.\vphantom{\begin{matrix}
							\\
							\\
							\\
					\end{matrix}}\right)$};
			\end{scope}
	\end{tikzpicture}}
\caption{\label{fig:aboveandbelow} Figures showing a substring $\sigma$ of $\rho$ on top of $\rho$ (above), and a substring $\sigma$ of $\rho$ at the bottom of $\rho$ (below).}
\end{figure}

These combinatorial notions on strings are sufficient to describe the homomorphism space between any pair of indecomposable representations of $(Q,R)$.

\begin{theorem}[\cite{CB89}] \label{thm:CB}
	Let $(Q,R)$ be a representation-finite gentle quiver. Consider $\rho$ and $\rho'$ two strings of $(Q,R)$. Then $\mathsf{Hom}(\MM(\rho),\MM(\rho')) \cong \mathbb{K}^{\#[\rho, \rho']}$ where $[\rho, \rho']$ is the set of pairs $(\sigma, \sigma')$ such that  $\sigma$ is a substring on top of $\rho$, $\sigma'$ is a substring at the bottom of $\rho'$ and $\sigma' = \sigma^{\pm 1}$. 
	
	More precisely, we can describe a basis $(\varphi_{(\sigma, \sigma')})_{(\sigma, \sigma') \in [\rho,\rho']}$ of $\mathsf{Hom}(\MM(\rho),\MM(\rho'))$ where $\varphi_{(\sigma, \sigma')}$ is defined as follows:
	\begin{enumerate}[label=$\bullet$, itemsep=1mm]
    \item Let $(x_0,\ldots,x_p)$ and $(y_0,\ldots, y_q)$ be collections of linearly independent variables defining $\MM(\rho)$ and $\MM(\rho')$ as in \cref{def:stringrep};
    \item Let $x_i,\ldots, x_{j+1}$, for some $0 \leqslant i \leqslant j+1 \leqslant p$, be the ordered variables corresponding to the arrows of $\sigma$ and  $y_k, \ldots, y_{\ell+1}$  for some $0 \leqslant k \leqslant \ell+1\leqslant q$ be the ordered variables corresponding to the arrows of $\sigma'$;
    \item If $\sigma = \sigma'$, we define $\varphi_{(\sigma, \sigma')}$ on the basis $(x_0,\ldots,x_p)$ by:
    \[\forall u \in \{0,\ldots,p\},\ \varphi_{(\sigma, \sigma')}(x_u) = \begin{cases} 
    y_{k+u-i} & \text{if } i \leqslant u \leqslant j+1,\\
    0 & \text{otherwise.}
    \end{cases}\]
    A similar definition can be given if $\sigma = (\sigma')^{-1}$.
	\end{enumerate}
\end{theorem}
\begin{figure}[!ht]
\centering
	\scalebox{0.9}{\begin{tikzpicture}[yshift = 1cm, xshift = 5cm, ->,line width=0.2mm,>= angle 60,color=black, scale=0.8]
			\node (rho) at (-1,.5){$\rho =$};
			\node (12) at (0,0){$\bullet$};
			\node (22) at (2,0){$\bullet$};
			\node (32) at (3,1){$\bullet$};
			\node (42) at (6,1){$\bullet$};
			\node (52) at (7,0){$\bullet$};
			\node (62) at (9,0){$\bullet$};
			\draw[-,decorate, decoration={snake,amplitude=.4mm}] (12) -- (22);
			\draw (32) -- node[above left]{$\alpha_{i-1}$} (22);
			\draw[-,line width=0.7mm,decorate, decoration={snake,amplitude=.4mm}] (32) -- node[above]{$\pmb \sigma$} (42);
			\draw (42) -- node[above right]{$\alpha_{j+1}$} (52);
			\draw[-,decorate, decoration={snake,amplitude=.4mm}]  (52) -- (62);
			\node at (-.3,0.5){$\left(\vphantom{\begin{matrix}
						\\
						\\
						\\
				\end{matrix}}\right.$};
			\node at (2.8,0.5){$\left.\vphantom{\begin{matrix}
						\\
						\\
						\\
				\end{matrix}}\right)$};
			\node at (6.2,0.5){$\left(\vphantom{\begin{matrix}
						\\
						\\
						\\
				\end{matrix}}\right.$};
			\node at (9.3,0.5){$\left.\vphantom{\begin{matrix}
						\\
						\\
						\\
				\end{matrix}}\right)$};
			
			\begin{scope}[yshift=-2cm]
				\node (rho) at (-1,-.5){$\rho' =$};
				\node (1) at (0,0){$\bullet$};
				\node (2) at (2,0){$\bullet$};
				\node (3) at (3,-1){$\bullet$};
				\node (4) at (6,-1){$\bullet$};
				\node (5) at (7,0){$\bullet$};
				\node (6) at (9,0){$\bullet$};
				\draw[-,decorate, decoration={snake,amplitude=.4mm}] (1) --  (2);
				\draw (2) -- node[below left]{$\beta_{k-1}$} (3);
				\draw[-,line width=0.7mm,decorate, decoration={snake,amplitude=.4mm}] (3) -- node[below]{$\pmb \sigma'$} (4);
				\draw (5) -- node[below right]{$\beta_{\ell+1}$} (4);
				\draw[-,decorate, decoration={snake,amplitude=.4mm}]  (5) -- (6);
				\node at (-.3,-0.5){$\left(\vphantom{\begin{matrix}
							\\
							\\
							\\
					\end{matrix}}\right.$};
				\node at (2.8,-0.5){$\left.\vphantom{\begin{matrix}
							\\
							\\
							\\
					\end{matrix}}\right)$};
				\node at (6.2,-0.5){$\left(\vphantom{\begin{matrix}
							\\
							\\
							\\
					\end{matrix}}\right.$};
				\node at (9.3,-0.5){$\left.\vphantom{\begin{matrix}
							\\
							\\
							\\
					\end{matrix}}\right)$};
			\end{scope}
			\draw[line width=0.5mm,blue] (32) -- node[left]{$\varphi_{(\sigma, \sigma')}$}(3);
			\draw[line width=0.5mm,blue] (4.5,0.7) -- node[above,rotate=-90]{$x_u \longmapsto y_{k+u-i}$} (4.5,-2.7);
			\draw[line width=0.5mm,blue] (42) -- (4);
	\end{tikzpicture}}
\caption{\label{fig:basiselementmorphisms} Figure illustrating combinatorially the map $\varphi_{(\sigma, \sigma')}$ defined in \cref{thm:CB}, for $\sigma = \sigma'$ a substring on the top of $\rho$ and at the bottom of $\rho'$.}
\end{figure}

\begin{remark} \label{rem:generalmorphGentle} One can describe the homomorphism space between indecomposable representations, including band representations, for any gentle quiver $(Q,R)$. We refer the reader to \cite{K88} for more details. As we focus on representation-finite gentle quivers, \cref{thm:CB} is enough for our purposes.
\end{remark}

To compute $\rep(Q,R)$, its indecomposable objects and its morphisms, we often use the \new{Auslander--Reiten quiver} of $(Q,R)$. This is a quiver, denoted by $\AR(Q,R)$, whose vertices are the (isomorphism classes of) indecomposable representations of $(Q,R)$, and whose arrows correspond bijectively to the basis elements of the $\mathbb{K}$-vector spaces of irreducible morphisms between indecomposable representations (see \cite{ASS06} for more details). In \cref{sec:TreePart1}, we will state useful results on $\AR(Q,R)$ in the case where $Q$ is a tree.

\subsection{Extensions between indecomposable representations}
\label{ss:extindec}

Thanks to the description of all the homomorphism spaces, extensions between indecomposable representations can be explained combinatorially.

\begin{definition} \label{def:arrowandoverlap}
Let $\rho$ and $\mu$ be two strings of $(Q,R)$.
\begin{enumerate}[label=$\bullet$,itemsep=1mm]
    \item An \new{arrow extension of $\rho$ by $\mu$} is a string $\nu$ in $(Q,R)$ such that there exists an arrow $a \in Q_1$ such that $\nu = \mu a \rho$. We denote by \new{$\ArExt(\rho,\mu)$} the set of arrow extensions of $\rho$ by $\mu$.
    
\begin{figure}[!ht]
\centering
	\scalebox{0.9}{\begin{tikzpicture}[yshift = 1cm, xshift = 5cm, ->,line width=0.2mm,>= angle 60,color=black, scale=0.8]
	        \node (rho) at (-1,-1){$\nu =$};
			\node (1) at (0,0){$\bullet$};
			\node (2) at (3,0){$\bullet$};
			\node (3) at (5,-2){$\bullet$};
			\node (4) at (8,-2){$\bullet$};
			\draw[-,decorate, decoration={snake,amplitude=.4mm}] (1) -- node[below]{$\rho$} (2);
			\draw[line width=0.7mm] (2) -- node[below left]{$\pmb a$} (3);
			\draw[-,decorate, decoration={snake,amplitude=.4mm}] (3) -- node[above]{$\mu$} (4);
	\end{tikzpicture}}
\caption{\label{fig:arrowext} Illustration of an arrow extension of $\rho$ by $\mu$.}
\end{figure}
    
    \item Assume that there exist five strings $\rho_\ell, \rho_r, \mu_\ell, \mu_r, \sigma$ such that $\rho = \rho_r \sigma \rho_\ell$ and $\mu = \mu_r \sigma \mu_\ell$ with the following conditions:
    \begin{enumerate}[label=$(\roman*)$, itemsep=1mm]
    \item If $\rho_\ell$ is lazy, then there exists $c \in Q_1$ and $\mu'_r$ a string of $(Q,R)$ such that $\mu_r = c^{-1} \mu'_r$
    \item If $\rho_r$ is lazy, then there exist $d \in Q_1$ and $\mu'_r$ a string of $(Q,R)$, such that $\mu_r = \mu'_r d$,
    \item If $\sigma$ is lazy, then:
    \begin{enumerate}[label*=$(\alph*)$,itemsep=1mm]
    \item if $\rho_r = \rho'_r b^{-1}$ and $\mu_r = \mu'_r d$, for some $b,d \in Q_1$ and some strings $\rho'_r$ and $\mu'_r$ of $(Q,R)$, then $db \in R$, and
    \item if $\rho_\ell = a \rho'_\ell$ and $\mu_\ell = c^{-1} \mu'_\ell $, for some $a,c \in Q_1$ and some strings $\rho'_r$ and $\mu'_\ell$ of $(Q,R)$ then $ca \in R$.
    \end{enumerate}
    \end{enumerate}
    An \new{overlap extension of $\rho$ by $\mu$} is a string $\nu$ of $(Q,R)$ such that either $\nu = \mu_r \sigma \rho_\ell$ or $\nu = \rho_r \sigma \mu_\ell$. We say that $\mu_r \sigma \rho_\ell$ and $\rho_r \sigma \mu_\ell$ are \emph{conjugated overlap extensions} of $\rho$ by $\mu$. We denote by \new{$\OvExt(\rho,\mu)$} the set of overlap extensions of $\rho$ by $\mu$. If $\rho$ and $\mu$ do not satisfy the above conditions, we set $\OvExt(\rho, \mu) = \varnothing$. 
    
\begin{figure}[!ht]
\centering
	\scalebox{0.9}{\begin{tikzpicture}[->,line width=0.2mm,>= angle 60,color=black, scale=0.8]
			\node (rho) at (-1,.5){$\mu =$};
			\node (1) at (0,0){$\bullet$};
			\node (2) at (2,0){$\bullet$};
			\node (3) at (3,1){$\bullet$};
			\node (4) at (6,1){$\bullet$};
			\node (5) at (7,0){$\bullet$};
			\node (6) at (9,0){$\bullet$};
			\draw[-,decorate, decoration={snake,amplitude=.4mm}] (1) -- node[below]{$\mu'_\ell$} (2);
			\draw (3) -- node[below right]{$c$} (2);
			\draw[-,line width=0.7mm,decorate, decoration={snake,amplitude=.4mm}] (3) -- (4);
			\draw (4) -- node[below left]{$d$} (5);
			\draw[-,decorate, decoration={snake,amplitude=.4mm}]  (5) -- node[below]{$\mu'_r$}  (6);
			\begin{scope}[yshift=3cm]
				\node (rho) at (-1,-.5){$\rho =$};
				\node (1) at (0,0){$\bullet$};
				\node (2) at (2,0){$\bullet$};
				\node (3) at (3,-1){$\bullet$};
				\node (4) at (6,-1){$\bullet$};
				\node (5) at (7,0){$\bullet$};
				\node (6) at (9,0){$\bullet$};
				\draw[-,decorate, decoration={snake,amplitude=.4mm}] (1) --  node[above]{$\rho'_\ell$} (2);
				\draw (2) -- node[above right]{$a$} (3);
				\draw[-,line width=0.7mm,decorate, decoration={snake,amplitude=.4mm}] (3) -- node[above]{$\pmb \sigma$} (4);
				\draw (5) -- node[above left]{$b$} (4);
				\draw[-,decorate, decoration={snake,amplitude=.4mm}]  (5) -- node[above]{$\rho'_r$} (6);
			\end{scope}
			\draw[-,decorate, decoration={snake,amplitude=.4mm}, dotted,orange, line width=0.5mm] (0,0.3) --  (1.8,0.3) -- (2.8,1.3) -- (6.2,1.7) -- (7.2,2.65) -- (9,2.65);
			\draw[-,decorate, decoration={snake,amplitude=.4mm},dashed,cyan, line width=0.5mm] (0,2.7) --  (1.8,2.7) -- (2.8,1.7) -- (6.2,1.3) -- (7.2,0.3) -- (9,0.3);
	\end{tikzpicture}}
\caption{\label{fig:overlapext} Illustration of the general configuration of an overlap extension of $\rho$ by $\mu$. The dashed and dotted colored lines correspond to the overlap extensions of $\rho$ by $\mu$.}
\end{figure}
\end{enumerate}
\end{definition}

\begin{prop}[\cite{CPS21}] \label{prop:arrowandoverlapinterpretext}
Let $\rho$ and $\mu$ be two strings in $(Q,R)$.
\begin{enumerate}[label=$\bullet$, itemsep=1mm]
    \item Assume that there is an arrow extension of $\rho$ by $\mu$.  Let $\nu = \mu a \rho$ in the notation of \cref{def:arrowandoverlap}. We have the following nonsplit short exact sequence. \[\begin{tikzcd}
	\MM(\mu) & \MM(\nu) & \MM(\rho)
	\arrow[tail, from=1-1, to=1-2]
	\arrow[two heads, from=1-2, to=1-3]
\end{tikzcd}\]

\begin{figure}[!ht]
\centering
	\[\begin{tikzcd}[>= angle 60]
	\vcenter{\hbox{\scalebox{0.4}{\begin{tikzpicture}[->,line width=0.2mm,>= angle 60,color=black, scale=0.8]
			\node (1) at (0,0){$\bullet$};
			\node (2) at (3,0){$\bullet$};
			\draw[-,decorate, decoration={snake,amplitude=.4mm}] (1) -- node[below]{$\mu$} (2);
	\end{tikzpicture}}}} & \vcenter{\hbox{\scalebox{0.4}{\begin{tikzpicture}[->,line width=0.2mm,>= angle 60,color=black, scale=0.8]
			\node (1) at (0,0){$\bullet$};
			\node (2) at (3,0){$\bullet$};
			\node (3) at (5,-2){$\bullet$};
			\node (4) at (8,-2){$\bullet$};
			\draw[-,decorate, decoration={snake,amplitude=.4mm}] (1) -- node[below]{$\rho$} (2);
			\draw[line width = 0.7mm] (2) -- node[below left]{$\pmb a$} (3);
			\draw[-,decorate, decoration={snake,amplitude=.4mm}] (3) -- node[above]{$\mu$} (4);
	\end{tikzpicture}}}} & \vcenter{\hbox{\scalebox{0.4}{\begin{tikzpicture}[->,line width=0.2mm,>= angle 60,color=black, scale=0.8]
			\node (1) at (0,0){$\bullet$};
			\node (2) at (3,0){$\bullet$};
			\draw[-,decorate, decoration={snake,amplitude=.4mm}] (1) -- node[below]{$\rho$} (2);
	\end{tikzpicture}}}}
	\arrow[line width=0.25mm, tail, from=1-1, to=1-2]
	\arrow[line width=0.25mm, two heads, from=1-2, to=1-3]
\end{tikzcd}\]
\caption{\label{fig:arrowextseq} The nonsplit exact sequence corresponding to the arrow extension of $\rho$ by $\mu$ pictured with their associated strings.}
\end{figure}

\item Assume that there is an overlap extension of $\rho$ by $\mu$. Let $\nu_1 = \mu_r \sigma \rho_\ell$ and $\nu_2 = \rho_r \sigma \mu_\ell $ in the notation of \cref{def:arrowandoverlap}. We have the following nonsplit short exact sequence. \[\begin{tikzcd}
	\MM(\mu) & \MM(\nu_1) \oplus \MM(\nu_2) & \MM(\rho)
	\arrow[tail, from=1-1, to=1-2]
	\arrow[two heads, from=1-2, to=1-3]
\end{tikzcd}\]

\begin{figure}[!ht]
\centering
	\[\begin{tikzcd}[>= angle 60]
	\vcenter{\hbox{\scalebox{0.4}{\begin{tikzpicture}[->,line width=0.2mm,>= angle 60,color=black, scale=0.8]
			\node (1) at (0,0){$\bullet$};
			\node (2) at (2,0){$\bullet$};
			\node (3) at (3,1){$\bullet$};
			\node (4) at (6,1){$\bullet$};
			\node (5) at (7,0){$\bullet$};
			\node (6) at (9,0){$\bullet$};
			\draw[-,decorate, decoration={snake,amplitude=.4mm}] (1) -- node[below]{$\mu'_\ell$} (2);
			\draw (3) -- node[above left]{$c$} (2);
			\draw[-,line width=0.7mm,decorate, decoration={snake,amplitude=.4mm}] (3) -- node[above]{$\pmb \sigma$} (4);
			\draw (4) -- node[above right]{$d$} (5);
			\draw[-,decorate, decoration={snake,amplitude=.4mm}]  (5) -- node[below]{$\mu'_r$}  (6);
	\end{tikzpicture}}}} & \vcenter{\hbox{\scalebox{0.4}{\begin{tikzpicture}[->,line width=0.2mm,>= angle 60,color=black, scale=0.8]
	        \node (0) at (4.5,3){$\bigoplus$};
			\node (1) at (0,0){$\bullet$};
			\node (2) at (2,0){$\bullet$};
			\node (3) at (3,1){$\bullet$};
			\node (4) at (6,1){$\bullet$};
			\node (5) at (7,2){$\bullet$};
			\node (6) at (9,2){$\bullet$};
			\draw[-,decorate, decoration={snake,amplitude=.4mm}] (1) -- node[below]{$\mu'_\ell$} (2);
			\draw (3) -- node[above left]{$c$} (2);
			\draw[-,line width=0.7mm,decorate, decoration={snake,amplitude=.4mm}] (3) -- node[above]{$\pmb \sigma$} (4);
			\draw (5) -- node[below right]{$b$} (4);
				\draw[-,decorate, decoration={snake,amplitude=.4mm}]  (5) -- node[above]{$\rho'_r$} (6);
			\begin{scope}[yshift=6cm]
				\node (1) at (0,0){$\bullet$};
				\node (2) at (2,0){$\bullet$};
				\node (3) at (3,-1){$\bullet$};
				\node (4) at (6,-1){$\bullet$};
				\node (5) at (7,-2){$\bullet$};
				\node (6) at (9,-2){$\bullet$};
				\draw[-,decorate, decoration={snake,amplitude=.4mm}] (1) --  node[above]{$\rho'_\ell$} (2);
				\draw (2) -- node[below left]{$a$} (3);
				\draw[-,line width=0.7mm,decorate, decoration={snake,amplitude=.4mm}] (3) -- node[below]{$\pmb \sigma$} (4);
				\draw (4) -- node[below left]{$d$} (5);
				\draw[-,decorate, decoration={snake,amplitude=.4mm}]  (5) -- node[below]{$\mu'_r$} (6);
			\end{scope}
	\end{tikzpicture}}}} & \vcenter{\hbox{\scalebox{0.4}{\begin{tikzpicture}[->,line width=0.2mm,>= angle 60,color=black, scale=0.8]
			\node (1) at (0,0){$\bullet$};
				\node (2) at (2,0){$\bullet$};
				\node (3) at (3,-1){$\bullet$};
				\node (4) at (6,-1){$\bullet$};
				\node (5) at (7,0){$\bullet$};
				\node (6) at (9,0){$\bullet$};
				\draw[-,decorate, decoration={snake,amplitude=.4mm}] (1) --  node[above]{$\rho'_\ell$} (2);
				\draw (2) -- node[below left]{$a$} (3);
				\draw[-,line width=0.7mm,decorate, decoration={snake,amplitude=.4mm}] (3) -- node[below]{$\pmb \sigma$} (4);
				\draw (5) -- node[below right]{$b$} (4);
				\draw[-,decorate, decoration={snake,amplitude=.4mm}]  (5) -- node[above]{$\rho'_r$} (6);
	\end{tikzpicture}}}}
	\arrow[line width=0.25mm, tail, from=1-1, to=1-2]
	\arrow[line width=0.25mm, two heads, from=1-2, to=1-3]
\end{tikzcd}\]
\caption{\label{fig:overlapextseq} The nonsplit exact sequence corresponding to the overlap extension of $\rho$ by $\mu$ pictured with their associated strings.}
\end{figure}
\end{enumerate}
\end{prop}

The theorem below justifies the terminology of arrow and overlap extensions, demonstrating that they are sufficient to describe any extension space.

\begin{theorem}[\cite{BDMTY19,CPS21}] \label{thm:ExplicitExtGentle}
Let $(Q,R)$ be a representation-finite gentle quiver. Let $\rho$ and $\mu$ be two strings of $(Q,R)$. A basis of $\Ext^1(\MM(\mu), \MM(\rho))$ is given by  nonsplit short exact sequences \[\begin{tikzcd}[>= angle 60]
	\MM(\mu) & E & \MM(\rho),
	\arrow[line width=0.25mm, tail, from=1-1, to=1-2]
	\arrow[line width=0.25mm, two heads, from=1-2, to=1-3]
\end{tikzcd}\] Such that, either:
\begin{enumerate}[label=$\bullet$,itemsep=1mm]
    \item $E \cong \MM(\nu)$ for some $\nu \in \ArExt(\rho, \mu)$; or,
    \item $E \cong \MM(\nu_1) \oplus \MM(\nu_2)$  with $\nu_1, \nu_2 \in\OvExt(\rho,\mu)$ that are conjugated.
\end{enumerate}
\end{theorem}

\begin{cor} \label{cor:GentleMid} The category $\rep(Q,R)$ satisfies \ref{MidHyp} for $p_0 = 2$. More precisely, for any short exact sequence 
\[\begin{tikzcd}[>= angle 60]
	X & E & Y
	\arrow[line width=0.25mm, tail, from=1-1, to=1-2]
	\arrow[line width=0.25mm, two heads, from=1-2, to=1-3]
\end{tikzcd}\] where $X,Y \in \ind(Q,R)$, the representation $E$ admits at most two indecomposable summands.
\end{cor}

\begin{remark} \label{rem:WellknownGentleMid} \cref{cor:GentleMid} is a well-known result for string algebras of which gentle algebras are a subclass, which follows from results due to \cite{BS80, WW85}, with methods given in \cite{GP68}. An independent proof using covering theory appeared in \cite{SW83,DS87}. 
\end{remark}

\subsection{Kernels of epimorphisms}
\label{ss:GentleKernelEpi}

In this subsection, we compute kernels of a certain family of minimal epimorphisms in $\rep(Q,R)$, which will be enough to study all possible epimorphisms in $\rep(Q,R)$, under additional assumptions on $(Q,R)$ (see \cref{sec:TreePart1}).

\begin{definition} \label{def:pcoverofstrings}
Let $(Q,R)$ be a gentle quiver. Let $\rho$ be a string of $(Q,R)$. Given some $p \in \mathbb{N}^*$, the $p$-tuple $(\sigma_1,\ldots,\sigma_p)$ of strings is said to be a \new{$p$-cover of $\rho$} whenever both of the following properties hold:
\begin{itemize}
    \item For each $i \in \{1,\ldots,p\}$, $\sigma_i$ is a substring of $\rho$; and,
    \item There exists a $p$-tuple of strings $(\widetilde{\sigma_1}, \ldots, \widetilde{\sigma_p})$ such that:
    \begin{itemize}
        \item for all $i \in \{1,\ldots,p\}$, $\widetilde{\sigma_i}$ is a substring of $\sigma_i$; and,
        \item $\rho = \widetilde{\sigma_p} \cdots \widetilde{\sigma_1}$.
    \end{itemize}
\end{itemize}
Such a $p$-cover of $\rho$ is reduced if there are no $j \in \{1,\ldots,p-1\}$ such that for some $(i_1,\ldots,i_j) \in \{1,\ldots,p\}$, $(\sigma_{i_1},\ldots, \sigma_{i_j})$ is a $j$-cover of $\rho$.
\end{definition}

\begin{ex} \label{ex:cover}
In \cref{fig:cover}, the strings $\sigma_1,\ldots,\sigma_7$ are substrings of the string $\mu$ in some ambient gentle quiver $(Q,R)$.
\begin{enumerate}[label=$\bullet$,itemsep=1mm]
    \item The $7$-tuple $(\sigma_1,\sigma_2,\sigma_3, \sigma_4, \sigma_5, \sigma_6, \sigma_7)$ is a $7$-cover of $\mu$ but it is not reduced.
    \item The $4$ tuple $(\sigma_1,\sigma_4,\sigma_6,\sigma_7)$ is a reduced $4$-cover of $\mu$.
\end{enumerate} 
\begin{figure}[!ht]
    \centering
    	\scalebox{0.9}{\begin{tikzpicture}[yshift = 1cm, xshift = 5cm, -,line width=0.2mm,>= angle 60,color=black, scale=1]
			\node[blue] at (0,0){$\bullet$};
			\node[blue] at (11,0){$\bullet$};
			\node[orange] at (0,0.3){$\bullet$};
			\node[orange] at (2,0.3){$\bullet$};
			\node[mypurple] at (5,0.3){$\bullet$};
			\node[orange] at (6,0.3){$\bullet$};
			\node[orange] at (10,0.3){$\bullet$};
			\node[mypurple] at (1,0.6){$\bullet$};
			\node[mypurple] at (3,0.6){$\bullet$};
			\node[mypurple] at (4,0.6){$\bullet$};
			\node[mypurple] at (8,0.6){$\bullet$};
			\node[orange] at (9,0.6){$\bullet$};
			\node[orange] at (11,0.6){$\bullet$};
			\node[orange] at (2,0.9){$\bullet$};
			\node[orange] at (7,0.9){$\bullet$};
			\draw[decorate, decoration={snake,amplitude=.4mm},blue] (0,0) to node[below]{$\mu$} (11,0);
			\draw[line width=0.5mm,decorate, decoration={snake,amplitude=.4mm},dotted,orange] (0,0.3) to node[above left]{$\sigma_1$}(2,0.3);
			\draw[line width=0.5mm,decorate, decoration={snake,amplitude=.4mm},dash pattern={on 10pt off 2pt on 5pt off 2pt},mypurple] (1,0.6) to node[above left]{$\sigma_2$}(3,0.6);
			\draw[line width=0.5mm,decorate, decoration={snake,amplitude=.4mm},dotted,orange] (2,0.9) to node[above]{$\sigma_4$}(7,0.9);
			\draw[line width=0.5mm,decorate, decoration={snake,amplitude=.4mm},dash pattern={on 10pt off 2pt on 5pt off 2pt},mypurple] (2,0.3) to node[above]{$\sigma_3$}(5,0.3);
			\draw[line width=0.5mm,decorate, decoration={snake,amplitude=.4mm},dash pattern={on 10pt off 2pt on 5pt off 2pt},mypurple] (4,0.6) to node[below left]{$\sigma_5$}(8,0.6);
			\draw[line width=0.5mm,decorate, decoration={snake,amplitude=.4mm},dotted,orange] (6,0.3) to node[above right]{$\sigma_6$}(10,0.3);
			\draw[line width=0.5mm,decorate, decoration={snake,amplitude=.4mm},dotted,orange] (9,0.6) to node[above]{$\sigma_7$}(11,0.6);
	\end{tikzpicture}}
    \caption{\label{fig:cover} An example of a $p$-cover of a string $\mu$}
\end{figure}
\end{ex}

\begin{lemma} \label{lem:covertripletdecomp}
Let $(Q,R)$ be a gentle quiver, and let $p \geqslant 2$ be an integer. Let $\mu$ be a string of $(Q,R)$. For any reduced $p$-cover $(\sigma_1,\ldots,\sigma_p)$, and for any $i \in \{1,\ldots,p-1\}$, there exists a unique triplet of strings $(\sigma_i^\ell, \sigma_{i,i+1}, \sigma_{i+1}^r)$ such that:
\begin{enumerate}[label=$\bullet$, itemsep=1mm]
    \item $\sigma_i = \sigma_{i,i+1} \sigma_i^\ell$; 
    \item $\sigma_{i+1} = \sigma_{i+1}^r \sigma_{i,i+1}$; and,
    \item $\sigma_{i,i+1}$ is the maximal common substring of $\sigma_{i}$ and $\sigma_{i+1}$ satisfying the equalities above.
\end{enumerate}
\end{lemma}
\begin{proof}
As $(\sigma_1,\ldots, \sigma_p)$ is a $p$-cover of $\mu$, there exists a substring $\widetilde{\sigma_i}$ of $\sigma_i$ and a substring $\widetilde{\sigma_{i+1}}$ of $\sigma_{i+1}$, such that $\widetilde{\sigma_{i+1}} \widetilde{\sigma_{i}}$ is a string. Moreover, as $(\sigma_1,\ldots, \sigma_p)$ is reduced, $\sigma_i$ cannot be a substring of $\sigma_{i+1}$ and vice-versa. Therefore, there exists at least one triplet of strings $(\sigma_i^{\ell}, \sigma_{i,i+1}, \sigma_{i+1}^r)$ such that $\sigma_i = \sigma_{i,i+1} \sigma_i^\ell$ and $\sigma_{i+1} = \sigma_{i+1}^r \sigma_{i,i+1}$. 

Note that if  $(\nu_i^\ell, \nu_{i,i+1}, \nu_{i+1}^r)$  is another triplet such that $\sigma_i = \nu_{i,i+1} \nu_i^\ell$ and $\sigma_{i+1} = \nu_{i+1}^r \nu_{i,i+1}$, and either $\nu_{i,i+1}$ is a substring of $\sigma_{i,i+1}$ or $\sigma_{i,i+1}$ is a substring of $\nu_{i,i+1}$. Therefore, the  triplet $(\sigma_i^{\ell}, \sigma_{i,i+1}, \sigma_{i+1}^r)$ satisfies the required properties.
\end{proof}

\begin{figure}[!ht]
    \centering
    	\scalebox{0.9}{\begin{tikzpicture}[yshift = 1cm, xshift = 5cm, -,line width=0.2mm,>= angle 60,color=black, scale=1]
			\node[blue] at (0,0){$\bullet$};
			\node[blue] at (11,0){$\bullet$};
			\node[orange] at (0,0.3){$\bullet$};
			\node[orange] at (2,0.3){$\bullet$};
			\node[orange] at (6,0.3){$\bullet$};
			\node[orange] at (10,0.3){$\bullet$};
			\node[orange] at (9,0.6){$\bullet$};
			\node[orange] at (11,0.6){$\bullet$};
			\node[orange] at (2,0.9){$\bullet$};
			\node[orange] at (7,0.9){$\bullet$};
			\draw[decorate, decoration={snake,amplitude=.4mm},blue] (0,0) to node[below]{$\mu$} (11,0);
			\draw[line width=0.5mm,decorate, decoration={snake,amplitude=.4mm},dotted,orange] (0,0.3) to node[above left]{$\sigma_1$}(2,0.3);
			\draw[line width=0.5mm,decorate, decoration={snake,amplitude=.4mm},dotted,orange] (2,0.9) to node[above]{$\sigma_2$}(7,0.9);
			\draw[line width=0.5mm,decorate, decoration={snake,amplitude=.4mm},dotted,orange] (6,0.3) to node[above right]{$\sigma_3$}(10,0.3);
			\draw[line width=0.5mm,decorate, decoration={snake,amplitude=.4mm},dotted,orange] (9,0.6) to node[above right]{$\sigma_4$}(11,0.6);
			\filldraw[gray,opacity=0.5] (1.95,1) to (1.95,-0.1) to node[below]{$\sigma_{1,2}$} (2.05,-.1) to (2.05,1) to cycle;
			\filldraw[gray,opacity=0.5] (6,1) to (6,-0.1) to node[below]{$\sigma_{2,3}$} (7,-.1) to (7,1) to cycle;
			\filldraw[gray,opacity=0.5] (9,1) to (9,-0.1) to node[below]{$\sigma_{3,4}$} (10,-.1) to (10,1) to cycle;
	\end{tikzpicture}}
    \caption{\label{fig:redcoveranddecomp} Illustration of \cref{lem:covertripletdecomp} on the reduced $4$-cover seen in \cref{ex:cover}.}
\end{figure}

\begin{prop} \label{prop:EpiandCover} Let $(Q,R)$ be a gentle quiver. For some $p \in \mathbb{N}^*$, consider $p+1$ strings, $\rho_1,\ldots, \rho_p$, and $\mu$, of $(Q,R)$ such that there exists a minimal epimorphism \[ \begin{tikzcd}
	{\displaystyle f: \bigoplus_{i=1}^p \MM(\rho_i)} & \MM(\mu)
	\arrow[two heads, from=1-1, to=1-2]
\end{tikzcd}\] such that $f_{|\MM(\rho_i)} = k_i \varphi_{(\sigma_i, \sigma_i')}$, where $\varphi_{(\sigma_i, \sigma_i')}\in\Hom(\MM(\rho_i),\MM(\mu))$, for some $k_i \in \mathbb{K}^\times$ and $(\sigma_i,\sigma_i') \in [\rho_i,\mu]$.
Then there exists a reduced $p$-cover $(\upsilon_1,\ldots,\upsilon_p)$ of $\mu$ and a permutation $\tau \in \mathfrak{S}_p$ such that, for all $i \in \{1,\ldots,p\}$, $\upsilon_i = \sigma_{\tau(i)}^{\pm 1}$.
\end{prop}

\begin{proof}
The result holds since $f$ is an epimorphism. Indeed, up to reordering the strings $\rho_1,\ldots, \rho_p$, and up to reversing the strings $\sigma_i$, we have that $(\sigma_1, \ldots, \sigma_p)$ is a $p$-cover of $\rho$. As $f$ is a minimal epimorphism, this $p$-cover is reduced.
\end{proof}

\begin{theorem}[\cite{CB89,CB18,BS21}] \label{thm:StringKerEpi} Let $(Q,R)$ be a representation-finite gentle quiver. For some $p \in \mathbb{N}^*$, consider $p+1$ strings, $\rho_1,\ldots, \rho_p$, and $\mu$, of $(Q,R)$ such that there exists a minimal epimorphism \[ \begin{tikzcd}
	{\displaystyle f: \bigoplus_{i=1}^p \MM(\rho_i)} & \MM(\mu)
	\arrow[two heads, from=1-1, to=1-2]
\end{tikzcd},\]
such that $f_{|\MM(\rho_i)} = k_i \varphi_{(\sigma_i,\sigma_i')}$, where $\varphi_{(\sigma_i, \sigma_i')}\in\Hom(\MM(\rho_i),\MM(\mu))$, for some $k_i \in \mathbb{K}^\times$ and $(\sigma_i,\sigma_i') \in [\rho_i,\mu]$. Then we have that $\Ker(f) = \bigoplus_{i=0}^p \MM(\nu_i)$, where $\nu_0,\ldots,\nu_p$ are strings of $(Q,R)$ constructed from $\rho_1,\ldots,\rho_p$, and $\mu$ as follows:
\begin{itemize}
    \item consider $(\upsilon_1,\ldots,\upsilon_p)$ to be the reduced $p$-cover of $\mu$ obtained  in \cref{prop:EpiandCover}, and order the strings $\rho_1,\ldots,\rho_p$ according to this $p$-cover;
    \item for $i \in \{1,\ldots, p\}$, decompose $\rho_i$ with respect to $\upsilon_i$: write $\rho_i = \rho_i^r \upsilon_i \rho_i^\ell$;
    \item  consider the unique triplet $(\upsilon_i^\ell, \upsilon_{i,i+1}, \upsilon_{i+1}^r)$ satisfying the conditions given in \cref{lem:covertripletdecomp}.
    \item Then we define $\nu_0,\ldots, \nu_p$ as follows:
    \begin{itemize}
    \item $\nu_0 =  \rho_{1}^{\ell}{}'$, where $\rho_{1}^{\ell}{}'$ is the string obtained from $\rho_{1}^{\ell}$ by deleting its last arrow (if $\rho_{1}^{\ell}$ is lazy, then there is no string $\nu_0$);
    \item $\nu_p = \rho_{p}^{r}{}'$, where $\rho_{p}^{r}{}'$ is the string obtained from $\rho_{p}^{r}$ by deleting its first arrow (if $\rho_{p}^{r}$ is lazy, then there is no string $\nu_p$); and,
    \item for $i \in \{1,\ldots,p-1\}$, $\nu_i = \rho_i^r \upsilon_{i,i+1} \rho_{i+1}^{\ell}$.
    \end{itemize}
\end{itemize}

\end{theorem}

\begin{figure}[!ht]
    \centering
    	\scalebox{0.9}{\begin{tikzpicture}[yshift = 1cm, xshift = 5cm, -,line width=0.2mm,>= angle 60,color=black, scale=1]
			\node[blue] at (0,0){$\bullet$};
			\node[blue] at (0,-0.4){$s(\mu)$};
			\draw[decorate, decoration={snake,amplitude=.4mm},blue] (0,0) to node[above]{$\mu$} (2,0);
			\draw[line width=0.5mm,decorate, decoration={snake,amplitude=.4mm},dash pattern={on 10pt off 2pt on 5pt off 2pt},mypurple] (-1,1.5) to [bend right=30] node[right]{$\rho_1$} (0,.1) to [bend left=0] (1,-.1) to [bend left=30]  (2,-1.5);
			\draw[line width=0.7mm,decorate, decoration={snake,amplitude=.4mm},dotted,orange] (-1.2,1.5) to [bend right=30] node[left]{$\nu_0$} (0,0);
			
			\filldraw[gray,opacity=0.4] (5.5,.3) to (5.5,-.4) to (7,-.4) to (7,.3) to cycle;
			\node[gray] at (6.25,-.6){$\upsilon_{i,i+1}$};
			\draw[decorate, decoration={snake,amplitude=.4mm},blue] (3,0) to (7,0) to node[above right]{$\mu$}  (9.5,0);
			\draw[line width=0.5mm,decorate, decoration={snake,amplitude=.4mm},dash pattern={on 10pt off 2pt on 5pt off 2pt},mypurple] (3,1.5) to [bend right=30] node[left]{$\rho_i$} (4,.2) to [bend left=0] (7,-.2) to [bend left=30]  (8,-1.5);
			\draw[line width=0.5mm,decorate, decoration={snake,amplitude=.4mm},dash pattern={on 10pt off 2pt on 5pt off 2pt},mypurple] (4.5,1.5) to [bend right=30] node[right]{$\rho_{i+1}$}  (5.5,.2) to [bend left=0] (8.5,-.2) to [bend left=30]  (9.5,-1.5);
			\draw[line width=0.5mm,decorate, decoration={snake,amplitude=.4mm},dashed,orange] (4.3,1.5) to [bend right=30] node[left]{$\nu_i$} (5.5,0.1) to [bend left=0] (7,-.1) to [bend left=30]  (8.2,-1.5);
	\end{tikzpicture}}
    \caption{\label{fig:KerEpiConstruction} Illustration of the construction of the string $\nu_0$ (left), and the string $\nu_i$ for $i \in \{1,\ldots,p-1\}$ (right), pictured as orange dotted lines, for given strings $\rho_1,\ldots,\rho_p$ (purple dashed lines), and $\mu$ (blue line). The gray part corresponds to the common substring $\upsilon_{i,i+1}$ of $\mu, \rho_{i}, \rho_{i+1}$, $\zeta_i$ and $\nu_i$.}
\end{figure}

\subsection{Projective covers and projective resolutions}
\label{ss:ProjRes}
In this subsection, we recall some useful notions and results about projective objects of $\rep(Q,R)$. For more details, we refer the reader to \cite{ASS06}. 

Let $(Q,R)$ be a representation-finite gentle quiver, and $M \in \rep(Q,R)$. As $\rep(Q,R)$ is Krull--Schmidt, consider the minimal projective resolution of any object in $\rep(Q,R)$. We denote by $\gldim(Q,R)$ the global dimension of $\rep(Q,R)$.

Note that to determine $\gldim(Q,R)$, it suffices to compute $\pdim(Y)$ for $Y \in \ind(Q,R)$.

\begin{prop}[\cite{ASS06}] \label{prop:Gentlegldimfinite} Let $(Q,R)$ be a representation-finite gentle quiver. Then $\gldim(Q,R) = \infty$ if and only if there exists a finite sequence of arrows $\alpha_1,\ldots, \alpha_p$, for some $p \in \mathbb{N}^*$ such that:
\begin{enumerate}[label=$\bullet$, itemsep=1mm]
    \item $s(\alpha_1) = t(\alpha_p)$, and, for all $i \in \{1,\ldots,p-1\}$, $s(\alpha_{i+1}) = t(\alpha_i)$; and,
    \item $\alpha_1 \alpha_p \in R$, and, for all $i \in \{1,\ldots, p-1\}$, $\alpha_{i+1} \alpha_i \in R$.
\end{enumerate}
\end{prop}
\begin{remark}
    The previous result will be relevant after a translation within the geometric model to characterize dissected surfaces associated with a gentle algebra of finite global dimension.
\end{remark}

    \section{Resolving posets for gentle trees}
	\label{sec:TreePart1}
	\pagestyle{plain}

From now on we restrict ourselves to certain quivers called gentle trees, and continue our study by examining properties of gentle trees in \cref{ss:Trees}. \cref{ss:VocabPoset} provides standard tools and notation for posets. Gentle trees enable us to create an order on non-projective indecomposable modules in \cref{ss:ResOrder} that conveys the data of monogeneous resolving subcategories. In \cref{ss:ResposetAndSubcat} we study the interaction between our poset of non-projective indecomposable modules and the poset of resolving subcategories ordered by inclusion.

\begin{conv}
In this section, all epimorphisms are assumed to be \new{minimal}. We note that an additively closed subcategory of $\rep(Q,R)$ is closed under kernels of epimorphisms if and only if it is closed under kernels of minimal epimorphisms.
\end{conv}
\subsection{Generalities on gentle trees} 
\label{ss:Trees}

\begin{definition} \label{def:gentle trees}
A \new{gentle tree} is a gentle quiver $(Q,R)$ such that the underlying graph is a tree.
\end{definition}

\begin{ex} See in \cref{fig:treeex} for an example of a gentle tree. \qedhere
\begin{figure}[!ht]
    \centering
     \scalebox{0.6}{\begin{tikzpicture}[->]
		\node (a) at (0,0) {$1$};
		\node (b) at (1,0) {$2$};
		\node (c) at (2,0) {$3$};
		\node (d) at (3,1) {$4$};
		\node (e) at (4,1) {$5$};
		\node (f) at (5,1) {$6$};
		\node (g) at (6,0) {$7$};
		\node (h) at (5,-1) {$8$};
		\node (i) at (7,0) {$9$};
		\node (j) at (8,0) {$10$};
		\node (k) at (9,0) {$11$};
		\node (l) at (10,0) {$12$};
		\node (m) at (11,-1) {$13$};
		\node (n) at (12,-1) {$14$};
		\node (o) at (13,0){$15$};
		\node (p) at (14,-1){$16$};
		\node (q) at (14,1){$17$};
		\node (r) at (15,1){$18$};
		
		\draw (a) -- (b);
		\draw (b) -- (c);
		\draw ([yshift=-1mm]d.west)--([yshift=1mm]c.east);
		\draw (f) -- (e);
		\draw (e) -- (d);
		\draw ([yshift=1mm]g.west)--([yshift=-1mm]f.east);
		\draw ([yshift=1mm]h.east) -- ([yshift=-1mm]g.west);
		\draw (i) to (g);
		\draw (i) -- (j);
		\draw (j) -- (k);
		\draw (l) to (k);
		\draw[<-] ([yshift=1mm]m.west) -- ([yshift=-1mm]l.east);
		\draw[<-] (n) to (m);
		\draw[<-]([yshift=-1mm]o.west) -- ([yshift=1mm]n.east);
		\draw ([yshift=-1mm]q.west)--([yshift=1mm]o.east);
		\draw[<-] ([yshift=1mm]p.west)--([yshift=-1mm]o.east);
		\draw (q)--(r);

		\draw[dashed,-] ([xshift=-.3cm,yshift=-.15cm]d.south) arc[start angle = -110, end angle = -15, x radius=.6cm, y radius =.6cm];
		\draw[dashed,-] ([xshift=-.3cm,yshift=.1cm]e.south) arc[start angle = -165, end angle = -15, x radius=.3cm, y radius =.3cm];
		\draw[dashed,-] ([xshift=-.5cm,yshift=0.2cm]f.south) arc[start angle = -165, end angle = -60, x radius=.6cm, y radius =.6cm];
		\draw[dashed,-] ([xshift=-.2cm,yshift=-.3cm]g.west) arc[start angle = 225, end angle = 135, x radius=.5cm, y radius =.5cm];
		\draw[dashed,-] ([xshift=.3cm,yshift=.1cm]n.north) arc[start angle = 45, end angle = 190, x radius=.4cm, y radius =.4cm];
		\draw[dashed,-] ([xshift=.3cm,yshift=-.15cm]m.north) arc[start angle = 10, end angle = 110, x radius=.5cm, y radius =.5cm];
		\draw[dashed,-] ([xshift=.2cm,yshift=-.3cm]o.east) arc[start angle = -45, end angle = 45, x radius=.5cm, y radius =.5cm];
		\end {tikzpicture}}
        \caption{ \label{fig:treeex} An example of a gentle tree.}
\end{figure}
\end{ex}

 The following theorem gives a complete description of the marked surface associated to a gentle tree.

The following two propositions can be seen as corollaries of the geometric model that will be introduced in \cref{sec:Surface}, even though they can be shown independently.

\begin{prop}\label{prop:vdimtree} Let $(Q,R)$ be a gentle tree. For any $M \in \ind (Q,R)$, and for any $q \in Q_0$, we have that $\dim(M_q) \leqslant 1$.
\end{prop}

\begin{prop}\label{prop:homtree} Let $(Q,R)$ be a gentle tree. For any $M,N \in \ind(Q,R)$, then $\dim(\Hom(M,N)) \leqslant 1$.
\end{prop}

The following notion will be useful to introduce a poset structure on the set of non-projective indecomposable modules in the next section.

\begin{definition}
Let $(Q,R)$ be a gentle algebra. A representation $X \in \rep(Q,R)$ is said to be \new{directing} whenever, for any $t \geqslant 1$, and for any sequence 
  \[\begin{tikzcd}
	X & M_1 & \cdots & M_t,
	\arrow["f_1", from=1-1, to=1-2]
	\arrow["f_2", from=1-2, to=1-3]
	\arrow["f_t", from=1-3, to=1-4]
\end{tikzcd},\]
where $f_1,\ldots,f_t$ are nonzero nonisomorphisms between $M_1,\ldots, M_t \in \ind_\mathbb{K}(Q,R)$, we have that $M_t$ and $X$ are not isomorphic.
\end{definition}

\begin{theorem}[{\cite[Chapter IX]{ASS06}}] \label{thm:acyclicAR} Let $(Q,R)$ be a gentle tree. Any $X \in \rep(Q,R)$ is directing, and $\AR(Q,R)$ is acyclic.
\end{theorem}

We also give a useful corollary of \cref{prop:homtree} which allows us to focus only on the basis elements of extensions between indecomposable representations.

\begin{cor}\label{cor:Exttree} Let $(Q,R)$ be a gentle tree. For any pair $(\rho, \mu)$ of strings of $(Q,R)$, we have that $\dim(\Ext^1(M(\rho), M(\mu))) \leqslant 1$.
\end{cor}

\begin{cor} \label{cor:numberArOvExt}
Let $(Q,R)$ be a gentle tree.
\begin{enumerate}[label=$(\roman*)$,itemsep=1mm]
    \item All nonsplit exact sequences are isomorphic to arrow and overlap extensions, as explained in \cref{prop:arrowandoverlapinterpretext};
    \item For any pair $(\rho, \mu)$ of strings of $(Q,R)$, $\#\ArExt(\rho, \mu) \in \{0,1\}$ and \\ 
    $\#\OvExt(\rho, \mu) \in \{0,2\}$; and,
    \item For any pair $(\rho, \mu)$ of strings of $(Q,R)$, the sets $\ArExt(\rho, \mu)$ and\\
    $\OvExt(\rho, \mu)$, cannot be both nonempty.
\end{enumerate}
\end{cor}

\subsection{Vocabulary on posets}
\label{ss:VocabPoset}
In this subsection, we recall some relevant notions on finite posets. We refer the reader to \cite{DP02,S11,CLM12} for further details.

Let $(\mathfrak{A}, \leqslant)$ be a finite poset. A \new{cover relation} is a pair $(x,y) \in \mathfrak{A}^2$ such that $x \leqslant y$, $x \neq y$, and, for all $z \in \mathfrak{A}$, whenever $x \leqslant z \leqslant y$, we have $x=z$ or $y= z$. The \new{Hasse diagram} of $(\mathfrak{A}, \leqslant)$ is a directed graph whose vertices are elements in $\mathfrak{A}$, and arrows $y \longrightarrow x$ are given by cover relations $x \leqslant y$.

Given $x,y \in \mathfrak{A}$, we say that $x$ and $y$ admit a \new{least upper bound} if there exists $z \in \mathfrak{A}$ such that \[ \forall t \in \mathfrak{A},\ t \geqslant z \Longleftrightarrow t \geqslant x \text{ and } t \geqslant y. \] In such a case, it is called the \new{join} of $x$ and $y$, and we denote it by $x \vee y$. We define the \new{greatest lower bound} dually, called the \new{meet} of $x$ and $y$, and denoted by $x \wedge y$, whenever it exists. We say that $(\mathfrak{A}, \leqslant)$ is a \new{lattice} whenever every pair $(x,y) \in \mathfrak{A}^2$ admits a join and a meet. Note that, as $\mathfrak{A}$ is finite, if $(\mathfrak{A},\leqslant)$ is a lattice, then it is  \new{complete}: $(\mathfrak{A}, \leqslant)$ admits a unique minimal element, denoted by $\widehat{0}_\mathfrak{A}$, and a unique maximal element, denoted by $\widehat{1}_\mathfrak{A}$. We recall the following useful characterization of a lattice.
\begin{prop}\label{prop:latticechar}
Let $(\mathfrak{A},\leqslant)$ be a finite poset. If:
\begin{enumerate}
    \item $(\mathfrak{A},\leqslant)$ admits a unique maximal element $\widehat{1}_\mathfrak{A}$, and
    \item every pair $(x,y) \in \mathfrak{A}$ admits a greatest lower bound,
\end{enumerate}
then $(\mathfrak{A},\leqslant)$ is a (complete) lattice. 
\end{prop}
\begin{remark}
The dual result also holds.
\end{remark}

A lattice $(\mathfrak{A},\leqslant)$ is \new{distributive} whenever the meet and join operations are distributive; meaning that, for all $x,y,z \in \mathfrak{A}$:
\begin{enumerate}[label=$\bullet$,itemsep=1mm]
    \item $x \wedge (y \vee z) = (x \wedge y) \vee (x \wedge z)$, and
    \item $(x\wedge y) \vee z = (x \vee z) \wedge (y \vee z)$.
\end{enumerate}

An \new{order ideal} of $(\mathfrak{A}, \leqslant)$ is a set $I$ such that, for all $x \in I$, $\langle x \rangle := \{y \in \mathfrak{A} \mid y \leqslant x\} \subseteq I$. Note that $\varnothing$ is an order ideal following our definition. We denote by $\mathscr{J}(\mathfrak{A}, \leqslant)$ the set of order ideals of $(\mathfrak{A}, \leqslant)$. We recall that $(\mathscr{J}(\mathfrak{A},\leqslant),\subseteq)$ is a complete distributive lattice, where:
\begin{enumerate}[label=$\bullet$, itemsep=1mm]
\item the join operation is the union,
\item the meet operation is the intersection, 
\item $\widehat{0}_{\mathscr{J}(\mathfrak{A}, \leqslant)} = \varnothing$, and $\widehat{1}_{\mathscr{J}(\mathfrak{A}, \leqslant)} = \mathcal{P}$.
\end{enumerate}

We say that $I \in \mathscr{J}(\mathfrak{A}, \leqslant)$ is \new{generated by a subset $\mathfrak{B}$} if $I = \bigcup_{x \in \mathfrak{B}} \langle x \rangle$. Given a subset $\mathfrak{B}$ of $\mathfrak{A}$,  we set $\langle \mathfrak{B} \rangle =  \bigcup_{x \in \mathfrak{B}} \langle x \rangle $ which is the order ideal generated by $\mathfrak{B}$.

Two elements $x,y \in \mathfrak{A}$ are said to be \new{comparable} whenever $x \leqslant y$ or $y \leqslant x$. Otherwise, we say that they are \new{non-comparable}. An \new{antichain} in $(\mathfrak{A}, \leqslant)$ is a collection $\mathcal{X}$ of elements in $\mathfrak{A}$ that are pairwise non-comparable.  We recall that any order ideal is generated by a unique antichain $\mathcal{X}_I$: if $I = \varnothing$, then $\mathcal{X}_I = \varnothing$; otherwise $\mathcal{X}_I$ is the set of maximal elements in $I$.

A non-minimal element $x \in \mathfrak{A}$ is said to be \new{join-irreducible} if, for all $(y,z) \in \mathfrak{A}^2$, whenever we have $x = y \vee z$, we get $x = y$ or $x = z$. 

\begin{remark} \label{rem:completlyjoinirr}
   If $(\mathfrak{A}, \leqslant)$ is a finite lattice, a join-irreducible element $x \in \mathfrak{A}$ is \emph{completely join-irreducible}; meaning that, for all collection $(x_i)_{i \in I}$ of elements in $\mathfrak{A}$, if $x = \bigvee_{i \in I} x_i$ then there exists $i \in I$ such that $x_i = x$. 
\end{remark}

\begin{prop} \label{prop:JoinIrreducible} Let $(\mathfrak{A}, \leqslant)$ be a finite lattice, and $x \in \mathfrak{A}$. Then $x$ is join-irreducible if and only if there exists exactly one element in $\mathfrak{A}$ covered by $x$.
\end{prop}

\subsection{Resolving order on non-projective indecomposable representations} \label{ss:ResOrder}

In this section, given a gentle tree $(Q,R)$, we will construct a poset structure on the non-projective indecomposable representations of $(Q,R)$ based on the monogeneous resolving subcategories they generate.

\begin{definition}\label{def:orders} Let $(Q,R)$ be a gentle quiver. We define the \new{$\Hom$-relation}, denoted by \new{$\Homleq$}, on $\ind(Q,R)$ as follows: \[ \forall X, Y \in \ind(Q,R),\  X \Homleq Y \Longleftrightarrow \Hom(X,Y) \neq 0. \] We write $\THomleq$ for the transitive closure of $\Homleq$. We also define the \new{$\Res$-relation}, denoted by \new{$\Resleq$}, on $\pmb{\ind \setminus \proj}(Q,R)$ as follows:  \[\forall X,Y \in\pmb{ \ind \setminus \proj}(Q,R),\ X \Resleq Y \Longleftrightarrow \Res(X) \subseteq \Res(Y). \]
\end{definition}

We can easily check that the $\Hom$-relation is reflexive on $\ind(Q,R)$, and the $\Res$-relation is a pre-order on $\pmb{\ind \setminus \proj}(Q,R)$: it is reflexive and transitive. However, in general, neither of them is an order relation.

\begin{ex}
	\label{ex:ResKro} Let $(Q,R)$ be the \emph{Kronecker quiver}, which is given by the quiver with relation :
	\[\begin{tikzpicture}[>= angle 60,<-]
	    \node (QR) at (-1,0){$(Q,R) =$};
		\node (a) at (0,0) {$1$};
		\node (b) at (2,0) {$2.$};
		\draw ([yshift=1mm]b.west)--node[above]{$\alpha$}([yshift=1mm]a.east);
		\draw  ([yshift=-1mm]b.west)--node[below]{$\beta$}([yshift=-1mm]a.east);
		\end {tikzpicture} \]
	Fix $\lambda \in \mathbb{K}^\times$, and consider the following non-isomorphic non-projective indecomposable representations.
	\[\begin{tikzpicture}[>= angle 60,<-]
	    \node (QR) at (-.7,0){$X =$};
		\node (a) at (0,0) {$\mathbb{K}$};
		\node (b) at (2,0) {$\mathbb{K}$};
		\node (c) at (3,0) {$\text{and}$};
		\draw ([yshift=1mm]b.west)--node[above]{$1$}([yshift=1mm]a.east);
		\draw  ([yshift=-1mm]b.west)--node[below]{$\lambda$}([yshift=-1mm]a.east);
		\begin{scope}[xshift=5cm]
		\node (QR) at (-.7,0){$Y =$};
		\node (a) at (0,0) {$\mathbb{K}^2$};
		\node (b) at (2,0) {$\mathbb{K}^2$};
		\draw ([yshift=1mm]b.west)--node[above]{$\left[ \begin{matrix}
		1 & 0 \\
		0 & 1 \end{matrix}\right]$}([yshift=1mm]a.east);
		\draw  ([yshift=-1mm]b.west)--node[below]{$\left[ \begin{matrix}
		\lambda & 1 \\
		0 & \lambda \end{matrix}\right]$}([yshift=-1mm]a.east);
		\end{scope}
		\end {tikzpicture} \]
		Then $X \Resleq Y$ since $X$ appears as a summand of a representation in $\Ext^1(Y,Y)$. For an analogous reason, we also have $Y \Resleq X$.
\end{ex}


\begin{ex}
	\label{ex:respdiminfinite} Let $(Q,R)$ be the following gentle quiver.
	\[\begin{tikzpicture}[>= angle 60,<-]
	    \node (QR) at (-1,0){$(Q,R) =$};
		\node (a) at (0,0) {$1$};
		\node (b) at (2,0) {$2$};
		\draw[bend right=20] ([yshift=1mm]b.west) to node[above]{$\alpha$} ([yshift=1mm]a.east);
		\draw[bend right=20]  ([yshift=-1mm]a.east) to node[below]{$\beta$}([yshift=-1mm]b.west);
		\draw[dashed,-] ([yshift=-.15cm, xshift=-.2cm]b.west) arc[start angle = -150, end angle = -210, x radius=.4cm, y radius =.4cm];
		\draw[dashed,-] ([yshift=-.15cm,xshift=.2cm]a.east) arc[start angle = 330, end angle = 390, x radius=.4cm, y radius =.4cm];
		\end {tikzpicture} \]
	Consider the following two non-isomorphic non-projective indecomposable representations.
	\[\begin{tikzpicture}[>= angle 60,<-]
	    \node (QR) at (-.7,0){$X =$};
		\node (a) at (0,0) {$\mathbb{K}$};
		\node (b) at (2,0) {$0$};
		\draw[bend right=20] ([yshift=1mm]b.west) to node[above]{$0$} ([yshift=1mm]a.east);
		\draw[bend right=20]  ([yshift=-1mm]a.east) to node[below]{$0$}([yshift=-1mm]b.west);
		\draw[dashed,-] ([yshift=-.15cm, xshift=-.2cm]b.west) arc[start angle = -150, end angle = -210, x radius=.4cm, y radius =.4cm];
		\draw[dashed,-] ([yshift=-.15cm,xshift=.2cm]a.east) arc[start angle = 330, end angle = 390, x radius=.4cm, y radius =.4cm];
		\begin{scope}[xshift=5cm]
		\node (QR) at (-.7,0){$Y =$};
		\node (a) at (0,0) {$0$};
		\node (b) at (2,0) {$\mathbb{K}$};
		\draw[bend right=20] ([yshift=1mm]b.west) to node[above]{$0$} ([yshift=1mm]a.east);
		\draw[bend right=20]  ([yshift=-1mm]a.east) to node[below]{$0$}([yshift=-1mm]b.west);
		\draw[dashed,-] ([yshift=-.15cm, xshift=-.2cm]b.west) arc[start angle = -150, end angle = -210, x radius=.4cm, y radius =.4cm];
		\draw[dashed,-] ([yshift=-.15cm,xshift=.2cm]a.east) arc[start angle = 330, end angle = 390, x radius=.4cm, y radius =.4cm];
		\end{scope}
		\end {tikzpicture} \]
		Then $X \Resleq Y$, as $Y \cong \Ker(\PP_{X}^0 \twoheadrightarrow X)$, and $Y \Resleq X$ for a similar reason. 
\end{ex}

In the previous two examples, the Auslander--Reiten quiver of each gentle quiver is not acyclic. We will prove that,  for gentle trees $(Q,R)$,  $\THomleq$ and $\Resleq$ are order relations.

\begin{lemma} \label{lem:THomOrderTrees}
	Let $(Q,R)$ be a gentle tree. Then $\THomleq$ is an order on $\ind(Q,R)$. Moreover, $\AR(Q,R)$ is the Hasse diagram of $(\ind (Q,R),\THomleq)$.
\end{lemma}

\begin{proof}
First of all, as $\Homleq$ is reflexive, we already know that $\THomleq$ is reflexive and transitive. Moreover, $\THomleq$ is antisymmetric by \cref{thm:acyclicAR}.

By \cref{prop:homtree},  for any $X,Y \in \ind(Q,R)$, there exists at most one nonzero nonisomorphism, up to a multiplication by a scalar, in $\Hom(X,Y)$; meaning that, in $\AR(Q,R)$ there is at most one arrow from $X$ to $Y$.

Assume that we have $X,Y \in \ind(Q,R)$ such that $X\ncong Y$ and $X \THomleq Y$ is a cover relation. Then $X \Homleq Y$, and the nonzero nonisomorphism $f: X \longrightarrow Y$ must be irreducible. So we have an arrow from $X$ to $Y$ in the Auslander-Reiten quiver of $(Q,R)$. Conversely, if we have an arrow from $X$ to $Y$ in the Auslander-Reiten quiver, it is straightforward to verify that $X \THomleq Y$ and that this relation is a cover. We get the desired result.
\end{proof}

\begin{lemma} \label{lem:THominAlgo}
    Let $(Q,R)$ be a gentle tree. Let $\mathcal{X}^0 = \mathcal{X} \subset \pmb{\ind \setminus \proj}(Q,R)$, and, for $i \in \mathbb{N}$, denote by $\mathcal{X}^i$ the subset of $\pmb{\ind \setminus \proj}(Q,R)$ obtained at the $i$th iteration of \cref{algo:inductivRes}. Then, for any $Y \in \mathcal{X}^{i+1} \setminus \proj(Q,R)$, there exists $X \in \mathcal{X}^{i} \setminus \proj(Q,R)$ such that $Y \THomleq X$.
\end{lemma}

\begin{proof}
	Fix $i \in \mathbb{N}$. Consider $Y \in \mathcal{X}^{i+1}$. If $Y \in \mathcal{X}^i$, then we are done. Otherwise, we have that $Y \in \mathcal{X}^{i+1} \setminus \mathcal{X}^i$. Following \cref{algo:inductivRes}, $Y$ can appear in at most two different cases. 
	\begin{enumerate}
	    \item Assume that $Y$ is a summand of $F$ which is in a short exact sequence \[\begin{tikzcd}
	A & F & B
	\arrow["f",tail, from=1-1, to=1-2]
	\arrow["g",two heads, from=1-2, to=1-3]
\end{tikzcd},\] with $A,B \in \mathcal{X}^i$. If $B \in \proj(Q,R)$, then the short exact sequence splits and we get that $F \cong A \oplus B$. So $Y \cong A \in \mathcal{X}^i$ or $Y \cong B \in \mathcal{X}^i$. Otherwise, $B \in \mathcal{X}^i \setminus \proj(Q,R)$, and as $g$ is a nonzero minimal epimorphism, we know that $Y \Homleq B$ and, so $Y \THomleq B$;
	    \item Assume that $Y$ is a summand of $K$ in a short exact sequence \[\begin{tikzcd}
	K & {\displaystyle \bigoplus_{\ell=1}^p A_\ell} & B
	\arrow["f",tail, from=1-1, to=1-2]
	\arrow["g",two heads, from=1-2, to=1-3]
\end{tikzcd},\] for some $p \in \mathbb{N}^*$, with $B,A_1,\ldots, A_p \in \mathcal{X}^i$. If $B \in \proj(Q,R)$, then the short exact sequence splits, and $K$ is isomorphic to a summand of $\bigoplus_{i=1}^p A_p$. So $Y \cong A_j \in \mathcal{X}^{i}$ for some $j \in \{1,\ldots,p\}$, which is a contradiction. Otherwise, as $g$ is a minimal epimorphism and not an isomorphism, we have that $f_{|Y}$ is a nonzero morphism. Therefore, $Y \Homleq A_j$ for some $j \in \{1,\ldots,p\}$. 
	\end{enumerate}

\end{proof}

\begin{cor} \label{cor:THominRes}
	Let $(Q,R)$ be a gentle tree. Let $\mathcal{X} \subset \pmb{\ind \setminus \proj}(Q,R)$. Then, for any $Y \in \mathcal{X}^{\Res} \setminus \proj(Q,R)$, there exists $X \in \mathcal{X} \setminus \proj(Q,R)$ such that $Y \THomleq X$. 
\end{cor}

\begin{proof}
    By induction using \cref{lem:THominAlgo} and by transitivity of $\THomleq$.
\end{proof}

We can state a more precise result which highlights a compatibility between $\Resleq$ and $\THomleq$ on $\pmb{\ind \setminus \proj}(Q,R)$.

\begin{cor} \label{cor:ResIntoTHom}
	Let $(Q,R)$ be a gentle tree. Let $X,Y \in \pmb{\ind \setminus \proj}(Q,R)$. If we have $X \Resleq Y$, then $X \THomleq Y$.  
\end{cor}

\begin{prop} \label{prop:ResOrderTrees}
Let $(Q,R)$ be a gentle tree. Then the $\Res$-relation is an order relation on $\pmb{\ind \setminus \proj}(Q,R)$.
\end{prop}

\begin{proof}
	We only have to check that $\Resleq$ is antisymmetric.
	By \cref{cor:ResIntoTHom}, whenever we have $X \Resleq Y$ for some $X,Y \in \pmb{\ind \setminus \proj}(Q,R)$, we have $X \THomleq Y$. By \cref{lem:THomOrderTrees}, $\THomleq$ is an order on $\ind (Q,R)$, and so it is antisymmetric. Thus $\Resleq$ is antisymmetric.
\end{proof}

\begin{remark} \label{rem:charactgentlequivResorder}
This is not a characterization of gentle quivers $(Q,R)$ such that $\Resleq$ is an order on $\pmb{\ind\setminus \proj}(Q,R)$. See \cref{fig:nontrivialexResord} for an example.
\end{remark}
\begin{figure}[!ht]
    \centering
    \begin{tikzpicture}[scale=.8,xshift=-1cm,yshift=-4cm,->, >= angle 60]
        \node at (-1.5,0) {$(Q,R) =$};
            \node (a) at (0,0) {$1$};
		\node (b) at (1,1) {$2$};
		\node (c) at (2,0) {$3$};
        \node (d) at (1,-1) {$4$};
		\draw (a) -- node[above left]{$\alpha$}(b);
        \draw (b) -- node[above right]{$\beta$} (c);
        \draw (c) -- node[below right]{$\gamma$} (d);
        \draw (d) -- node[below left]{$\delta$} (a);
		\draw[dashed,-] ([xshift=-.3cm,yshift=-.2cm]b.south) arc[start angle = 225, end angle = 315, x radius=.5cm, y radius =.5cm];
        \draw[dashed,-] ([xshift=.3cm,yshift=.2cm]d.north) arc[start angle = 45, end angle = 135, x radius=.5cm, y radius =.5cm];

        \begin{scope}[scale=1.2,-,xshift=4cm,yshift=-1cm]
            \node  (a) at (0,0) {\scalebox{.6}{$\begin{tikzpicture}[baseline={(0,-.1)},scale=0.2]
									\node at (0,0){$2$};
								\end{tikzpicture}$}};
						\node (b) at (0,1) {\scalebox{.6}{$\begin{tikzpicture}[baseline={(0,-.1)},scale=0.2]
									\node at (0,0){$1$};
								\end{tikzpicture}$}};
						\node (c) at (0,2)  {\scalebox{.6}{$\begin{tikzpicture}[baseline={(0,-.1)},scale=0.2]
									\node at (0,0){$4$};
                                    \node at (1,-1){$1$};
								\end{tikzpicture}$}};
						\node (d) at (1,0)  {\scalebox{.6}{$\begin{tikzpicture}[baseline={(0,-.1)},scale=0.2]
									\node at (0,0){$4$};
								\end{tikzpicture}$}};
						\node (g) at (1,1) {\scalebox{.6}{$\begin{tikzpicture}[baseline={(0,-.1)},scale=0.2]
									\node at (0,0){$3$};
								\end{tikzpicture}$}};
						\node (h) at (1,2) {\scalebox{.6}{$\begin{tikzpicture}[baseline={(0,-.1)},scale=0.2]
									\node at (0,0){$2$};
                                    \node at (-1,-1){$3$};
								\end{tikzpicture}$}};
						
						\draw[line width=0.7mm,black] (a) -- (b);
						\draw[line width=0.7mm,black] (b) -- (c);
						\draw[line width=0.7mm,black] (d) -- (g);
						\draw[line width=0.7mm,black] (g) -- (h);

						\node at (.5,-.5) {$(\pmb{\ind \setminus \proj}(Q,R), \Resleq)$};
        \end{scope}
		\end{tikzpicture} 
    \caption{An example of a gentle quiver $(Q,R)$ which is not a gentle tree (left) for which $\protect\Resleq$ is an order relation on $\protect\pmb{\ind \setminus \proj}(Q,R)$ (Hasse diagram on the right).}
    \label{fig:nontrivialexResord}
\end{figure}

\begin{definition}
Let $(Q,R)$ be a gentle quiver. Whenever $(\pmb{\ind \setminus \proj}(Q,R), \Resleq)$ is a poset, we call it \new{the resolving poset of $(Q,R)$}.
\end{definition}

We will see, in what follows and in future work, that this poset plays a crucial role in constructing the lattice of resolving subcategories of a gentle tree (see \cite{DS252}).

\subsection{Resolving posets and posets of resolving subcategories}
\label{ss:ResposetAndSubcat}

Given a gentle quiver $(Q,R)$, we denote by $\ResOrd(Q,R)$ the set of all resolving subcategories of $\rep(Q,R)$. 

\begin{prop} \label{prop:ResOrdPoset} For any representation-finite gentle quiver $(Q,R)$, the pair $(\ResOrd(Q,R), \subseteq)$ is a finite complete lattice.
\end{prop}

\begin{proof}
    It follows from the fact $\rep(Q,R)$ is Krull-Schmidt and representation finite.
\end{proof}

\begin{prop} \label{prop:Posetmorph} Let $(Q,R)$ be a gentle tree. Then the map \[\underline{\Res}: \left\{ \begin{matrix}
(\ResOrd(Q,R), \subseteq) & \longrightarrow & (\mathscr{J}(\pmb{\ind \setminus \proj}(Q,R), \Resleq), \subseteq) \\
\mathscr{R} & \longmapsto & \pmb{\ind \setminus \proj}(\mathscr{R})
\end{matrix} \right.\]
is an injective poset morphism between the poset of resolving subcategories and the ideals of the resolving poset of (Q,R).
\end{prop}

\begin{proof}
    It follows from \cref{lem:resclosop},  since,for any $\mathscr{R} \in \ResOrd(Q,R)$, we have that $(\pmb{\ind \setminus \proj}(\mathscr{R}))^{\Res} = \pmb{\ind \setminus \proj}(\mathscr{R})$.
\end{proof}

\begin{lemma} \label{lem:coveredbyResX}
    Let $(Q,R)$ be a gentle tree, and $X \in \pmb{\ind \setminus \proj}(Q,R)$. Set $\mathcal{Y}_X = \pmb{\ind \setminus \proj}(\Res(X)) \setminus \{X\}$. Then $\add(\mathcal{Y}_X \cup \proj(Q,R)) \in \ResOrd(Q,R)$.
\end{lemma}

\begin{proof}
    The case $\mathcal{Y}_X = \varnothing$ is obvious.
    Assume that $\mathcal{Y}_X \neq \varnothing$. By \cref{prop:Posetmorph}, we get that $X^{\Res}$ is the ideal generated by $X$ in $(\pmb{\ind \setminus \proj}(Q,R), \Resleq)$. Therefore, $\mathcal{Y}_X$ is also an ideal of $(\pmb{\ind \setminus \proj}(Q,R), \Resleq)$. We must show that $(\mathcal{Y}_X)^{\Res} = \mathcal{Y}_X$.
    Assume, by contradiction, that $(\mathcal{Y}_X)^{\Res} \supsetneq \mathcal{Y}_X$. There exists $Z \in (\mathcal{Y}_X)^{\Res} \setminus \mathcal{Y}_X$, and, by \cref{lem:resclosop}, we have that $\{X\} \cup \mathcal{Y}_X = X^{\Res} \supseteq (\mathcal{Y}_X)^{\Res}$. Therefore, $Z \cong X$.
    
    Moreover, using \cref{cor:THominRes}, as $Z \in (\mathcal{Y}_X)^{\Res}$, there exists $Y \in \mathcal{Y}_X$ such that $Z \THomleq Y$. By \cref{cor:ResIntoTHom} since $Y \Resleq X$, we get that $X \cong Z \THomleq Y \THomleq X$. Therefore, $X \cong Y \in \mathcal{Y}_X$, which raises a contradiction.
\end{proof}

\begin{theorem} \label{thm:Monoareallthejoinirred} Let $(Q,R)$ be a gentle tree. Then join-irreducible elements in $(\ResOrd(Q,R), \subseteq)$ are exactly given by monogeneous resolving subcategories.
\end{theorem}

\begin{proof}
    Let $\mathscr{R} \in \ResOrd(Q,R)$ be a join-irreducible element. We consider $\mathcal{X}$ a subset of non-projective indecomposable representations in $\mathscr{R}$. If $\mathcal{X} = \varnothing$, then $\mathscr{R} = \proj(Q,R)$ which is absurd as it is minimal. Thus, we can write that $\mathscr{R} = \bigvee_{X \in \mathcal{X}} \Res(X)$ as $\mathscr{R} = \Res(\mathcal{X})$. By hypothesis, $\mathscr{R}$ is monogeneous.
    
    We still have to show that monogeneous resolving subcategories are join-irreducible. Let $X \in \pmb{\ind \setminus \proj}(Q,R)$. Set $\mathcal{Y}_X = \pmb{\ind \setminus \proj}(\Res(X)) \setminus \{X\}$. By \cref{lem:coveredbyResX}, we showed that $\Res(X)$ covers $\Res(\mathcal{Y}_X) = \add(\mathcal{Y}_X \cup \proj(Q,R))$ in $(\ResOrd(Q,R), \subseteq)$. Indeed, this is the unique $\add(\mathcal{Y}_X \cup \proj(Q,R))$ resolving subcategory covered by $\Res(X)$ by construction. Thus, by \cref{prop:JoinIrreducible}, we get that $\Res(X)$ is join-irreducible in $(\ResOrd(Q,R), \subseteq)$.
\end{proof}

	\section{Surface model for the representations of gentle algebras}
	\label{sec:Surface}
	\pagestyle{plain}

In this section, we introduce a variant of the topological model developed in \cite{BCS21,OPS18,PPP18}, which is useful to study gentle quivers and their representations. This geometric model will make it easier to compute resolving closures for any subset of $\ind \setminus \proj(Q,R)$. The section is organized as follows: \cref{ss:surfdissec} defines the combinatorial shell of dissected surfaces, which is linked bijectively in \cref{ss:gentlesurf} to gentle quivers. Modules are seen as accordions through \cref{ss:Accordionsandindec}, allowing us to translate morphisms and extensions between indecomposables in \cref{ss:GeomKerandExt} into topological notions. To compute resolving subcategories, we need to consider slight changes to the geometric model detailed in \cref{ss:ProjVariant}, and we define a first set of tools that gives us the combinatorial insight we needed to improve Takahashi's algorithm.


\subsection{Surfaces and dissections} \label{ss:surfdissec}

\begin{conv}
Whenever we refer to \new{surfaces}, we will always mean an oriented compact surface $\pmb{\Sigma}$ with a finite number (which may be zero) of removed open discs .  Denote the boundary of $\pmb{\Sigma}$ by $\partial \pmb{\Sigma}$.  Such a surface is determined, up to homeomorphism, by its genus and by the number of connected components of $\partial \pmb{\Sigma}$; we will refer to these connected components as the \new{boundary components} of $\pmb{\Sigma}$.
\end{conv}

\begin{definition} \label{def:marksurf}
A \new{marked surface} is a pair $(\pmb{\Sigma},\mathcal{M})$, where $\pmb{\Sigma}$ is a surface and $\mathcal{M}$ is a finite set of \new{marked points} of $\pmb{\Sigma}$ such that $\mathcal{M}$ admits a bipartition $\{\mathcal{M}_{\gpoint},\mathcal{M}_{\rpoint}\}$ satisfying that for each boundary component of $\pmb{\Sigma}$: 
\begin{enumerate}[label = $\bullet$, itemsep=0.1em]
	\item there is at least one marked point of each color;
	\item marked points in $\mathcal{M}_{\gpoint}$ and $\mathcal{M}_{\rpoint}$ alternate.
\end{enumerate}
\end{definition}

Note that marked points are also allowed in the interior of the surface.

\begin{definition} \label{def:arcsanddissec} Let $(\pmb{\Sigma},\mathcal{M})$ be a marked surface.
\begin{enumerate}[label=$\bullet$,itemsep=1mm]
\item A \new{$\gpoint$-arc} is a curve on $(\pmb{\Sigma},\mathcal{M})$ joining two points in $\mathcal{M}_{\gpoint}$; it is a continuous map from the closed interval $[0,1]$ to $\pmb{\Sigma}$ with endpoints in $\mathcal{M}_{\gpoint}$, and such that the image of its interior is disjoint from $\mathcal{M}$ up to homotopy.

\item A $\gpoint$-arc is said to be \new{simple} if it does not intersect itself (except perhaps at its endpoints). 

\item A \new{$\gpoint$-dissection of $(\pmb{\Sigma}, \mathcal{M})$} is a collection $\Delta^{\gpoint}$ of pairwise non-intersecting simple $\gpoint$-arcs which cut the surface into polygons (that is to say, into simply-connected regions), called the \new{cells} of the dissection. The triplet $(\pmb{\Sigma}, \mathcal{M}, \Delta^{\gpoint})$  is called a  \new{$\gpoint$-dissected marked surface}. 

\item A cell is said to be \new{large} if the polygon has more more than three edges.

\item A dissection is said to be \new{dualizable} if each cell of the dissection contains exactly one marked point in $\mathcal{M}_{\rpoint}$.

\item We define \new{$\rpoint$-arcs} and \new{$\rpoint$-dissections} dually.

\item Given a $\gpoint$-dissection or a $\rpoint$-dissection $\Delta$ of a marked surface $(\pmb{\Sigma}, \mathcal{M})$, we denote by $\pmb{\Gamma}(\Delta)$ the set of cells of $(\pmb{\Sigma}, \mathcal{M},\Delta)$ and for $C\in\pmb{\Gamma}(\Delta)$ we denote by $C_0$ the endpoints of the arcs of $\partial C$. 

\end{enumerate}
\end{definition}


\begin{ex}\label{exam::three-surfaces}
	Here are examples of $\gpoint$-dissections of an annulus (or a sphere with two boundary components), and a torus with one boundary component.
	
	\begin{figure}[!ht]
		\begin{center}
		     \begin{tikzpicture}[mydot/.style={
					circle,
					thick,
					fill=white,
					draw,
					outer sep=0.5pt,
					inner sep=1pt
				}]
				\tkzDefPoint(0,0){O}\tkzDefPoint(0,1.5){1}
				\tkzDefPointsBy[rotation=center O angle 180](1){2}
				\tkzDefPoint(0.7,0){3} \tkzDefPoint(1.1,0){8} \tkzDefPoint(-1.1,0){9}
				\tkzDefPointsBy[rotation=center O angle 90](3,4,5){4,5,6}
				\tkzDrawCircle[line width=0.7mm,black](O,1)
				\tkzDrawCircle[line width=0.7mm,black,fill=gray!40](O,3)
				\tkzDrawPoints[fill =red,size=4,color=red](2,3,5)

				\draw[line width=0.5mm,dark-green](1) edge (4);
				\draw[line width=0.5mm,bend left = 60,dark-green](4) edge (8);
				\draw[line width=0.5mm,bend left = 60,dark-green](8) edge (6);
				\draw[line width=0.5mm,bend right = 60,dark-green](4) edge (9);
				\draw[line width=0.5mm,bend right =60,dark-green](9) edge (6);
				\tkzDrawPoints[size=4,color=dark-green,mydot](1,4,6)
			\end{tikzpicture}
			\qquad \begin{tikzpicture}[mydot/.style={
					circle,
					thick,
					fill=white,
					draw,
					outer sep=0.5pt,
					inner sep=1pt
				}]
				\tkzDefPoint(0,0){O}\tkzDefPoint(1.4,1.4){1}
				\tkzDefPointsBy[rotation=center O angle 90](1,2,3){2,3,4}
				\tkzDefPoint(0.9,0.9){5}
				\tkzDefPointsBy[rotation=center O angle 90](5,6,7){6,7,8}
				
				\draw[line width=0.7mm] (-1.4142,-1.4142) rectangle (1.4142,1.4142);
				\tkzDefPointsBy[rotation=center 1 angle 45](5){9}
				\tkzDrawSector[rotate,line width=0.5mm, black, fill = gray!40](1,9)(-90)
				\tkzDefPointsBy[rotation=center 2 angle 45](6){10}
				\tkzDrawSector[rotate,line width=0.5mm, black, fill = gray!40](2,10)(-90)
				\tkzDefPointsBy[rotation=center 3 angle 45](7){11}
				\tkzDrawSector[rotate,line width=0.5mm, black, fill = gray!40](3,11)(-90)
				\tkzDefPointsBy[rotation=center 4 angle 45](8){12}
				\tkzDrawSector[rotate,line width=0.5mm, black, fill = gray!40](4,12)(-90)
				
				\tkzDefPoint(1.4,0){13}
				\tkzDefPoint(-1.4,0){14}
				\tkzDefPoint(0,1.4){15}
				\tkzDefPoint(0,-1.4){16}
				
				\tkzDrawPoints[fill =red,size=4,color=red](5,7)
				\draw[line width=0.5mm,dark-green](6) edge (8);
				\draw[line width=0.5mm,bend left=40,dark-green](6) edge (14);
				\draw[line width=0.5mm,bend right=40,dark-green](6) edge (15);
				\draw[line width=0.5mm,bend left=40,dark-green](8) edge (13);
				\draw[line width=0.5mm,bend right=40,dark-green](8) edge (16);
				\tkzDrawPoints[size=4,color=dark-green,mydot](6,8)
			\end{tikzpicture}
			\caption{\label{fig:dissections} Examples of a $\gpoint$-dissection of an annulus (left), and of a torus with one boundary component (right).These dissections are dualizable}
		\end{center}
	\end{figure}
	
\end{ex}

\subsection{The gentle quiver of a dissection} \label{ss:gentlesurf}
Let $(\pmb{\Sigma}, \mathcal{M})$ be a marked surface.

\begin{definition} \label{def:GentleQuivfromsurf}
Let $\Delta^{\gpoint}$ be a $\gpoint$-dissection of $(\pmb{\Sigma},\mathcal{M})$. We define a quiver $\opQ(\Delta^{\gpoint})$ and a set of relations $\opR(\Delta^{\gpoint})$ as follows:
\begin{enumerate}[label = $(\roman*)$,itemsep=0.25em]
	\item vertices of $\opQ(\Delta^{\gpoint})$ correspond bijectively to arcs in $\Delta^{\gpoint}$; 
	\item \label{ii} for each configuration where there exists a common endpoint of a pair of arcs $i$ and $j$ around which $i$ comes immediately after $j$ in clockwise orientation in $\Delta^{\gpoint}$, as drawn on the left in \cref{fig:rulesquiverdissections}, there is an arrow $j\to i$ in $\opQ(\Delta^{\gpoint})$;
	\item \label{iii}for each configuration as drawn on the right in \cref{fig:rulesquiverdissections} in $\Delta^{\gpoint}$ with corresponding arrows $j\xrightarrow{\alpha} i$ and $k\xrightarrow{\beta} j$, the path $\alpha \beta$ is a relation in $\opR(\Delta^{\gpoint})$.
\end{enumerate}
\begin{figure}[!ht] 
	\centering
	\begin{tikzpicture}[mydot/.style={
				circle,
				thick,
				fill=white,
				draw,
				outer sep=0.5pt,
				inner sep=1pt
			}, fl/.style={->,>=latex}]
			\tkzDefPoint(0,0){O}\tkzDefPoint(0,1.5){1}
			\tkzDefPointsBy[rotation=center O angle 120](1,2){2,3}
			\tkzDefMidPoint(1,2)
			\tkzGetPoint{i}
			\tkzDefMidPoint(1,3)
			\tkzGetPoint{j}
			
			\draw[line width=0.5mm,dark-green](1) edge (2);
			\draw[line width=0.5mm,dark-green](1) edge (3);
			\draw[->, >= angle 60,line width=0.7mm,black] (j) -- (i);
			\tkzDrawPoints[size=4,color=black](i,j);\tkzDrawPoints[size=4,color=dark-green,mydot](1);
			\tkzLabelPoints[left](i)
			\tkzLabelPoints[right](j)
		\end{tikzpicture} \quad
		\begin{tikzpicture}[mydot/.style={
				circle,
				thick,
				fill=white,
				draw,
				outer sep=0.5pt,
				inner sep=1pt
			}]
			\tkzDefPoint(0,0){O}\tkzDefPoint(-1,1.5){1} \tkzDefPoint(1,1.5){2} \tkzDefPoint(-3,-1){3} \tkzDefPoint(3,-1){4}
			\tkzDefPoint(0,1.1){5}
			
			\tkzDefMidPoint(1,3)
			\tkzGetPoint{i}
			\tkzDefMidPoint(1,2)
			\tkzGetPoint{j}
			\tkzDefMidPoint(2,4)
			\tkzGetPoint{k}
			\tkzDefMidPoint(i,j)
			\tkzGetPoint{a}
			\tkzDefMidPoint(j,k)
			\tkzGetPoint{b}
			\draw[line width=0.5mm,dark-green](1) edge (2);
			\draw[line width=0.5mm,dark-green](1) edge (3);
			\draw[line width=0.5mm,dark-green](2) edge (4);
			\draw[->, >= angle 60,line width=0.7mm,black] (j) -- (i);
			\draw[->, >= angle 60,line width=0.7mm,black] (k) -- (j);
			\tkzDrawArc[angles,line width=0.7mm,black,dashed](j,5)(-150,-30)
			\tkzDrawPoints[size=4,color=black](i,j,k);\tkzDrawPoints[size=4,color=dark-green,mydot](1,2);
			\tkzLabelPoints[left](i)
			\tkzLabelPoints[above](j)
			\tkzLabelPoints[right](k)
			\tkzLabelPoint[below](a){$\alpha$}
			\tkzLabelPoint[below](b){$\beta$}
		\end{tikzpicture} 
		 \caption{\label{fig:rulesquiverdissections} Drawings representing the construction rules of $\opQ(\Delta^{\gpoint})$ and $\opR(\Delta^{\gpoint})$ as in \cref{def:GentleQuivfromsurf}; \ref{ii} on the left and \ref{iii} on the right.}
\end{figure}
The pair $(\opQ(\Delta^{\gpoint}), \opR(\Delta^{\gpoint}))$ is the \new{bounded quiver associated} to $(\pmb{\Sigma}, \mathcal{M}, \Delta^{\gpoint})$.
\end{definition}

\begin{conv} \label{conv:clockwise}
We orient our marked surface clockwise around the green marked points. 
\end{conv}

\begin{ex}\label{ex:quiverofsurfaces}
	See in \cref{fig:twosurfquiver} the quivers with relations associated with the dissections shown in  \cref{fig:dissections}.
	
	\begin{figure}[!ht]
		\centering 
		\begin{tikzpicture}[mydot/.style={
				circle,
				thick,
				fill=white,
				draw,
				outer sep=0.5pt,
				inner sep=1pt,
			},scale = 1.2]
			\tkzDefPoint(0,0){O}\tkzDefPoint(0,1.5){1}
			\tkzDefPointsBy[rotation=center O angle 180](1){2}
			\tkzDefPoint(0.7,0){3} \tkzDefPoint(1.1,0){8} \tkzDefPoint(-1.1,0){9}
			\tkzDefPointsBy[rotation=center O angle 90](3,4,5){4,5,6}
			
			\tkzDrawCircle[line width=0.7mm,black](O,1)
			\tkzDrawCircle[line width=0.7mm,black,fill=gray!40](O,3)

			\draw[line width=0.5mm,dark-green](1) edge (4);
			\draw[line width=0.5mm,bend left = 60,dark-green](4) edge (8);
			\draw[line width=0.5mm,bend left = 60,dark-green](8) edge (6);
			\draw[line width=0.5mm,bend right = 60,dark-green](4) edge (9);
			\draw[line width=0.5mm,bend right =60,dark-green](9) edge (6);
			\tkzDrawPoints[size=4,color=dark-green,mydot](1,4,6)
			\tkzDrawPoints[fill =red,size=4,color=red](2,3,5)
			
			\tkzDefPoint(0,1.1){a}
            \node at (0.2,1.3){$1$};
			\tkzDefPoint(1.1,0){b}
            \node at (.9,0){$2$};
			\tkzDefPoint(-1.1,0){c}
            \node at (-.9,0){$3$};
			\tkzDrawPoints[size=4,fill=black,color=black](a,b,c)

			\tkzDefPoint(0,-1.2){D}
			\tkzCircumCenter(b,c,D)\tkzGetPoint{O1}
			\tkzDrawArc[<-, >= angle 60,line width=0.7mm,black](O1,c)(b)
			\tkzDefPoint(0.9,0.9){E}
			\tkzCircumCenter(a,b,E)\tkzGetPoint{O1}
			\tkzDrawArc[<-, >= angle 60,line width=0.7mm,black](O1,b)(a)
			\tkzDefPoint(-0.9,0.9){F}
			\tkzCircumCenter(a,c,F)\tkzGetPoint{O1}
			\tkzDrawArc[<-, >= angle 60,line width=0.7mm,black](O1,a)(c)
			
			\tkzDefPoint(-1.15,0.4){A}
			\tkzDefPoint(-1.1,-0.4){B}
			\tkzDefPoint(-1,0){C}
			\tkzDrawArc[line width=0.7mm,black,dashed](C,A)(B)
			\tkzDefPoint(1.15,0.4){A}
			\tkzDefPoint(1.1,-0.4){B}
			\tkzDefPoint(1,0){C}
			\tkzDrawArc[line width=0.7mm,black,dashed](C,B)(A)

            \begin{scope}[scale=.8,xshift=-1cm,yshift=-4cm,->, >= angle 60]
            \node (a) at (1,1.7) {$1$};
		\node (b) at (2,0) {$2$};
		\node (c) at (0,0) {$3$};
		\draw (a) -- (b);
        \draw (b) -- (c);
        \draw (c) -- (a);
		\draw[dashed,-] ([xshift=-.4cm,yshift=.3cm]b.north) arc[start angle = 120, end angle = 200, x radius=.4cm, y radius =.4cm];
		\draw[dashed,-] ([xshift=.6cm,yshift=-.2cm]c.north) arc[start angle = -20, end angle = 60, x radius=.4cm, y radius =.4cm];
        \end{scope}
		\begin{scope}[xshift=4cm]
		    \tkzDefPoint(0,0){O}\tkzDefPoint(1.4,1.4){1}
			\tkzDefPointsBy[rotation=center O angle 90](1,2,3){2,3,4}
			\tkzDefPoint(0.9,0.9){5}
			\tkzDefPointsBy[rotation=center O angle 90](5,6,7){6,7,8}
			
			\draw[line width=0.7mm] (-1.4142,-1.4142) rectangle (1.4142,1.4142);
			\tkzDefPointsBy[rotation=center 1 angle 45](5){9}
			\tkzDrawSector[rotate,line width=0.5mm, black, fill = gray!40](1,9)(-90)
			\tkzDefPointsBy[rotation=center 2 angle 45](6){10}
			\tkzDrawSector[rotate,line width=0.5mm, black, fill = gray!40](2,10)(-90)
			\tkzDefPointsBy[rotation=center 3 angle 45](7){11}
			\tkzDrawSector[rotate,line width=0.5mm, black, fill = gray!40](3,11)(-90)
			\tkzDefPointsBy[rotation=center 4 angle 45](8){12}
			\tkzDrawSector[rotate,line width=0.5mm, black, fill = gray!40](4,12)(-90)
			
			\tkzDefPoint(1.4,0){13}
			\tkzDefPoint(-1.4,0){14}
			\tkzDefPoint(0,1.4){15}
			\tkzDefPoint(0,-1.4){16}

			\draw[line width=0.5mm,dark-green](6) edge (8);
			\draw[line width=0.5mm,bend left=40,dark-green](6) edge (14);
			\draw[line width=0.5mm,bend right=40,dark-green](6) edge (15);
			\draw[line width=0.5mm,bend left=40,dark-green](8) edge (13);
			\draw[line width=0.5mm,bend right=40,dark-green](8) edge (16);
			\tkzDrawPoints[fill =red,size=4,color=red](5,7)
			\tkzDrawPoints[size=4,color=dark-green,mydot](6,8)
			
			\tkzDefPoint(0,0){a}
            \node at (-.3,.2){$2$};
			\tkzDefPoint(1.4,0){b}
            \node[blue] at (1.6,0){$3$};
			\tkzDefPoint(-1.4,0){c}
            \node[blue] at (-1.6,0){$3$};
			\tkzDefPoint(0,1.4){d}
            \node[orange] at (0,1.7){$1$};
			\tkzDefPoint(0,-1.4){e}
            \node[orange] at (0,-1.7){$1$};
			\tkzDrawPoint[size=4,fill=black,color=black](a)

			\draw[<-, >= angle 60,line width=0.7mm,black](b) -- (a);
			\draw[<-, >= angle 60,line width=0.7mm,black](c) -- (a);
			\draw[<-, >= angle 60,line width=0.7mm,black](a) -- (d);
			\draw[<-, >= angle 60,line width=0.7mm,black](a) -- (e);
			\tkzDrawPoints[rectangle,size=4,fill=blue,color=blue](b,c)
			\tkzDrawPoints[cross,size=4,fill=orange,color=orange,line width=0.8mm](d,e)
			
			\tkzDefPoint(0.5,0){B}
			\tkzDefPoint(0,0.5){D}
			\tkzDefPoint(-0.5,0){C}
			\tkzDefPoint(0,-0.5){E}
			\tkzDrawArc[line width=0.7mm,black,dashed](a,B)(D)
			\tkzDrawArc[line width=0.7mm,black,dashed](a,C)(E)

        \begin{scope}[scale=.8,xshift=-2cm,yshift=-3.5cm,<-, >= angle 60]
            \node (a) at (0,0) {${\color{orange}{1}}$};
		\node (b) at (2,0) {$2$};
		\node (c) at (4,0) {${\color{blue}{3}}$};
		\draw ([yshift=1mm]b.west)--node[above]{$\beta_{1}$}([yshift=1mm]a.east);
		\draw  ([yshift=-1mm]b.west)--node[below]{$\alpha_{1}$}([yshift=-1mm]a.east);
		\draw ([yshift=1mm]c.west)--node[above]{$\alpha_{2}$}([yshift=1mm]b.east);
		\draw  ([yshift=-1mm]c.west)--node[below]{$\beta_{2}$}([yshift=-1mm]b.east);
		\draw[dashed,-] ([xshift=.4cm]b.north) arc[start angle = 0, end angle = 180, x radius=.4cm, y radius =.2cm];
		\draw[dashed,-] ([xshift=-.4cm]b.south) arc[start angle = 180, end angle = 360, x radius=.4cm, y radius =.2cm];
        \end{scope}
		\end{scope}
		\end{tikzpicture} 
		\caption{\label{fig:twosurfquiver} Quivers with relations arising from the dissected surfaces of \cref{fig:dissections}}
	\end{figure} 
	
	To construct the quiver with relations associated with the example on the right, we have to identify the two square points (blue) and the two cross points (orange).
	\end{ex}
	The following result motivates our interest in the geometric model.
	\begin{theorem}[\cite{BCS21, OPS18, PPP18}] $ $
		\begin{enumerate}[label = $(\roman*)$, itemsep=0.5em]
			\item If $\Delta^{\gpoint}$ is a $\gpoint$-dissection of the marked surface $(\pmb{\Sigma},\mathcal{M})$ such that there are no $\gpoint$-marked points in the interior of $\pmb{\Sigma}$, then $\left(\opQ(\Delta^{\gpoint}), \opR(\Delta^{\gpoint})\right)$ is a gentle quiver.
			
			\item For any gentle quiver $(Q,R)$, there exists a marked surface $(\pmb{\Sigma},\mathcal{M})$ with no $\gpoint$-marked point in the interior of $\pmb{\Sigma}$ and with a $\gpoint$-dissection $\Delta^{\gpoint}$ such that $(Q,R)$ is isomorphic to $\left(\opQ(\Delta^{\gpoint}), \opR(\Delta^{\gpoint})\right)$.  Moreover, $(\pmb{\Sigma},\mathcal{M})$ and $\Delta^{\gpoint}$ are unique up to oriented homeomorphism of surfaces and homotopy of arcs.
		\end{enumerate}
	\end{theorem}
	An explicit recipe to construct $(\pmb{\Sigma},\mathcal{M}, \Delta^{\gpoint})$ from the gentle quiver $(Q,R)$ is given in \cite[Section 4.2]{PPP18}. In what follows, we set $\Surf(Q,R) = (\pmb{\Sigma},\mathcal{M},\Delta^{\gpoint})$.
	
	\subsection{Accordions and indecomposable representations}
	\label{ss:Accordionsandindec}
	
	Let $(\pmb{\Sigma},\mathcal{M}, \Delta^{\gpoint})$ be a $\gpoint$-dissected marked surface such that $\mathcal{M}_{\gpoint} \subset \partial \pmb{\Sigma}$.

Representations of gentle quivers are related to certain arcs on the associated $\gpoint$-dissected marked surface: these are known as \emph{accordions}.
    
	\begin{definition}[{\cite{OPS18}}] \label{def:accordions}
	    Let $(\pmb{\Sigma},\mathcal{M}, \Delta^{\gpoint})$ be a $\gpoint$-dissected marked surface such that $\mathcal{M}_{\gpoint} \subset \partial \pmb{\Sigma}$. A $\rpoint$-arc $\delta$ in $(\pmb{\Sigma}, \mathcal{M},\Delta^{\gpoint})$ is an \new{accordion} if it satisfies the following conditions: whenever~$\delta$ enters a cell of~$\Delta^{\gpoint}$ by crossing an $\gpoint$-arc~$\eta$,
		    \begin{enumerate}[label=$(\alph*)$,itemsep=0.5em]
			\item if it leaves the cell, it leaves it by crossing an $\gpoint$-arc~$\zeta$ adjacent to~$\eta$;
			\item the relevant segments of the arcs~$\eta$,~$\zeta$ and~$\delta$ bound a disk that does not contain the unique marked point in~$\mathcal{M}_{\rpoint}$ inside the cell.
		    \end{enumerate}
		See \cref{fig:rulesclosedaccord} for a picture illustrating the rules $(a)$ and $(b)$. 
		\begin{figure}[!ht] 
			\begin{center}
			    \begin{tikzpicture}[mydot/.style={
						circle,
						thick,
						fill=white,
						draw,
						outer sep=0.5pt,
						inner sep=1pt
					}]
					
				\begin{scope}[xshift=0cm]
					\tkzDefPoint(0,0){O}\tkzDefPoint(0,1.5){1}
					\tkzDefPointsBy[rotation=center O angle 120](1,2){2,3}
					\tkzDefMidPoint(1,2)
					\tkzGetPoint{i}
					\tkzDefMidPoint(1,3)
					\tkzGetPoint{j}
					\tkzDefPoint(-1.5,1.5){4}
					\tkzDefPoint(1.5,1.5){5}
					
					\draw[line width=0.5mm,dark-green](1) edge (2);
					\draw[line width=0.5mm,dark-green](1) edge (3);
					\draw[line width=0.9mm,red] (i) edge (j);
					\draw[line width=0.9mm, bend right=50,red] (4) edge (i);
					\draw[line width=0.9mm,bend left=50,red] (5) edge (j);
					\tkzDrawPoints[size=4,color=dark-green,mydot](1);
					
					\tkzDefPoint(0,-0.5){6}
					\tkzDrawPoints[size=4,color=red](6);
					
					\tkzLabelPoint[below left=0.3](4){{\large $\color{red}{\delta}$}}
					\tkzLabelPoint[below](3){{\large $\color{dark-green}{\eta}$}}
					\tkzLabelPoint[below](2){{\large $\color{dark-green}{\zeta}$}}
				\end{scope} 
				\begin{scope}[xshift=5cm]
					\tkzDefPoint(0,0){O}\tkzDefPoint(0,1.5){1}
					\tkzDefPointsBy[rotation=center O angle 120](1,2){2,3}
					\tkzDefMidPoint(1,2)
					\tkzGetPoint{i}
					\tkzDefMidPoint(1,3)
					\tkzGetPoint{j}
					\tkzDefPoint(-1.5,1.5){4}
					\tkzDefPoint(1.5,1.5){5}
					
					\draw[line width=0.5mm,dark-green](1) edge (2);
					\draw[line width=0.5mm,dark-green](1) edge (3);
					\draw[line width=0.9mm,red] (i) edge (j);
					\draw[line width=0.9mm, bend right=50,red] (4) edge (i);
					\draw[line width=0.9mm,bend left=50,red] (5) edge (j);
					
					\tkzDrawPoints[size=4,color=dark-green,mydot](1);
					
					\tkzDefPoint(0,0.75){6}
					\tkzDrawPoints[size=4,color=red](6);
					
					\tkzLabelPoint[below left=0.3](4){{\large $\color{red}{\delta}$}}
					\tkzLabelPoint[below](3){{\large $\color{dark-green}{\eta}$}}
					\tkzLabelPoint[below](2){{\large $\color{dark-green}{\zeta}$}}
				\end{scope}
				\end{tikzpicture}
				\caption{\label{fig:rulesclosedaccord} Drawings representing the rules that an accordion must satisfy: on the left, $\delta$ satisfies the rule $(a)$ and $(b)$; on the right, $\delta$ does not satisfy the rule $(b)$.} \end{center} \end{figure}
	\end{definition} 

    \begin{remark} \label{rem:BCSOPS} In \cite{OPS18}, accordions are called  \emph{admissible curves}. However, in this paper, we use the accordion terminology adopted in \cite{PPP18}.
    \end{remark}
	
		\begin{figure}[!ht]
			\centering
			 \begin{tikzpicture}[mydot/.style={
					circle,
					thick,
					fill=white,
					draw,
					outer sep=0.5pt,
					inner sep=1pt
				}, scale = 1]
			
			\tkzDefPoint(0,0){O}\tkzDefPoint(0,1.5){1}
				\tkzDefPointsBy[rotation=center O angle 180](1){2}
				\tkzDefPoint(0.7,0){3} \tkzDefPoint(1.1,0){8} \tkzDefPoint(-1.1,0){9}
				\tkzDefPointsBy[rotation=center O angle 90](3,4,5){4,5,6}
				\tkzDrawCircle[line width=0.7mm,black](O,1)
				\tkzDrawCircle[line width=0.7mm,black,fill=gray!40](O,3)
				\tkzDrawPoints[fill =red,size=4,color=red](2,3,5)

				\draw[line width=0.5mm,dark-green](1) edge (4);
				\draw[line width=0.5mm,bend left = 60,dark-green](4) edge (8);
				\draw[line width=0.5mm,bend left = 60,dark-green](8) edge (6);
				\draw[line width=0.5mm,bend right = 60,dark-green](4) edge (9);
				\draw[line width=0.5mm,bend right =60,dark-green](9) edge (6);
				\tkzDrawPoints[size=4,color=dark-green,mydot](1,4,6)
				
				\draw[line width=0.9mm,red,bend right=40] (2) edge (1.35,0);
				\draw[line width=0.9mm,red,bend right=40] (1.35,0) edge (0,1.25);
				\draw[line width=0.9mm,bend right=60,red](0,1.25) edge (5);
				
			\begin{scope}[xshift=5cm]	\tkzDefPoint(0,0){O}\tkzDefPoint(1.4,1.4){1}
				\tkzDefPointsBy[rotation=center O angle 90](1,2,3){2,3,4}
				\tkzDefPoint(0.9,0.9){5}
				\tkzDefPointsBy[rotation=center O angle 90](5,6,7){6,7,8}
				
				\draw[line width=0.7mm] (-1.4142,-1.4142) rectangle (1.4142,1.4142);
				\tkzDefPointsBy[rotation=center 1 angle 45](5){9}
				\tkzDrawSector[rotate,line width=0.5mm, black, fill = gray!40](1,9)(-90)
				\tkzDefPointsBy[rotation=center 2 angle 45](6){10}
				\tkzDrawSector[rotate,line width=0.5mm, black, fill = gray!40](2,10)(-90)
				\tkzDefPointsBy[rotation=center 3 angle 45](7){11}
				\tkzDrawSector[rotate,line width=0.5mm, black, fill = gray!40](3,11)(-90)
				\tkzDefPointsBy[rotation=center 4 angle 45](8){12}
				\tkzDrawSector[rotate,line width=0.5mm, black, fill = gray!40](4,12)(-90)
				
				\tkzDefPoint(1.4,0){13}
				\tkzDefPoint(-1.4,0){14}
				\tkzDefPoint(0,1.4){15}
				\tkzDefPoint(0,-1.4){16}
				
				\tkzDrawPoints[fill =red,size=4,color=red](5,7)
				\draw[line width=0.5mm,dark-green](6) edge (8);
				\draw[line width=0.5mm,bend left=40,dark-green](6) edge (14);
				\draw[line width=0.5mm,bend right=40,dark-green](6) edge (15);
				\draw[line width=0.5mm,bend left=40,dark-green](8) edge (13);
				\draw[line width=0.5mm,bend right=40,dark-green](8) edge (16);
				\tkzDrawPoints[size=4,color=dark-green,mydot](6,8)
				
				\draw[line width=0.9mm,red](-1.4,-0.5) edge (1.4,0.5);
				\draw[line width=0.9mm,bend right=20,red](-1.4,0.5) edge (5);
				\draw[line width=0.9mm,bend left=20,red](7) edge (1.4,-0.5);
			\end{scope}
			\end{tikzpicture}
			\caption{\label{fig:accordandclosedaccord} An example of an accordion over each \mbox{$\gpoint$-dissection} seen in \cref{exam::three-surfaces}.}
		\end{figure}
	
\begin{theorem}[\cite{BCS21,OPS18,PPP182}] \label{thm:GeomandRep} Consider $(\pmb{\Sigma}, \mathcal{M}, \Delta^{\gpoint})$ a $\gpoint$-dissected marked surface such that $\mathcal{M}_{\gpoint} \subset \partial \pmb{\Sigma}$, and let $(\opQ(\Delta^{\gpoint}),\opR(\Delta^{\gpoint}))$ be its associated gentle quiver. Then we have a one-to-one correspondence between accordions on $(\pmb{\Sigma}, \mathcal{M}, \Delta^{\gpoint})$ and string representations of $(\opQ(\Delta^{\gpoint}),\opR(\Delta^{\gpoint}))$.
Explicitly, the string associated with the accordion is built starting from an endpoint and concatenating the arrows associated with each pair of consecutive dissection arcs it crosses.
\end{theorem}

\begin{remark} \label{rmk:Geomandband} We adapted \cref{thm:GeomandRep} as we restrict ourselves to representation-finite gentle quivers. Indeed, \cite{BCS21,OPS18,PPP182} states a result that allows one to interpret all elements of $\ind(Q,R)$ geometrically. In the general setup, there exist non-contractible closed curves, and these curves give rise to bands.
\end{remark}
	
\begin{ex} See \cref{fig:accordandrep} for an illustration of the bijection from accordions in a  $\gpoint$-dissected marked surface $(\pmb{\Sigma}, \mathcal{M}, \Delta^{\gpoint})$ to strings in the associated gentle quiver $(\opQ(\Delta^{\gpoint}), \opR(\Delta^{\gpoint}))$.
\begin{figure}[!ht]
\centering
\begin{tikzpicture}[mydot/.style={
				circle,
				thick,
				fill=white,
				draw,
				outer sep=0.5pt,
				inner sep=1pt
			}, scale = 1.2]
			\tkzDefPoint(0,0){O}\tkzDefPoint(1.4,1.4){1}
			\tkzDefPointsBy[rotation=center O angle 90](1,2,3){2,3,4}
			\tkzDefPoint(0.9,0.9){5}
			\tkzDefPointsBy[rotation=center O angle 90](5,6,7){6,7,8}
			
			\draw[line width=0.7mm] (-1.4142,-1.4142) rectangle (1.4142,1.4142);
			\tkzDefPointsBy[rotation=center 1 angle 45](5){9}
			\tkzDrawSector[rotate,line width=0.5mm, black, fill = gray!40](1,9)(-90)
			\tkzDefPointsBy[rotation=center 2 angle 45](6){10}
			\tkzDrawSector[rotate,line width=0.5mm, black, fill = gray!40](2,10)(-90)
			\tkzDefPointsBy[rotation=center 3 angle 45](7){11}
			\tkzDrawSector[rotate,line width=0.5mm, black, fill = gray!40](3,11)(-90)
			\tkzDefPointsBy[rotation=center 4 angle 45](8){12}
			\tkzDrawSector[rotate,line width=0.5mm, black, fill = gray!40](4,12)(-90)
			
			\tkzDefPoint(1.4,0){13}
			\tkzDefPoint(-1.4,0){14}
			\tkzDefPoint(0,1.4){15}
			\tkzDefPoint(0,-1.4){16}

			\draw[line width=0.5mm,dark-green](6) edge (8);
			\draw[line width=0.5mm,bend left=40,dark-green](6) edge (14);
			\draw[line width=0.5mm,bend right=40,dark-green](6) edge (15);
			\draw[line width=0.5mm,bend left=40,dark-green](8) edge (13);
			\draw[line width=0.5mm,bend right=40,dark-green](8) edge (16);
			\tkzDrawPoints[fill =red,size=4,color=red](5,7)
			\tkzDrawPoints[size=4,color=dark-green,mydot](6,8)
			
			\tkzDefPoint(0,0){a}
			\tkzDefPoint(1.4,0){b}
			\tkzDefPoint(-1.4,0){c}
			\tkzDefPoint(0,1.4){d}
			\tkzDefPoint(0,-1.4){e}
			\tkzDrawPoint[size=4,fill=black,color=black](a)

			\draw[<-, >= angle 60,line width=0.7mm,black](b) -- (a);
			\draw[<-, >= angle 60,line width=0.7mm,black](c) -- (a);
			\draw[<-, >= angle 60,line width=0.7mm,black](a) -- (d);
			\draw[<-, >= angle 60,line width=0.7mm,black](a) -- (e);
			\tkzDrawPoints[rectangle,size=4,fill=blue,color=blue](b,c)
			\tkzDrawPoints[cross,size=4,fill=orange,color=orange,line width=0.8mm](d,e)
			
			\tkzDefPoint(0.5,0){B}
			\tkzDefPoint(0,0.5){D}
			\tkzDefPoint(-0.5,0){C}
			\tkzDefPoint(0,-0.5){E}
			\tkzDrawArc[line width=0.7mm,black,dashed](a,B)(D)
			\tkzDrawArc[line width=0.7mm,black,dashed](a,C)(E)
			\draw[line width=0.9mm,red](-1.4,-0.5) edge (1.4,0.5);
			\draw[line width=0.9mm,bend right=20,red](-1.4,0.5) edge (5);
			\draw[line width=0.9mm,bend left=20,red](7) edge (1.4,-0.5);
			\draw[line width=0.9mm,red,opacity=0.4](5,0.8) edge (7,0.9);
			\draw[line width=0.9mm,red,opacity=0.4](5,1) edge (7,0.9);
			\draw[line width=0.9mm,red,opacity=0.4](5,1) edge (7,1.1);
			\draw[line width=0.9mm,red,opacity=0.4](5,1.2) edge (7,1.1);
			\node at (5,1){$\begin{tikzpicture}[>= angle 60, <-]
		\node (a) at (0,0) {${\color{orange}{1}}$};
		\node (b) at (2,0) {$2$};
		\node (c) at (4,0) {${\color{blue}{3}}$};
		\draw ([yshift=1mm]b.west)--node[above]{$\beta_{1}$}([yshift=1mm]a.east);
		\draw  ([yshift=-1mm]b.west)--node[below]{$\alpha_{1}$}([yshift=-1mm]a.east);
		\draw ([yshift=1mm]c.west)--node[above]{$\alpha_{2}$}([yshift=1mm]b.east);
		\draw  ([yshift=-1mm]c.west)--node[below]{$\beta_{2}$}([yshift=-1mm]b.east);
		\draw[dashed,-] ([xshift=.4cm]b.north) arc[start angle = 0, end angle = 180, x radius=.4cm, y radius =.2cm];
		\draw[dashed,-] ([xshift=-.4cm]b.south) arc[start angle = 180, end angle = 360, x radius=.4cm, y radius =.2cm];
		\end {tikzpicture} $};
		\node at (5,-1){$\begin{tikzpicture}[>= angle 60, ->,red,line width=0.5mm]
		\node (a) at (0,0) {$2$};
		\node (b) at (1,-1) {$3$};
		\node (c) at (2,0) {$2$};
		\node (d) at (3,-1) {$3$};
		\node (e) at (4,0) {$2$};
		\draw (a) -- node[above right]{$\beta_{2}$} (b);
		\draw (c) -- node[below right]{$\alpha_{2}$}  (b);
		\draw (c) -- node[above right]{$\beta_{2}$}  (d);
		\draw (e) -- node[below right]{$\alpha_{2}$}  (d);
		\end{tikzpicture}$};
	\end{tikzpicture}
	\caption{\label{fig:accordandrep} Illustration of the one-to-one correspondence between accordions $\Surf(Q,R)$ and strings of $(Q,R)$ via an example.}
	\end{figure}
\end{ex}

	From now on, we denote by $\Accord(\pmb{\Sigma}, \mathcal{M},\Delta^{\gpoint}) = \Accord$ the set of accordions on a given $\gpoint$-dissected marked surface $(\pmb{\Sigma}, \mathcal{M},\Delta^{\gpoint})$, by $\Accord'$ the subset of $\Accord$ composed of accordions associated with non-projective objects and by $\Accord_{[e]}$ the subset of $\Accord$ composed of accordions admitting $e\in\mathcal{M}_{\rpoint}$ as an endpoint. For any $\delta \in \Accord$, by abuse of notation, we set $\MM(\delta)$ as its associated string representation of $(\opQ(\Delta^{\gpoint}), \opR (\Delta^{\gpoint}))$. Conversely, for any $M \in \ind(Q,R)$, we write $\gamma_{(M)}$ for its associated accordion on $\Surf(Q,R)$.

\begin{definition}
\label{def:counterclockorder}
     Let $e \in \mathcal{M}_{\rpoint}$. The \new{counterclockwise order around e} is denoted by $\preccurlyeq_e$ and defined  for any $(\delta_1,\delta_2) \in \Accord_{[e]}^2$ by  $\delta_1 \preccurlyeq_e \delta_2$ if $\delta_2$ if the oriented angle formed by $\delta_1$ and $\delta_2$ at vertex $e$, in this order, is counterclockwise.
\end{definition}

\subsection{Morphisms and extensions}
\label{ss:GeomKerandExt}
We recall some results allowing us to read morphisms and extensions between indecomposable representations of $(Q,R)$ directly via the associated accordions on $\Surf(Q,R)$.
First, we recall the result that describes the homomorphism space between the string modules associated with the accordions. The following proposition is a direct consequence of \cref{thm:CB} applied to $\Accord$.

\begin{prop}[\cite{BR87}] \label{prop:HomCrossing} Let $(\pmb{\Sigma},\mathcal{M},\Delta^{\gpoint})$ be a $\gpoint$-dissected marked surface with $\mathcal{M}_{\gpoint} \subset \partial{\pmb{\Sigma}}$. For any pair $(\delta, \eta) \in \Accord$, a basis of $\Hom(\MM(\delta), \MM(\eta))$ is given by crossings as depicted in \cref{fig:HomCrossing}.
\begin{figure}[ht!]
    \centering
    \begin{tikzpicture}[mydot/.style={
				circle,
				thick,
				fill=white,
				draw,
				outer sep=0.5pt,
				inner sep=1pt
			}, fl/.style={->,>=latex}]
			\tkzDefPoint(0,0){O}\tkzDefPoint(0,1.5){1} 
			\tkzDefPoint(-1.5,1.5){2} 
			\tkzDefPoint(-1.5,0){3} 
			
			\tkzDefPoint(3,0){4}\tkzDefPoint(3,1.5){5} 
			\tkzDefPoint(4.5,1.5){6} 
			\tkzDefPoint(4.5,0){7} 
			
			\filldraw [fill=gray,opacity=0.1] (O) to (1) to [bend right=10] (5) to (4) to [bend right=10] cycle;
			
			\draw[line width=0.5mm,dark-green](1) edge (O);
			\draw[line width=0.5mm,dark-green, bend left=30](1) edge (2);
			\draw[line width=0.5mm,dark-green,bend right=30](O) edge (3);
			
			\draw[line width=0.5mm,dark-green](4) edge (5);
			\draw[line width=0.5mm,dark-green, bend left=30](6) edge (5);
			\draw[line width=0.5mm,dark-green,bend right=30](7) edge (4);
			
			\draw[line width=0.7mm,orange,densely dotted](-0.75,1.5) to [bend right=40] (0,0.85) to [bend left=0] node[above]{$\delta$}  (3,0.85) to [bend left=40] (3.75,0);
			\draw[line width=0.7mm,blue, bend left=30, loosely dotted](-0.75,0) to [bend left=30]  (0,0.65) to [bend left=0] node[below]{$\eta$} (3,0.65) to [bend right=30] (3.75,1.5);
			
			\tkzDrawPoints[size=4,color=dark-green,mydot](O,1,2,3,4,5,6,7);
		\end{tikzpicture}
    \caption{\label{fig:HomCrossing} Illustration of a crossing corresponding to a basis element of $\Hom(\MM(\delta), \MM(\eta))$.  The segments of $\delta$ and $\eta$ with $\pmb{\Sigma}\setminus\pmb{\Gamma}(\Delta^{\gpoint})$ are homotopic.}
\end{figure}
\end{prop}

\begin{remark} \label{rem:orientationandmorph}
Note that, given two accordions that are crossing, the orientation of the surface induces the directions of arrows and thus the direction of the morphism between the two associated indecomposable representations.
\end{remark}

We recall a result that allows us to read extensions between indecomposable representations in the geometric model.

\begin{prop}\label{prop:geom_ext}
Let $(\pmb{\Sigma}, \mathcal{M},\Delta^{\gpoint})$ be a $\gpoint$-dissected marked surface with $\mathcal{M}_{\gpoint} \subset \partial \pmb{\Sigma}$. Let $\delta, \eta \in \Accord$. We distinguish two types of extensions of $\MM(\delta)$ by $\MM(\eta)$ which are:

\begin{enumerate}[label = $\bullet$, itemsep =0.5em]
    \item \new{Overlap extensions}: whenever we have a non-split short exact sequence \[\begin{tikzcd}
	 \MM(\eta) & {E_1 \oplus E_2} & \MM(\delta) 
	\arrow[tail, from=1-1, to=1-2]
	\arrow[two heads, from=1-2, to=1-3]
\end{tikzcd},\] where $E_1,E_2 \in \ind_\mathbb{K}(Q,R)$, then $\delta$ and $\eta$ intersect, and $\gamma_{(E_1)}$ and $\gamma_{(E_2)}$ can be obtained from $\delta$ and $\eta$ as pictured in \cref{fig:overlapextaccord}. By abuse of notation, we write $\OvExt(\delta,\eta)$ for the union of all sets $\{\gamma_{(E_1)}, \gamma_{(E_2)}\}$ of curves, up to homotopy, over all the isomorphism classes of short exact sequences of the above shape.
    \begin{figure}[!ht]
\centering 
    \begin{tikzpicture}[mydot/.style={
					circle,
					thick,
					fill=white,
					draw,
					outer sep=0.5pt,
					inner sep=1pt
				}, scale = 1]
		\tikzset{
		osq/.style={
        rectangle,
        thick,
        fill=white,
        append after command={
            node [
                fit=(\tikzlastnode),
                orange,
                line width=0.3mm,
                inner sep=-\pgflinewidth,
                cross out,
                draw
            ] {}}}}
		\draw [line width=0.7mm,domain=50:130] plot ({4*cos(\x)}, {1.5*sin(\x)});
        \draw [line width=0.7mm,domain=230:310] plot ({4*cos(\x)}, {1.5*sin(\x)});
		\foreach \X in {0}
        {
		\tkzDefPoint(4*cos(pi/6*\X +pi/3),1.5*sin(pi/6*\X + pi/3)){\X};
		};
        \foreach \X in {1}
		{
		\tkzDefPoint(4*cos(pi/6*\X +pi/2),1.5*sin(pi/6*\X + pi/2)){\X};
		};
        \foreach \X in {2}
        {
		\tkzDefPoint(4*cos(pi/6*(\X-2.1) +3*pi/2.2),1.5*sin(pi/6*(\X-2.1) + 3*pi/2.2)){\X};
		}
		\foreach \X in {3}
		{
		\tkzDefPoint(4*cos(pi/6*(\X-2.1) +3*pi/2),1.5*sin(pi/6*(\X-2.1) + 3*pi/2)){\X};
		};
        
        \tkzDefPoint(4*cos(pi/6*0 +pi/3.5),1.5*sin(pi/6*0 + pi/3.5)){e0};
        \tkzDefPoint(4*cos(pi/6*0 +pi/2.6),1.5*sin(pi/6*0 + pi/2.6)){e0'};
         \tkzDefPoint(4*cos(pi/6*1 +pi/1.85),1.5*sin(pi/6*1 + pi/1.85)){e1};
        \tkzDefPoint(4*cos(pi/6*1 +pi/2.2),1.5*sin(pi/6*1 + pi/2.2)){e1'};
        \tkzDefPoint(4*cos(pi/6*(2-2.1) +3*pi/2.28),1.5*sin(pi/6*(2-2.1) + 3*pi/2.28)){e2};
        \tkzDefPoint(4*cos(pi/6*(2-2.1) +3*pi/2.15),1.5*sin(pi/6*(2-2.1) + 3*pi/2.15)){e2'};
        \tkzDefPoint(4*cos(pi/6*(3-2.1) +3*pi/1.94),1.5*sin(pi/6*(3-2.1) + 3*pi/1.94)){e3};
        \tkzDefPoint(4*cos(pi/6*(3-2.1) +3*pi/2.05),1.5*sin(pi/6*(3-2.1) + 3*pi/2.05)){e3'};

        \filldraw [fill=gray,opacity=0.1] (e1') to (e2') to [bend left=14] (e3') to (e0') to [bend left=15] cycle;

        \draw[line width=0.5mm,dark-green](e0') edge (e3');
        \draw[line width=0.5mm,dark-green](e1') edge (e2');

        \draw[line width=0.5mm,dark-green, bend left=30](e0) edge (e0');
        \draw[line width=0.5mm,dark-green, bend left=30](e1') edge (e1);
        \draw[line width=0.5mm,dark-green, bend left=30](e2) edge (e2');
        \draw[line width=0.5mm,dark-green, bend left=30](e3') edge (e3);

		\draw[line width=0.7mm,bend right=5,orange, dotted](0) edge (2);
		
		\draw[line width=0.7mm,bend right=10,blue, loosely dotted](1) edge (3);
		
		\draw [line width=0.7mm, mypurple,dash pattern={on 10pt off 2pt on 5pt off 2pt}, bend right=20] (1) edge (0);
		\draw [line width=0.7mm, mypurple,dash pattern={on 10pt off 2pt on 5pt off 2pt}, bend right=20] (3) edge (2);
		\foreach \X in {0,...,3}
		{
		\tkzDrawPoints[circle,fill =red,size=4,color=red](\X);
		};

        \foreach \X in {e0,e0',e1,e1',e2,e2',e3,e3'}
		{
		\tkzDrawPoints[mydot,size=6,color=dark-green,thick,fill=white](\X);
		};
		
		\begin{scope}[xshift=3ex]
		\tkzDefPoint(-0.1,0.3){gammaM};
		\tkzLabelPoint[orange](gammaM){\Large $\delta$}
		\tkzDefPoint(-1.6,0){gammaP};
		\tkzLabelPoint[blue](gammaP){\Large $\eta$}
		\tkzDefPoint(-1.2,1.0){deltav};
		\tkzLabelPoint[mypurple](deltav){\Large $\gamma_{(E_1)}$}
		\tkzDefPoint(-0.6,-0.4){deltaw};
		\tkzLabelPoint[mypurple](deltaw){\Large $\gamma_{(E_2)}$}
		\end{scope}
    \end{tikzpicture}
    
\caption{\label{fig:overlapextaccord} Illustration of an overlap extension. The crossing is resolved by selecting these two arcs as the four need to share the common part of $\delta$ and $\eta$. The green lines correspond to the dissection of the surface.}
\end{figure}
    \item \new{Arrow extension}: 
    whenever we have a nonsplit short exact sequence  \[\begin{tikzcd}
	 \MM(\eta) & {E} & \MM(\delta) 
	\arrow[tail, from=1-1, to=1-2]
	\arrow[two heads, from=1-2, to=1-3]
\end{tikzcd},\] where $E \in \ind(Q,R)$, then we can construct $\gamma_{(E)}$ from $\delta$ and $\eta$ as described in \cref{fig:arrowextaccord}. Similarly to overlap extensions, by abuse of notation, we write $\ArExt(\delta,\eta)$ for the set of all the accordions $\gamma_{(E)}$ over all the isomorphism classes of short exact sequences of the above shape.
 
\begin{figure}[!ht]
\centering 
    \begin{tikzpicture}[mydot/.style={
					circle,
					thick,
					fill=white,
					draw,
					outer sep=0.5pt,
					inner sep=1pt
				}, scale = 1]
		\tikzset{
		osq/.style={
        rectangle,
        thick,
        fill=white,
        append after command={
            node [
                fit=(\tikzlastnode),
                orange,
                line width=0.3mm,
                inner sep=-\pgflinewidth,
                cross out,
                draw
            ] {}}}}
            
        \draw [line width=0.7mm,domain=-10:10] plot ({4*cos(\x)}, {1.5*sin(\x)});
		\draw [line width=0.7mm,domain=20:50] plot ({4*cos(\x)}, {1.5*sin(\x)});
		\draw [line width=0.7mm,domain=80:100] plot ({4*cos(\x)}, {1.5*sin(\x)});
		\draw [line width=0.7mm,domain=110:130] plot ({4*cos(\x)}, {1.5*sin(\x)});
		\draw [line width=0.7mm,domain=170:190] plot ({4*cos(\x)}, {1.5*sin(\x)});
        \draw [line width=0.7mm,domain=250:290] plot ({4*cos(\x)}, {1.5*sin(\x)});
		\foreach \X in {0,1}
		{
		\tkzDefPoint(4*cos(pi/2*(\X-0.3) +pi/3),1.5*sin(pi/2*(\X-0.3) + pi/3)){\X};
		};
		\foreach \X in {2,3}
		{
		\tkzDefPoint(4*cos(pi/6*(\X-3) +3*pi/2),1.5*sin(pi/6*(\X-3) + 3*pi/2)){\X};
		};
		
		\tkzDefPoint(4,0){4};
		\tkzDefPoint(0,1.5){5};
		\tkzDefPoint(-4,0){6};
		\tkzDefPoint(4*cos(3*pi/2-pi/11),1.5*sin(3*pi/2-pi/11)){7};
		\tkzDefPoint(4*cos(3*pi/2+pi/11),1.5*sin(3*pi/2+pi/11)){8};
		
		\filldraw[gray,opacity=0.1] (4) to [bend left=20] (5) to [bend left=20] (6) to [bend left=10] (7) to [bend right=10] (8) to [bend left=10] cycle;
		
		\draw[line width=0.5mm,bend left = 20,dark-green] (4) to (5) to (6);
		\draw[line width=0.5mm,bend left = 10,dark-green,dashed] (6) to (7);
		\draw[line width=0.5mm,bend left = 10,dark-green,dashed] (8) to (4);
		
		\draw[line width=0.7mm,bend right=10,orange, dotted](0) edge (3);
		
		\draw[line width=0.7mm,bend left=10,blue, loosely dotted](1) edge (3);
		
		\draw [line width=0.7mm, mypurple,dash pattern={on 10pt off 2pt on 5pt off 2pt}, bend left=20] (0) edge (1);

		\foreach \X in {0,1,3}
		{
		\tkzDrawPoints[circle,fill =red,size=4,color=red](\X);
		};
		
		\foreach \X in {4,...,8}
		{
		\tkzDrawPoints[mydot,size=6,color=dark-green,thick,fill=white](\X);
		};
        
		\begin{scope}[xshift=2ex]
		\tkzDefPoint(0.6,0.1){gammaM};
		\tkzLabelPoint[orange](gammaM){\Large $\delta$}
		\tkzDefPoint(-1.6,0.3){gammaP};
		\tkzLabelPoint[blue](gammaP){\Large $\eta$}
		\tkzDefPoint(-0.4,1.1){deltav};
		\tkzLabelPoint[mypurple](deltav){\Large $\gamma_{(E)}$}
		\end{scope}
    \end{tikzpicture}
\caption{\label{fig:arrowextaccord} Illustration of an arrow extension. The gray part corresponds to one cell of $\pmb{\Gamma}(\Delta^{\gpoint})$.}
\end{figure}
\end{enumerate}
\end{prop}
\begin{proof}
This is a direct transcription of \cref{thm:ExplicitExtGentle} in the geometric model.
\end{proof}

For our purposes, it will be useful to give a geometric interpretation of the kernel of any epimorphism$\begin{tikzcd}
	M & Y
	\arrow[two heads, from=1-1, to=1-2]
\end{tikzcd}$involving $M \in \rep(Q,R)$ and $Y \in \ind(Q,R)$.
	
\begin{prop}\label{prop:geo_mor}
Let $(\pmb{\Sigma},\mathcal{M},\Delta^{\gpoint})$ be a $\gpoint$-dissected marked surface with $\mathcal{M}_{\gpoint} \subset \partial{\pmb{\Sigma}}$. Consider $p \in \mathbb{N}^*$ and $\delta_1,\ldots, \delta_p \in \Accord$. Let $\eta \in \Accord$ such that there exists a minimal epimorphism \[ \begin{tikzcd}
	{\displaystyle f: \bigoplus_{i=1}^p \MM(\delta_i)} & \MM(\eta)
	\arrow[two heads, from=1-1, to=1-2]
\end{tikzcd}\]
such that $f_{|\MM(\delta_i)}$ is given by exactly one crossing between $\delta_i$ and $\eta$. Then $\Ker(f) = \bigoplus_{i=0}^p \MM(\kappa_i)$ where $\kappa_0,\ldots,\kappa_p \in \Accord$  are constructed from $\delta_1, \ldots, \delta_p$, and $\eta$ as in \cref{fig:kerepi}. Note that we define $\kappa_0$ and $\kappa_p$ dually, and thus the drawing of $\kappa_p$ is dual to that of $\kappa_0$.
\begin{figure}[!ht]
\centering 
     \begin{tikzpicture}[mydot/.style={
					circle,
					thick,
					fill=white,
					draw,
					outer sep=0.5pt,
					inner sep=1pt
				}, scale = 1]
		\tikzset{
		osq/.style={
        rectangle,
        thick,
        fill=white,
        append after command={
            node [
                fit=(\tikzlastnode),
                orange,
                line width=0.3mm,
                inner sep=-\pgflinewidth,
                cross out,
                draw
            ] {}}}}
		\draw [line width=0.7mm,domain=50:130] plot ({4*cos(\x)}, {1.5*sin(\x)});
        \draw [line width=0.7mm,domain=230:310] plot ({4*cos(\x)}, {1.5*sin(\x)});
		\foreach \X in {0,1}
		{
		\tkzDefPoint(4*cos(pi/6*\X +pi/2),1.5*sin(pi/6*\X + pi/2)){\X};
		};
		\foreach \X in {2,3}
		{
		\tkzDefPoint(4*cos(pi/6*(\X-2) +3*pi/2),1.5*sin(pi/6*(\X-2) + 3*pi/2)){\X};
		};
		
		\draw[line width=0.7mm,bend right=20,blue, loosely dotted](-3,0.5) edge (0,0);
		
		\draw[line width=0.7mm,bend left=20,blue, loosely dotted](0,0) edge (3,-0.5);
		
		\draw[line width=0.7mm,bend right=10,orange,dotted](0) edge (2);

		\draw [line width=0.7mm, mypurple,dash pattern={on 10pt off 2pt on 5pt off 2pt}, bend right=20] (1) edge (2);
		\draw [line width=0.7mm, mypurple,dash pattern={on 10pt off 2pt on 5pt off 2pt}, bend right=20] (0) edge (3);
		
		\foreach \X in {0,1,2,3}
		{
		\tkzDrawPoints[fill =red,size=4,color=red](\X);
		};

		\begin{scope}[xshift=2ex,yshift=.6ex]
			\tkzDefPoint(-2.5,0){gammaM};
			\tkzLabelPoint[blue](gammaM){\Large $\eta$}
			\tkzDefPoint(0,-0.7){gammaP};
			\tkzLabelPoint[orange](gammaP){\Large $\kappa_i$}
			\tkzDefPoint(-2.4,0.9){deltav};
			\tkzLabelPoint[mypurple](deltav){\Large $\delta_i$}
			\tkzDefPoint(0.45,0.9){deltaw};
			\tkzLabelPoint[mypurple](deltaw){\Large $\delta_{i+1}$}
		\end{scope}

		\begin{scope}[xshift = 6.5cm]
		\draw [line width=0.7mm,domain=50:130] plot ({4*cos(\x)}, {1.5*sin(\x)});
        \draw [line width=0.7mm,domain=230:310] plot ({4*cos(\x)}, {1.5*sin(\x)});
		\foreach \X in {0,1}
		{
		\tkzDefPoint(4*cos(pi/6*\X +pi/2),1.5*sin(pi/6*\X + pi/2)){\X};
		};
		\foreach \X in {2,3}
		{
		\tkzDefPoint(4*cos(pi/6*(\X-2) +3*pi/2),1.5*sin(pi/6*(\X-2) + 3*pi/2)){\X};
		};
		
		\draw[line width=0.7mm,bend right=20,blue, loosely dotted](-3,-0.5) edge (0);
		
		\draw [line width=0.7mm, mypurple,dash pattern={on 10pt off 2pt on 5pt off 2pt}, bend right=20] (1) edge (3);
		
		\draw[line width=0.7mm,bend right=10,orange,dotted](0) edge (3);
		
		\foreach \X in {0,1,3}
		{
		\tkzDrawPoints[fill =red,size=4,color=red](\X);
		};

		\begin{scope}[xshift=2ex,yshift=.6ex]
		\tkzDefPoint(-2.8,0.4){gammaM};
		\tkzLabelPoint[blue](gammaM){\Large $\eta$}
		\tkzDefPoint(0.7,0.7){gammaP};
		\tkzLabelPoint[orange](gammaP){\Large $\kappa_p$}
		\tkzDefPoint(-0.3,0.1){deltav};
		\tkzLabelPoint[mypurple](deltav){\Large $\delta_p$}
		\end{scope}
		\end{scope}
    \end{tikzpicture}
\caption{\label{fig:kerepi} The two types of accordions $\kappa_i$ representing the indecomposable summands of $\Ker(f)$.}
\end{figure}
\end{prop}
\begin{proof} 
This is a direct transcription of \cref{thm:StringKerEpi} in the geometric model where $\kappa_i=\gamma(\MM(\nu_i))$.
\end{proof}

\subsection{Geometric model for gentle trees}

Previous results are simplified when the quiver is a gentle tree. This subsection will collect these.

\begin{theorem}[\cite{OPS18,APS19}]\label{thm:surftree}
Let $(Q,R)$ be a gentle tree. Set $(\pmb{\Sigma}, \mathcal{M}, \Delta^{\gpoint}) = \Surf(Q,R)$. Then $\mathcal{M} \subset \partial \pmb{\Sigma}$ and $\pmb{\Sigma}$ is homeomorphic to a disc.
\end{theorem}

\begin{remark} \label{rem:translatevdimtree} The previous theorem and \cref{prop:vdimtree} imply that, for a given gentle tree $(Q,R)$:
\begin{enumerate}[label=$\bullet$,itemsep=1mm]
    \item  strings of $(Q,R)$ are characterized by the vertices of $Q$ appearing as beginning or ends of arrows composing it; and,
    \item accordions of $\Surf(Q,R)$ are characterized by the $\gpoint$-arcs of the $\gpoint$-dissection they cross.
\end{enumerate}
\end{remark}

\begin{remark}\label{rem:translatehometree} We can infer from \cref{prop:homtree} that, given a gentle tree $(Q,R)$, the accordions of $\Surf(Q,R)$ cross at most once, up to homotopy.
\end{remark}

\subsection{Projective perspective}
\label{ss:ProjVariant}
As we will use minimal projective resolutions substantially in what follows, we first translate the condition $\gldim(Q,R)<\infty$ into the geometric model.

\begin{prop}[\cite{OPS18,LGH24}] \label{prop:gldimsurf}
Let $(Q,R)$ be a gentle quiver. Consider $\Surf(Q,R) = (\pmb{\Sigma}, \mathcal{M}, \Delta^{\gpoint})$ its associated $\gpoint$-dissected marked surface. Then we have that $\gldim(Q,R) < \infty$ if and only if $\mathcal{M} \subset \partial \pmb{\Sigma}$.
\end{prop}

\begin{remark}\label{rem:proofgldimsurf}
In particular, the fact that $\mathcal{M}_{\rpoint} \subset \partial \pmb{\Sigma}$ ensures that we avoid the pattern highlighted in \cref{prop:Gentlegldimfinite}.
\end{remark}

Under these conditions, we introduce a new modification of the standard geometric model based on the indecomposable projective objects of $\rep(Q,R)$.
	
\begin{definition} \label{def:projdissection} Let $(\pmb{\Sigma}, \mathcal{M}, \Delta^{\gpoint})$ be a $\gpoint$-dissected marked surface. We define the \new{projective dissection} of $(\pmb{\Sigma}, \mathcal{M}, \Delta^{\gpoint})$ as the set $\Prj(\Delta^{\gpoint})$ of accordions corresponding to the indecomposable representations of $\proj(\opQ(\Delta^{\gpoint}),\opR(\Delta^{\gpoint}))$.
\end{definition}

\begin{prop} \label{prop:projdissec}
For any $\gpoint$-dissected marked surface $(\pmb{\Sigma}, \mathcal{M}, \Delta^{\gpoint})$, the triplet $(\pmb{\Sigma},\mathcal{M}, \Prj(\Delta^{\gpoint}))$ gives a $\rpoint$-dissected marked surface. Moreover, if $\mathcal{M} \subset \partial \pmb{\Sigma}$, then there exists an oriented homeomorphism from $(\pmb{\Sigma}, \mathcal{M}, \Delta^{\gpoint})$ to $(\pmb{\Sigma},\mathcal{M}, \Prj(\Delta^{\gpoint}))$.
\end{prop}

\begin{proof}
    By \cite[Theorem 3.2.]{APS19}, we have that $\Prj(\Delta^{\gpoint})$ is a $\rpoint$-dissection of $(\pmb{\Sigma},\mathcal{M})$. It follows from \cref{thm:GeomandRep} that the accordions associated with projective objects are obtained by rotating the endpoints of the original $\gpoint$-dissection around the boundary components in the orientation induced by the counterclockwise orientation of the surface.
\end{proof}

See \cref{fig:projdissec} for examples of $\Prj(\Delta^{\gpoint})$ on our running example.
	\begin{figure}[!ht]
		\centering
		\begin{tikzpicture}[mydot/.style={
					circle,
					thick,
					fill=white,
					draw,
					outer sep=0.5pt,
					inner sep=1pt
				}, scale = 1]
			
			\tkzDefPoint(0,0){O}\tkzDefPoint(0,1.5){1}
				\tkzDefPointsBy[rotation=center O angle 180](1){2}
				\tkzDefPoint(0.7,0){3} \tkzDefPoint(1.1,0){8} \tkzDefPoint(-1.1,0){9}
				\tkzDefPointsBy[rotation=center O angle 90](3,4,5){4,5,6}
				\tkzDrawCircle[line width=0.7mm,black](O,1)
				\tkzDrawCircle[line width=0.7mm,black,fill=gray!40](O,3)
				\tkzDrawPoints[fill =red,size=4,color=red](2,3,5)

				\draw[line width=0.5mm,dark-green](1) edge (4);
				\draw[line width=0.5mm,bend left = 60,dark-green](4) edge (8);
				\draw[line width=0.5mm,bend left = 60,dark-green](8) edge (6);
				\draw[line width=0.5mm,bend right = 60,dark-green](4) edge (9);
				\draw[line width=0.5mm,bend right =60,dark-green](9) edge (6);
				\tkzDrawPoints[size=4,color=dark-green,mydot](1,4,6)
				
				\draw[line width=0.9mm,red,bend right=60] (5) edge (0,-1.1);
				\draw[line width=0.9mm,red,bend right=60] (0,-1.1) edge (3);
				\draw[line width=0.9mm,red,bend right=60] (3) edge (0,1.1);
				\draw[line width=0.9mm,red,bend right=60] (0,1.1) edge (5);
				\draw[line width=0.9mm,red,bend left=40] (2) edge (-1.35,0);
				\draw[line width=0.9mm,red,bend left=40] (-1.35,0) edge (0,1.25);
				\draw[line width=0.9mm,bend left=30,red](0,1.25) edge (0.9,0.7);
				\draw[line width=0.9mm,bend left=35,red](0.9,0.7) edge (3);
			\begin{scope}[xshift=-.5cm]	
			\begin{scope}[xshift=5cm,yshift=1cm, scale=0.7]	\tkzDefPoint(0,0){O}\tkzDefPoint(1.4,1.4){1}
				\tkzDefPointsBy[rotation=center O angle 90](1,2,3){2,3,4}
				\tkzDefPoint(0.9,0.9){5}
				\tkzDefPointsBy[rotation=center O angle 90](5,6,7){6,7,8}
				
				\draw[line width=0.3mm] (-1.4142,-1.4142) rectangle (1.4142,1.4142);
				\tkzDefPointsBy[rotation=center 1 angle 45](5){9}
				\tkzDrawSector[rotate,line width=0.5mm, black, fill = gray!40](1,9)(-90)
				\tkzDefPointsBy[rotation=center 2 angle 45](6){10}
				\tkzDrawSector[rotate,line width=0.5mm, black, fill = gray!40](2,10)(-90)
				\tkzDefPointsBy[rotation=center 3 angle 45](7){11}
				\tkzDrawSector[rotate,line width=0.5mm, black, fill = gray!40](3,11)(-90)
				\tkzDefPointsBy[rotation=center 4 angle 45](8){12}
				\tkzDrawSector[rotate,line width=0.5mm, black, fill = gray!40](4,12)(-90)
				
				\tkzDefPoint(1.4,0){13}
				\tkzDefPoint(-1.4,0){14}
				\tkzDefPoint(0,1.4){15}
				\tkzDefPoint(0,-1.4){16}
				
				\tkzDrawPoints[fill =red,size=4,color=red](5,7)
				\draw[line width=0.5mm,dark-green](6) edge (8);
				\draw[line width=0.5mm,bend left=40,dark-green](6) edge (14);
				\draw[line width=0.5mm,bend right=40,dark-green](6) edge (15);
				\draw[line width=0.5mm,bend left=40,dark-green](8) edge (13);
				\draw[line width=0.5mm,bend right=40,dark-green](8) edge (16);
				\tkzDrawPoints[size=4,color=dark-green,mydot](6,8)
				
				\draw[line width=0.9mm,red,bend right=40](5) edge (1.4,0);
				\draw[line width=0.9mm,bend left=40,red](-1.4,0) edge (7);
				
				\draw[line width=0.9mm,red,bend right=40](5) edge (1.4,0.3);
				\draw[line width=0.9mm,red](1.4,-0.3) edge (-1.4,0.3);
				\draw[line width=0.9mm,red,bend left=40] (-1.4,-0.3) edge (7);
				
				\draw[line width=0.9mm,red,bend right=40](5) edge (1.4,0.5);
				\draw[line width=0.9mm,red,bend right=40] (-1.4,0.5) edge (0,1.4);
				\draw[line width=0.9mm,red,bend left=40] (0,-1.4) edge (1.4,-0.5);
				\draw[line width=0.9mm,red,bend left=40] (-1.4,-0.5) edge (7);
			\end{scope}
			\begin{scope}[xshift=7cm,yshift=1cm, scale=0.7]	\tkzDefPoint(0,0){O}\tkzDefPoint(1.4,1.4){1}
				\tkzDefPointsBy[rotation=center O angle 90](1,2,3){2,3,4}
				\tkzDefPoint(0.9,0.9){5}
				\tkzDefPointsBy[rotation=center O angle 90](5,6,7){6,7,8}
				
				\draw[line width=0.3mm] (-1.4142,-1.4142) rectangle (1.4142,1.4142);
				\tkzDefPointsBy[rotation=center 1 angle 45](5){9}
				\tkzDrawSector[rotate,line width=0.5mm, black, fill = gray!40](1,9)(-90)
				\tkzDefPointsBy[rotation=center 2 angle 45](6){10}
				\tkzDrawSector[rotate,line width=0.5mm, black, fill = gray!40](2,10)(-90)
				\tkzDefPointsBy[rotation=center 3 angle 45](7){11}
				\tkzDrawSector[rotate,line width=0.5mm, black, fill = gray!40](3,11)(-90)
				\tkzDefPointsBy[rotation=center 4 angle 45](8){12}
				\tkzDrawSector[rotate,line width=0.5mm, black, fill = gray!40](4,12)(-90)
				
				\tkzDefPoint(1.4,0){13}
				\tkzDefPoint(-1.4,0){14}
				\tkzDefPoint(0,1.4){15}
				\tkzDefPoint(0,-1.4){16}
				
				\tkzDrawPoints[fill =red,size=4,color=red](5,7)
				\draw[line width=0.5mm,dark-green](6) edge (8);
				\draw[line width=0.5mm,bend left=40,dark-green](6) edge (14);
				\draw[line width=0.5mm,bend right=40,dark-green](6) edge (15);
				\draw[line width=0.5mm,bend left=40,dark-green](8) edge (13);
				\draw[line width=0.5mm,bend right=40,dark-green](8) edge (16);
				\tkzDrawPoints[size=4,color=dark-green,mydot](6,8)
				
				\draw[line width=0.9mm,red,bend right=40](5) edge (1.4,0);
				\draw[line width=0.9mm,bend left=40,red](-1.4,0) edge (7);
				
				\draw[line width=0.9mm,red,bend right=40](5) edge (1.4,0.3);
				\draw[line width=0.9mm,red](1.4,-0.3) edge (-1.4,0.3);
				\draw[line width=0.9mm,red,bend left=40] (-1.4,-0.3) edge (7);
				
				\draw[line width=0.9mm,red,bend right=40](5) edge (1.4,0.5);
				\draw[line width=0.9mm,red,bend right=40] (-1.4,0.5) edge (0,1.4);
				\draw[line width=0.9mm,red,bend left=40,dashdotted] (0,-1.4) edge (1.4,-0.5);
				\draw[line width=0.9mm,red,bend left=40] (-1.4,-0.5) edge (7);
			\end{scope}
			\begin{scope}[xshift=5cm,yshift=-1cm, scale=0.7]	\tkzDefPoint(0,0){O}\tkzDefPoint(1.4,1.4){1}
				\tkzDefPointsBy[rotation=center O angle 90](1,2,3){2,3,4}
				\tkzDefPoint(0.9,0.9){5}
				\tkzDefPointsBy[rotation=center O angle 90](5,6,7){6,7,8}
				
				\draw[line width=0.3mm] (-1.4142,-1.4142) rectangle (1.4142,1.4142);
				\tkzDefPointsBy[rotation=center 1 angle 45](5){9}
				\tkzDrawSector[rotate,line width=0.5mm, black, fill = gray!40](1,9)(-90)
				\tkzDefPointsBy[rotation=center 2 angle 45](6){10}
				\tkzDrawSector[rotate,line width=0.5mm, black, fill = gray!40](2,10)(-90)
				\tkzDefPointsBy[rotation=center 3 angle 45](7){11}
				\tkzDrawSector[rotate,line width=0.5mm, black, fill = gray!40](3,11)(-90)
				\tkzDefPointsBy[rotation=center 4 angle 45](8){12}
				\tkzDrawSector[rotate,line width=0.5mm, black, fill = gray!40](4,12)(-90)
				
				\tkzDefPoint(1.4,0){13}
				\tkzDefPoint(-1.4,0){14}
				\tkzDefPoint(0,1.4){15}
				\tkzDefPoint(0,-1.4){16}
				
				\tkzDrawPoints[fill =red,size=4,color=red](5,7)
				\draw[line width=0.5mm,dark-green](6) edge (8);
				\draw[line width=0.5mm,bend left=40,dark-green](6) edge (14);
				\draw[line width=0.5mm,bend right=40,dark-green](6) edge (15);
				\draw[line width=0.5mm,bend left=40,dark-green](8) edge (13);
				\draw[line width=0.5mm,bend right=40,dark-green](8) edge (16);
				\tkzDrawPoints[size=4,color=dark-green,mydot](6,8)
				
				\draw[line width=0.9mm,red,bend right=40,dashdotted](5) edge (1.4,0);
				\draw[line width=0.9mm,bend left=40,red](-1.4,0) edge (7);
				
				\draw[line width=0.9mm,red,bend right=40,dashdotted](5) edge (1.4,0.3);
				\draw[line width=0.9mm,red](1.4,-0.3) edge (-1.4,0.3);
				\draw[line width=0.9mm,red,bend left=40] (-1.4,-0.3) edge (7);
				
				\draw[line width=0.9mm,red,bend right=40,dashdotted](5) edge (1.4,0.5);
				\draw[line width=0.9mm,red,bend right=40] (-1.4,0.5) edge (0,1.4);
				\draw[line width=0.9mm,red,bend left=40] (0,-1.4) edge (1.4,-0.5);
				\draw[line width=0.9mm,red,bend left=40] (-1.4,-0.5) edge (7);
			\end{scope}
			\begin{scope}[xshift=7cm,yshift=-1cm, scale=0.7]	\tkzDefPoint(0,0){O}\tkzDefPoint(1.4,1.4){1}
				\tkzDefPointsBy[rotation=center O angle 90](1,2,3){2,3,4}
				\tkzDefPoint(0.9,0.9){5}
				\tkzDefPointsBy[rotation=center O angle 90](5,6,7){6,7,8}
				
				\draw[line width=0.3mm] (-1.4142,-1.4142) rectangle (1.4142,1.4142);
				\tkzDefPointsBy[rotation=center 1 angle 45](5){9}
				\tkzDrawSector[rotate,line width=0.5mm, black, fill = gray!40](1,9)(-90)
				\tkzDefPointsBy[rotation=center 2 angle 45](6){10}
				\tkzDrawSector[rotate,line width=0.5mm, black, fill = gray!40](2,10)(-90)
				\tkzDefPointsBy[rotation=center 3 angle 45](7){11}
				\tkzDrawSector[rotate,line width=0.5mm, black, fill = gray!40](3,11)(-90)
				\tkzDefPointsBy[rotation=center 4 angle 45](8){12}
				\tkzDrawSector[rotate,line width=0.5mm, black, fill = gray!40](4,12)(-90)
				
				\tkzDefPoint(1.4,0){13}
				\tkzDefPoint(-1.4,0){14}
				\tkzDefPoint(0,1.4){15}
				\tkzDefPoint(0,-1.4){16}
				
				\tkzDrawPoints[fill =red,size=4,color=red](5,7)
				\draw[line width=0.5mm,dark-green](6) edge (8);
				\draw[line width=0.5mm,bend left=40,dark-green](6) edge (14);
				\draw[line width=0.5mm,bend right=40,dark-green](6) edge (15);
				\draw[line width=0.5mm,bend left=40,dark-green](8) edge (13);
				\draw[line width=0.5mm,bend right=40,dark-green](8) edge (16);
				\tkzDrawPoints[size=4,color=dark-green,mydot](6,8)
				
				\draw[line width=0.9mm,red,bend right=40](5) edge (1.4,0);
				\draw[line width=0.9mm,bend left=40,red,dashdotted](-1.4,0) edge (7);
				
				\draw[line width=0.9mm,red,bend right=40](5) edge (1.4,0.3);
				\draw[line width=0.9mm,red,dashdotted](1.4,-0.3) edge (-1.4,0.3);
				\draw[line width=0.9mm,red,bend left=40] (-1.4,-0.3) edge (7);
				
				\draw[line width=0.9mm,red,bend right=40](5) edge (1.4,0.5);
				\draw[line width=0.9mm,red,bend right=40,dashdotted] (-1.4,0.5) edge (0,1.4);
				\draw[line width=0.9mm,red,bend left=40] (0,-1.4) edge (1.4,-0.5);
				\draw[line width=0.9mm,red,bend left=40] (-1.4,-0.5) edge (7);
			\end{scope}
			\begin{scope}[xshift=9cm,yshift=1cm, scale=0.7]	\tkzDefPoint(0,0){O}\tkzDefPoint(1.4,1.4){1}
				\tkzDefPointsBy[rotation=center O angle 90](1,2,3){2,3,4}
				\tkzDefPoint(0.9,0.9){5}
				\tkzDefPointsBy[rotation=center O angle 90](5,6,7){6,7,8}
				
				\draw[line width=0.3mm] (-1.4142,-1.4142) rectangle (1.4142,1.4142);
				\tkzDefPointsBy[rotation=center 1 angle 45](5){9}
				\tkzDrawSector[rotate,line width=0.5mm, black, fill = gray!40](1,9)(-90)
				\tkzDefPointsBy[rotation=center 2 angle 45](6){10}
				\tkzDrawSector[rotate,line width=0.5mm, black, fill = gray!40](2,10)(-90)
				\tkzDefPointsBy[rotation=center 3 angle 45](7){11}
				\tkzDrawSector[rotate,line width=0.5mm, black, fill = gray!40](3,11)(-90)
				\tkzDefPointsBy[rotation=center 4 angle 45](8){12}
				\tkzDrawSector[rotate,line width=0.5mm, black, fill = gray!40](4,12)(-90)
				
				\tkzDefPoint(1.4,0){13}
				\tkzDefPoint(-1.4,0){14}
				\tkzDefPoint(0,1.4){15}
				\tkzDefPoint(0,-1.4){16}
				
				\tkzDrawPoints[fill =red,size=4,color=red](5,7)
				\draw[line width=0.5mm,dark-green](6) edge (8);
				\draw[line width=0.5mm,bend left=40,dark-green](6) edge (14);
				\draw[line width=0.5mm,bend right=40,dark-green](6) edge (15);
				\draw[line width=0.5mm,bend left=40,dark-green](8) edge (13);
				\draw[line width=0.5mm,bend right=40,dark-green](8) edge (16);
				\tkzDrawPoints[size=4,color=dark-green,mydot](6,8)
				
				\draw[line width=0.9mm,red,bend right=40](5) edge (1.4,0);
				\draw[line width=0.9mm,bend left=40,red](-1.4,0) edge (7);
				
				\draw[line width=0.9mm,red,bend right=40](5) edge (1.4,0.3);
				\draw[line width=0.9mm,red](1.4,-0.3) edge (-1.4,0.3);
				\draw[line width=0.9mm,red,bend left=40] (-1.4,-0.3) edge (7);
				
				\draw[line width=0.9mm,red,bend right=40](5) edge (1.4,0.5);
				\draw[line width=0.9mm,red,bend right=40] (-1.4,0.5) edge (0,1.4);
				\draw[line width=0.9mm,red,bend left=40] (0,-1.4) edge (1.4,-0.5);
				\draw[line width=0.9mm,red,bend left=40,dashdotted] (-1.4,-0.5) edge (7);
			\end{scope}
			\begin{scope}[xshift=9cm,yshift=-1cm, scale=0.7]	\tkzDefPoint(0,0){O}\tkzDefPoint(1.4,1.4){1}
				\tkzDefPointsBy[rotation=center O angle 90](1,2,3){2,3,4}
				\tkzDefPoint(0.9,0.9){5}
				\tkzDefPointsBy[rotation=center O angle 90](5,6,7){6,7,8}
				
				\draw[line width=0.3mm] (-1.4142,-1.4142) rectangle (1.4142,1.4142);
				\tkzDefPointsBy[rotation=center 1 angle 45](5){9}
				\tkzDrawSector[rotate,line width=0.5mm, black, fill = gray!40](1,9)(-90)
				\tkzDefPointsBy[rotation=center 2 angle 45](6){10}
				\tkzDrawSector[rotate,line width=0.5mm, black, fill = gray!40](2,10)(-90)
				\tkzDefPointsBy[rotation=center 3 angle 45](7){11}
				\tkzDrawSector[rotate,line width=0.5mm, black, fill = gray!40](3,11)(-90)
				\tkzDefPointsBy[rotation=center 4 angle 45](8){12}
				\tkzDrawSector[rotate,line width=0.5mm, black, fill = gray!40](4,12)(-90)
				
				\tkzDefPoint(1.4,0){13}
				\tkzDefPoint(-1.4,0){14}
				\tkzDefPoint(0,1.4){15}
				\tkzDefPoint(0,-1.4){16}
				
				\tkzDrawPoints[fill =red,size=4,color=red](5,7)
				\draw[line width=0.5mm,dark-green](6) edge (8);
				\draw[line width=0.5mm,bend left=40,dark-green](6) edge (14);
				\draw[line width=0.5mm,bend right=40,dark-green](6) edge (15);
				\draw[line width=0.5mm,bend left=40,dark-green](8) edge (13);
				\draw[line width=0.5mm,bend right=40,dark-green](8) edge (16);
				\tkzDrawPoints[size=4,color=dark-green,mydot](6,8)
				
				\draw[line width=0.9mm,red,bend right=40](5) edge (1.4,0);
				\draw[line width=0.9mm,bend left=40,red](-1.4,0) edge (7);
				
				\draw[line width=0.9mm,red,bend right=40](5) edge (1.4,0.3);
				\draw[line width=0.9mm,red](1.4,-0.3) edge (-1.4,0.3);
				\draw[line width=0.9mm,red,bend left=40,dashdotted] (-1.4,-0.3) edge (7);
				
				\draw[line width=0.9mm,red,bend right=40](5) edge (1.4,0.5);
				\draw[line width=0.9mm,red,bend right=40] (-1.4,0.5) edge (0,1.4);
				\draw[line width=0.9mm,red,bend left=40] (0,-1.4) edge (1.4,-0.5);
				\draw[line width=0.9mm,red,bend left=40] (-1.4,-0.5) edge (7);
			\end{scope}
            \end{scope}
			\end{tikzpicture}
			\caption{\label{fig:projdissec} Examples of $\Prj(\Delta^{\gpoint})$ of the \mbox{$\gpoint$-dissections} seen in \cref{exam::three-surfaces}. For the torus example (right), we draw a part of its universal cover, and the dashed red lines correspond to the fundamental copy of all the projective accordions. Note that every red arc is obtained by moving every endpoint of a green arc to the closest red point in
the chosen orientation.}
		\end{figure}

In the follow-up sections, we will consider $\rpoint$-dissected marked surfaces by seeing them as the projective dissection of its associated $\gpoint$-dissected one. The following propositions are useful for considering accordions in the context of projective dissections.  

 Given $C \in \pmb{\Gamma}(\Prj(\Delta^{\gpoint}))$, a cell of the projective dissection. 
There exists a unique edge $e$ in $\partial C\cap\partial\pmb{\Sigma}$; denote by $m_C$ the $\rpoint$-point which is the counterclockwise endpoint of $e$ following the orientation induced on $\partial\pmb{\Sigma}$.

\begin{prop} \label{prop:arcsaccordions} Let $(\pmb{\Sigma}, \mathcal{M}, \Delta^{\gpoint})$ be a $\gpoint$-dissected marked surface such that $\mathcal{M} \subset \partial \pmb{\Sigma}$. A $\rpoint$-arc $\delta$ is an accordion if and only if either $\delta \in \Prj(\Delta^{\gpoint})$ or the following assertions hold: 
\begin{enumerate}[label=$(\alph*)$, itemsep=1mm]
    \item whenever $\delta$ enters a cell $C$ of $\Prj(\Delta^{\gpoint})$ by crossing $\mu \in \Prj(\Delta^{\gpoint})$,
    \begin{enumerate}[label=$(a \arabic*)$, itemsep=1mm]
        \item if $\delta$ leaves $C$, $\delta$ leaves it by crossing  a projective accordion adjacent to $\mu$;
        \item otherwise, its endpoint in $\partial C$ is either $m_C$ or the non-common endpoint of a projective accordion adjacent to $\mu$, or 
    \end{enumerate}
    \item if $\delta$ is contained in a cell $C$, then:
    \begin{enumerate}[label=$(b \arabic*)$, itemsep=1mm]
        \item its endpoints are the distinct endpoints of a pair of adjacent projective accordions; or,
        \item one of its endpoints must be $m_C$. 
    \end{enumerate}
\end{enumerate}
\end{prop}

\begin{remark} \label{rem:translationfromgreen}
    \cref{prop:arcsaccordions} is a translation of the definition of accordion in the $\gpoint$-dissection $\Delta^{\gpoint}$ into $\Prj(\Delta^{\gpoint})$.
\end{remark}

\begin{figure}[!ht]
			\centering
			 \begin{tikzpicture}[mydot/.style={
					circle,
					thick,
					fill=white,
					draw,
					outer sep=0.5pt,
					inner sep=1pt
				}, scale = 1]
			
			\tkzDefPoint(0,0){O}\tkzDefPoint(0,1.5){1}
				\tkzDefPointsBy[rotation=center O angle 180](1){2}
				\tkzDefPoint(0.7,0){3} \tkzDefPoint(1.1,0){8} \tkzDefPoint(-1.1,0){9}
				\tkzDefPointsBy[rotation=center O angle 90](3,4,5){4,5,6}
				\tkzDrawCircle[line width=0.7mm,black](O,1)
				\tkzDrawCircle[line width=0.7mm,black,fill=gray!40](O,3)
				\tkzDrawPoints[fill =red,size=4,color=red](2,3,5)

				\draw[line width=0.5mm,dark-green](1) edge (4);
				\draw[line width=0.5mm,bend left = 60,dark-green](4) edge (8);
				\draw[line width=0.5mm,bend left = 60,dark-green](8) edge (6);
				\draw[line width=0.5mm,bend right = 60,dark-green](4) edge (9);
				\draw[line width=0.5mm,bend right =60,dark-green](9) edge (6);
				\tkzDrawPoints[size=4,color=dark-green,mydot](1,4,6)

                \draw[line width=0.9mm,red,bend right=60] (5) edge (0,-1.1);
				\draw[line width=0.9mm,red,bend right=60] (0,-1.1) edge (3);
				\draw[line width=0.9mm,red,bend right=60] (3) edge (0,1.1);
				\draw[line width=0.9mm,red,bend right=60] (0,1.1) edge (5);
				\draw[line width=0.9mm,red,bend left=40] (2) edge (-1.35,0);
				\draw[line width=0.9mm,red,bend left=40] (-1.35,0) edge (0,1.25);
				\draw[line width=0.9mm,bend left=30,red](0,1.25) edge (0.9,0.7);
				\draw[line width=0.9mm,bend left=35,red](0.9,0.7) edge (3);

                \draw[line width=0.9mm,blue,bend right=40,loosely dotted] (2) edge (1.35,0);
				\draw[line width=0.9mm,blue,bend right=40, loosely dotted] (1.35,0) edge (0,1.25);
				\draw[line width=0.9mm,bend right=60,blue, loosely dotted](0,1.25) edge (5);
			\end{tikzpicture}
			\caption{\label{fig:accordionPrj} Example of the accordion drawn in \cref{fig:accordandclosedaccord} (left), represented in blue dotted line, in the projective dissection $\Prj(\Delta^{\gpoint})$ in \cref{fig:projdissec} (left). Note that the $\rpoint$-point $m_C$ for the cell $C$ whose boundaries contain the three projective accordions is the unique point on the outer boundary.}
		\end{figure}

In the remainder of this section, we introduce the notion of \emph{neighboring projective accordions} of a given accordion $\gamma$ in $(\pmb{\Sigma}, \mathcal{M}, \Delta^{\gpoint})$, which will prove to be useful in the next sections.
	
\begin{definition}\label{def:ExtprojCover}
Let $(Q,R)$ be a representation-finite gentle quiver. Let $\delta$ be an accordion of $\Surf(Q,R)$. A \new{neighboring projective accordion} of $\delta$ is an accordion $\gamma_{(P)}$ in $\Surf(Q,R)$ associated to an indecomposable projective representation $P \in \proj_\mathbb{K}(Q,R)$ such that either:
\begin{enumerate}[label=$\bullet$,itemsep=0.1em]
    \item $P$ and $\MM(\delta)$ have non trivial higher extension;
    \item $P$ is a direct summand of some projective representation appearing in the minimal projective resolution $\left(\PP_{\MM(\delta)}^\ast \right)$.
\end{enumerate}
For any accordion $\gamma$ of $\Surf(Q,R)$, write $\NP(\delta)$ for the set of all the neighboring projective accordions of $\delta$.
\end{definition}

We give an explicit computation of $\NP(\delta)$ in \cref{ex:Nproj1}, showing the generic behavior of the projective accordions that must appear in it. The following proposition describes the projective accordions in the geometric model corresponding to projective representaions in $\NP(\delta)$

\begin{prop}\label{prop:CombidescripNproj}
Let $(Q,R)$ be a representation-finite gentle quiver. Let $\delta$ be an accordion of $\Surf(Q,R)$. Then $\NP(\delta)$ is the subset of $\Prj(\Delta^{\gpoint})$ made of accordions $\eta$ satisfying at least one of the following conditions:
\begin{enumerate}[label=$(\roman*)$,itemsep=0.2em]
    \item the curve $\eta$ crosses $\delta$;
    \item the curve $\eta$ is included in the border of a cell of the projective dissection crossed by $\delta$ containing one of its endpoints, and is a part of the path of projective curves that links the last projective crossed to the endpoint of $\delta$.
    \item the curve $\eta$ and $\delta$ have a common endpoint $v$, and $\delta\preccurlyeq_v\eta$
\end{enumerate}
\end{prop}
\begin{proof}
    It follows directly from the description of the items of \cref{def:ExtprojCover}  done in \cite{CPS21} \cite{OPS18} and \cite{LGH24}.
\end{proof}
\begin{remark} \label{rem:condiaccordNproj}
The projective accordions $\eta$ from condition $(i)$ are the ones such that $\OvExt(\delta, \eta)\neq\emptyset$, those from condition $(ii)$ are some of the projectives coming from $\left( \PP_{\MM(\delta)}^\ast \right)$, and those from condition $(iii)$ are the ones having non-trivial higher extensions. 
\end{remark}

\begin{ex} \label{ex:Nproj1} In \cref{fig:Nprojex1}, we represent a projective dissection of a $\gpoint$-dissection of a marked disc (drawn as thick red curves), and we consider an accordion $\delta$ (drawn as a dotted blue curve). The densely dotted thick red lines correspond to accordions forming the set of its neighboring projective accordions $\NP(\delta)$. We pictured projective accordions appearing in  $(\PP_{\MM(\delta)}^\ast)$ in  \cref{fig:explicitcalculationNproj1}. \qedhere
\begin{figure}[!ht]
\centering 
    \begin{tikzpicture}[mydot/.style={
					circle,
					thick,
					fill=white,
					draw,
					outer sep=0.5pt,
					inner sep=1pt
				}, scale = .7]
		\draw[line width=0.7mm,black] (0,0) ellipse (4cm and 1.5cm);
		\foreach \X in {0,1,...,37}
		{
		\tkzDefPoint(4*cos(pi/19*\X),1.5*sin(pi/19*\X)){\X};
		};
        
		\draw[line width=0.9mm,bend left =60,red](1) edge (5);
		\draw[line width=0.9mm,bend left =60,red](3) edge (5);
		\draw[line width=0.9mm,bend right =30,red](5) edge (35);
		\draw[line width=0.9mm,bend left =30,red](33) edge (35);
		\draw[line width=0.9mm,bend left =30,red](35) edge (37);
		\draw[line width=0.9mm,bend left =30,red](31) edge (33);
		\draw[line width=0.9mm,bend left =30,red](7) edge (31);
		\draw[line width=0.9mm,bend left =30,red](7) edge (29);
		\draw[line width=0.9mm,bend left =30,red](9) edge (29);
		\draw[line width=0.9mm,bend left =30,red](11) edge (29);
		\draw[line width=0.9mm,bend left =30,red](11) edge (27);
		\draw[line width=0.9mm,bend right =30,red](13) edge (11);
		\draw[line width=0.9mm,bend left =30,red](13) edge (15);
		\draw[line width=0.9mm,bend left =30,red](17) edge (23);
		\draw[line width=0.9mm,bend left =30,red](15) edge (17);
		\draw[line width=0.9mm,bend left =30,red](19) edge (23);
		\draw[line width=0.9mm,bend left =30,red](21) edge (23);
		\draw[line width=0.9mm,bend left =30,red](25) edge (27);

		\draw[line width=0.7mm,bend right=10,blue, loosely dashed](5) edge (23);
		\foreach \X in {0,2,...,36}
		{
		\tkzDrawPoint[size=4,color=dark-green,mydot](\X);
		};
		\foreach \X in {1,3,...,37}
		{
		\tkzDrawPoints[fill =red,size=4,color=red](\X);
		};
		\begin{scope}[xshift=2ex,yshift=.6ex]
		\tkzDefPoint(-1.6,0.8){gamma};
		\tkzLabelPoint[blue](gamma){\Large $\delta$}
		\end{scope}

        \begin{scope}[xshift=9cm]
            \draw[line width=0.7mm,black] (0,0) ellipse (4cm and 1.5cm);
		\foreach \X in {0,1,...,37}
		{
		\tkzDefPoint(4*cos(pi/19*\X),1.5*sin(pi/19*\X)){\X};
		};
        
		\draw[line width=0.9mm,bend left =60,red, densely dashdotted](1) edge (5);
		\draw[line width=0.9mm,bend left =60,red, densely dashdotted](3) edge (5);
		\draw[line width=0.9mm,bend right =30,red, densely dashdotted](5) edge (35);
		\draw[line width=0.9mm,bend left =30,red,densely dashdotted](33) edge (35);
		\draw[line width=0.9mm,bend left =30,red, densely dashdotted](31) edge (33);
		\draw[line width=0.9mm,bend left =30,red, densely dashdotted](7) edge (31);
		\draw[line width=0.9mm,bend left =30,red, densely dashdotted](7) edge (29);
		\draw[line width=0.9mm,bend left =30,red, densely dashdotted](9) edge (29);
		\draw[line width=0.9mm,bend left =30,red, densely dashdotted](11) edge (29);
		\draw[line width=0.9mm,bend left =30,red, densely dashdotted](11) edge (27);
		\draw[line width=0.9mm,bend right =30,red, densely dashdotted](13) edge (11);
		\draw[line width=0.9mm,bend left =30,red, densely dashdotted](13) edge (15);
		\draw[line width=0.9mm,bend left =30,red, densely dashdotted](17) edge (23);
		\draw[line width=0.9mm,bend left =30,red, densely dashdotted](15) edge (17);
		\draw[line width=0.9mm,bend left =30,red, densely dashdotted](19) edge (23);
		\draw[line width=0.9mm,bend left =30,red,densely dashdotted](21) edge (23);

		\draw[line width=0.7mm,bend right=10,blue, loosely dashed](5) edge (23);
		\foreach \X in {0,2,...,36}
		{
		\tkzDrawPoint[size=4,color=dark-green,mydot](\X);
		};
		\foreach \X in {1,3,...,37}
		{
		\tkzDrawPoints[fill =red,size=4,color=red](\X);
		};
		\begin{scope}[xshift=2ex,yshift=.6ex]
		\tkzDefPoint(-1.6,0.8){gamma};
		\tkzLabelPoint[blue](gamma){\Large $\delta$}
		\end{scope}
        \end{scope}
    \end{tikzpicture}
\caption{\label{fig:Nprojex1} An example of $\NP(\delta)$.}
\end{figure}
\end{ex}

\begin{figure}[p]
    \centering
    \[\begin{tikzcd}
	{\begin{tikzpicture}[mydot/.style={
					circle,
					thick,
					fill=white,
					draw,
					outer sep=0.5pt,
					inner sep=1pt
				}, scale = 0.65]
		\draw[line width=0.7mm,black] (0,0) ellipse (4cm and 1.5cm);
		\foreach \X in {0,1,...,37}
		{
		\tkzDefPoint(4*cos(pi/19*\X),1.5*sin(pi/19*\X)){\X};
		};
		
		\draw[line width=0.9mm,bend left =60,red](1) edge (5);
		\draw[line width=0.9mm,bend left =60,red](3) edge (5);
		\draw[line width=0.9mm,bend right =30,red,densely dashdotted](5) edge (35);
		\draw[line width=0.9mm,bend left =30,red](35) edge (37);
		\draw[line width=0.9mm,bend left =30,red,densely dashdotted](33) edge (35);
		\draw[line width=0.9mm,bend left =30,red,densely dashdotted](31) edge (33);
		\draw[line width=0.9mm,bend left =30,red, densely dashdotted](7) edge (31);
		\draw[line width=0.9mm,bend left =30,red,densely dashdotted](7) edge (29);
		\draw[line width=0.9mm,bend left =30,red](9) edge (29);
		\draw[line width=0.9mm,bend left =30,red,densely dashdotted](11) edge (29);
		\draw[line width=0.9mm,bend left =30,red](11) edge (27);
		\draw[line width=0.9mm,bend right =30,red,densely dashdotted](13) edge (11);
		\draw[line width=0.9mm,bend left =30,red,densely dashdotted](13) edge (15);
		\draw[line width=1.5mm,bend left=30,red,opacity=0.3](17) edge (23);
		\draw[line width=0.9mm,bend left =30,red,densely dashdotted](17) edge (23);
		\draw[line width=0.9mm,bend left =30,red,densely dashdotted](15) edge (17);
		\draw[line width=0.9mm,bend left =30,red](19) edge (23);
		\draw[line width=0.9mm,bend left =30,red](21) edge (23);
		\draw[line width=0.9mm,bend left =30,red](25) edge (27);
		
		\draw[line width=0.7mm,bend right=10,blue, loosely dashed](5) edge (23);
	
		\foreach \X in {1,3,...,37}
		{
		\tkzDrawPoints[fill =red,size=4,color=red](\X);
		};
        \node[blue] at (-1.5,.2){$\delta$};
        \node[red] at (-2.7,0){$\rho_{10}$};
    \end{tikzpicture}} & {P_{10}} \\
	{\begin{tikzpicture}[mydot/.style={
					circle,
					thick,
					fill=white,
					draw,
					outer sep=0.5pt,
					inner sep=1pt
				}, scale = 0.65]
		\draw[line width=0.7mm,black] (0,0) ellipse (4cm and 1.5cm);
		\foreach \X in {0,1,...,37}
		{
		\tkzDefPoint(4*cos(pi/19*\X),1.5*sin(pi/19*\X)){\X};
		};
		
		\draw[line width=0.9mm,bend left =60,red](1) edge (5);
		\draw[line width=0.9mm,bend left =60,red](3) edge (5);
		\draw[line width=1.3mm,bend right=30,red,opacity=0.3](5) edge (35);
		\draw[line width=0.9mm,bend right =30,red,densely dashdotted](5) edge (35);
		\draw[line width=0.9mm,bend left =30,red](35) edge (37);
		\draw[line width=0.9mm,bend left =30,red,densely dashdotted](33) edge (35);
		\draw[line width=0.9mm,bend left =30,red,densely dashdotted](31) edge (33);
		\draw[line width=0.9mm,bend left =30,red, densely dashdotted](7) edge (31);
		\draw[line width=0.9mm,bend left =30,red,densely dashdotted](7) edge (29);
		\draw[line width=0.9mm,bend left =30,red](9) edge (29);
		\draw[line width=0.9mm,bend left =30,red,densely dashdotted](11) edge (29);
		\draw[line width=0.9mm,bend left =30,red](11) edge (27);
		\draw[line width=0.9mm,bend right =30,red,densely dashdotted](13) edge (11);
		\draw[line width=0.9mm,bend left =30,red,densely dashdotted](13) edge (15);
		\draw[line width=0.9mm,bend left =30,red](17) edge (23);
		\draw[line width=1.3mm,bend left=30,red,opacity=0.3](15) edge (17);
		\draw[line width=0.9mm,bend left =30,red,densely dashdotted](15) edge (17);
		\draw[line width=0.9mm,bend left =30,red](19) edge (23);
		\draw[line width=0.9mm,bend left =30,red](21) edge (23);
		\draw[line width=0.9mm,bend left =30,red](25) edge (27);
		
		\draw[line width=0.7mm,bend right=10,blue, loosely dashed](5) edge (23);
		
		\foreach \X in {1,3,...,37}
		{
		\tkzDrawPoints[fill =red,size=4,color=red](\X);
		};
        \node[blue] at (-1.5,.3){$\delta$};
        \node[red] at (-3.1,.3){$\rho_8$};
        \node[red] at (2.5,-.1){$\rho_9$};
    \end{tikzpicture}} & {P_8 \oplus P_9}\\
{\begin{tikzpicture}[mydot/.style={
					circle,
					thick,
					fill=white,
					draw,
					outer sep=0.5pt,
					inner sep=1pt
				}, scale = 0.65]
		\draw[line width=0.7mm,black] (0,0) ellipse (4cm and 1.5cm);
		\foreach \X in {0,1,...,37}
		{
		\tkzDefPoint(4*cos(pi/19*\X),1.5*sin(pi/19*\X)){\X};
		};
		
		\draw[line width=0.9mm,bend left =60,red](1) edge (5);
		\draw[line width=0.9mm,bend left =60,red](3) edge (5);
		\draw[line width=0.9mm,bend right =30,red](5) edge (35);
		\draw[line width=0.9mm,bend left =30,red](35) edge (37);
		\draw[line width=1.5mm,bend left=30,red,opacity=0.3](33) edge (35);
		\draw[line width=0.9mm,bend left =30,red,densely dashdotted](33) edge (35);
		\draw[line width=0.9mm,bend left =30,red,densely dashdotted](31) edge (33);
		\draw[line width=0.9mm,bend left =30,red, densely dashdotted](7) edge (31);
		\draw[line width=0.9mm,bend left =30,red,densely dashdotted](7) edge (29);
		\draw[line width=0.9mm,bend left =30,red](9) edge (29);
		\draw[line width=0.9mm,bend left =30,red,densely dashdotted](11) edge (29);
		\draw[line width=0.9mm,bend left =30,red](11) edge (27);
		\draw[line width=0.9mm,bend right =30,red,densely dashdotted](13) edge (11);
		\draw[line width=1.5mm,bend left=30,red,opacity=0.3](13) edge (15);
		\draw[line width=0.9mm,bend left =30,red,densely dashdotted](13) edge (15);
		\draw[line width=0.9mm,bend left =30,red](17) edge (23);
		\draw[line width=0.9mm,bend left =30,red](15) edge (17);
		\draw[line width=0.9mm,bend left =30,red](19) edge (23);
		\draw[line width=0.9mm,bend left =30,red](21) edge (23);
		\draw[line width=0.9mm,bend left =30,red](25) edge (27);
		
		\draw[line width=0.7mm,bend right=10,blue, loosely dashed](5) edge (23);

		\foreach \X in {1,3,...,37}
		{
		\tkzDrawPoints[fill =red,size=4,color=red](\X);
		};
        \node[blue] at (-1.5,.2){$\delta$};
        \node[red] at (-2.5,.6){$\rho_6$};
        \node[red] at (2.7,-.6){$\rho_7$};
    \end{tikzpicture}} & {P_6 \oplus P_7} \\
	{\begin{tikzpicture}[mydot/.style={
					circle,
					thick,
					fill=white,
					draw,
					outer sep=0.5pt,
					inner sep=1pt
				}, scale = 0.65]
		\draw[line width=0.7mm,black] (0,0) ellipse (4cm and 1.5cm);
		\foreach \X in {0,1,...,37}
		{
		\tkzDefPoint(4*cos(pi/19*\X),1.5*sin(pi/19*\X)){\X};
		};
	
		\draw[line width=0.9mm,bend left =60,red](1) edge (5);
		\draw[line width=0.9mm,bend left =60,red](3) edge (5);
		\draw[line width=0.9mm,bend right =30,red](5) edge (35);
		\draw[line width=0.9mm,bend left =30,red](35) edge (37);
		\draw[line width=0.9mm,bend left =30,red](33) edge (35);
		\draw[line width=1.5mm,bend left=30,red,opacity=0.3](31) edge (33);
		\draw[line width=0.9mm,bend left =30,red,densely dashdotted](31) edge (33);
		\draw[line width=0.9mm,bend left =30,red, densely dashdotted](7) edge (31);
		\draw[line width=1.5mm,bend left=30,red,opacity=0.3](7) edge (29);
		\draw[line width=0.9mm,bend left =30,red,densely dashdotted](7) edge (29);
		\draw[line width=0.9mm,bend left =30,red](9) edge (29);
		\draw[line width=0.9mm,bend left =30,red,densely dashdotted](11) edge (29);
		\draw[line width=0.9mm,bend left =30,red](11) edge (27);
		\draw[line width=1.5mm,bend left=30,red,opacity=0.3](11) edge (13);
		\draw[line width=0.9mm,bend right =30,red,densely dashdotted](13) edge (11);
		\draw[line width=0.9mm,bend left =30,red](13) edge (15);
		\draw[line width=0.9mm,bend left =30,red](17) edge (23);
		\draw[line width=0.9mm,bend left =30,red](15) edge (17);
		\draw[line width=0.9mm,bend left =30,red](19) edge (23);
		\draw[line width=0.9mm,bend left =30,red](21) edge (23);
		\draw[line width=0.9mm,bend left =30,red](25) edge (27);
		
		\draw[line width=0.7mm,bend right=10,blue, loosely dashed](5) edge (23);
	
		\foreach \X in {1,3,...,37}
		{
		\tkzDrawPoints[fill =red,size=4,color=red](\X);
		};
        \node[blue] at (-1.5,.2){$\delta$};
        \node[red] at (-1.7,.8){$\rho_3$};
        \node[red] at (1.2,.4){$\rho_4$};
        \node[red] at (2.5,-.7){$\rho_5$};
    \end{tikzpicture}} & {P_3 \oplus P_4 \oplus P_5}\\
	{\begin{tikzpicture}[mydot/.style={
					circle,
					thick,
					fill=white,
					draw,
					outer sep=0.5pt,
					inner sep=1pt
				}, scale = 0.65]
		\draw[line width=0.7mm,black] (0,0) ellipse (4cm and 1.5cm);
		\foreach \X in {0,1,...,37}
		{
		\tkzDefPoint(4*cos(pi/19*\X),1.5*sin(pi/19*\X)){\X};
		};
		
		\draw[line width=0.9mm,bend left =60,red](1) edge (5);
		\draw[line width=0.9mm,bend left =60,red](3) edge (5);
		\draw[line width=0.9mm,bend right =30,red](5) edge (35);
		\draw[line width=0.9mm,bend left =30,red](35) edge (37);
		\draw[line width=0.9mm,bend left =30,red](33) edge (35);
		\draw[line width=0.9mm,bend left =30,red](31) edge (33);
		\draw[line width=1.5mm,bend left=30,red,opacity=0.3](7) edge (31);
		\draw[line width=0.9mm,bend left =30,red, densely dashdotted](7) edge (31);
		\draw[line width=0.9mm,bend left =30,red](7) edge (29);
		\draw[line width=0.9mm,bend left =30,red](9) edge (29);
		\draw[line width=1.5mm,bend left=30,red,opacity=0.3](11) edge (29);
		\draw[line width=0.9mm,bend left =30,red,densely dashdotted](11) edge (29);
		\draw[line width=0.9mm,bend left =30,red](11) edge (27);
		\draw[line width=0.9mm,bend right =30,red](13) edge (11);
		\draw[line width=0.9mm,bend left =30,red](13) edge (15);
		\draw[line width=0.9mm,bend left =30,red](17) edge (23);
		\draw[line width=0.9mm,bend left =30,red](15) edge (17);
		\draw[line width=0.9mm,bend left =30,red](19) edge (23);
		\draw[line width=0.9mm,bend left =30,red](21) edge (23);
		\draw[line width=0.9mm,bend left =30,red](25) edge (27);
		
		\draw[line width=0.7mm,bend right=10,blue, loosely dashed](5) edge (23);
	
		\foreach \X in {1,3,...,37}
		{
		\tkzDrawPoints[fill =red,size=4,color=red](\X);
		};

        \node[blue] at (-1.5,.2){$\delta$};
        \node[red] at (0,-.6){$\rho_1$};
        \node[red] at (1.5,-.8){$\rho_2$};
    \end{tikzpicture}} & {P_{1} \oplus P_2}\\ 
	{\begin{tikzpicture}[mydot/.style={
					circle,
					thick,
					fill=white,
					draw,
					outer sep=0.5pt,
					inner sep=1pt
				}, scale = 0.65]
		\draw[line width=0.7mm,black] (0,0) ellipse (4cm and 1.5cm);
		\foreach \X in {0,1,...,37}
		{
		\tkzDefPoint(4*cos(pi/19*\X),1.5*sin(pi/19*\X)){\X};
		};
        
		\draw[line width=0.9mm,bend left =60,red](1) edge (5);
		\draw[line width=0.9mm,bend left =60,red](3) edge (5);
		\draw[line width=0.9mm,bend right =30,red](5) edge (35);
		\draw[line width=0.9mm,bend left =30,red](35) edge (37);
		\draw[line width=0.9mm,bend left =30,red](33) edge (35);
		\draw[line width=0.9mm,bend left =30,red](31) edge (33);
		\draw[line width=0.9mm,bend left =30,red](7) edge (31);
		\draw[line width=0.9mm,bend left =30,red](7) edge (29);
		\draw[line width=0.9mm,bend left =30,red](9) edge (29);
		\draw[line width=0.9mm,bend left =30,red](11) edge (29);
		\draw[line width=0.9mm,bend left =30,red](11) edge (27);
		\draw[line width=0.9mm,bend right =30,red](13) edge (11);
		\draw[line width=0.9mm,bend left =30,red](13) edge (15);
		\draw[line width=0.9mm,bend left =30,red](17) edge (23);
		\draw[line width=0.9mm,bend left =30,red](15) edge (17);
		\draw[line width=0.9mm,bend left =30,red](19) edge (23);
		\draw[line width=0.9mm,bend left =30,red](21) edge (23);
		\draw[line width=0.9mm,bend left =30,red](25) edge (27);
		
		\draw[line width=1mm,bend right=10,blue,opacity=0.2](5) edge (23);
		\draw[line width=0.7mm,bend right=10,blue, loosely dashed](5) edge (23);
		
		\foreach \X in {1,3,...,37}
		{
		\tkzDrawPoints[fill =red,size=4,color=red](\X);
		};

        \node[blue] at (-1.5,.2){$\delta$};
    \end{tikzpicture}} & {\MM(\delta)}
    \arrow[from=1-1, to=2-1]
	\arrow[from=2-1, to=3-1]
	\arrow[from=3-1, to=4-1]
	\arrow[from=4-1, to=5-1]
	\arrow[two heads, from=5-1, to=6-1]
    \arrow[from=1-2, to=2-2]
	\arrow[from=2-2, to=3-2]
	\arrow[from=3-2, to=4-2]
	\arrow[from=4-2, to=5-2]
	\arrow[two heads, from=5-2, to=6-2]
\end{tikzcd} \]
    \caption{\label{fig:explicitcalculationNproj1} Explicit calculation of the projective accordions $\rho_i \in \NP(\delta)$ of \cref{fig:Nprojex1} such that $\MM(\rho_i) = P_i$ is a summand of projective representations in $(\PP_{\MM(\delta)}^\ast)$.}
\end{figure}

\begin{cor}
\label{cor:geom_syzygy}
Let $(\pmb{\Sigma},\mathcal{M},\Delta^{\gpoint})$ be a $\gpoint$-dissected marked surface with $\mathcal{M}_{\gpoint} \subset \partial{\pmb{\Sigma}}$. Consider $p \in \mathbb{N}^*$ and $\delta_1,\ldots, \delta_p \in \Prj(\Delta^{\gpoint})$. Let $\eta \in \Accord$ satisfying that the following epimorphism is the projective cover of $M(\eta)$: \[ \begin{tikzcd}
	{\displaystyle f: \bigoplus_{i=1}^p \MM(\delta_i)} & \MM(\eta)
	\arrow[two heads, from=1-1, to=1-2]
\end{tikzcd}\]
Then the syzygy $\Ker(f)= \Omega(M(\eta)) = \bigoplus_{i=0}^p \MM(\kappa_i)$ where $\kappa_0,\ldots,\kappa_p \in \Accord$  are constructed from $\delta_1, \ldots, \delta_p$, and $\eta$ as shown in \cref{fig:kerepi}. The modules $M(\kappa_0)$ and $M(\kappa_p)$ are the only indecomposable summands of the syzygy that can be non-projective. Moreover, the accordions associated with the non-projective summands of the higher syzygies form the fan in the large cells of $\NP(\eta)$ as depicted in \cref{fig:syzygyl}. 
\begin{figure}[!ht]
\centering 
\begin{tikzpicture}[mydot/.style={
					circle,
					thick,
					fill=white,
					draw,
					outer sep=0.5pt,
					inner sep=1pt
				}, scale = 1]
		\tikzset{
		osq/.style={
        rectangle,
        thick,
        fill=white,
        append after command={
            node [
                fit=(\tikzlastnode),
                orange,
                line width=0.3mm,
                inner sep=-\pgflinewidth,
                cross out,
                draw
            ] {}}}}
		\draw [line width=0.7mm,domain=40:140] plot ({4*cos(\x)}, {2*sin(\x)});
        \draw [line width=0.7mm,domain=210:355] plot ({4*cos(\x)}, {2*sin(\x)});
		\foreach \X in {0,1}
		{
		\tkzDefPoint(4*cos(pi/3*\X +2*pi/5),2*sin(pi/3*\X + 2* pi/5)){\X};
		};
		\foreach \X in {2,3,4,5,6}
		{
		\tkzDefPoint(4*cos(pi/7*(\X-2) +4*pi/3),2*sin(pi/7*(\X-2) + 4*pi/3)){\X};
		};
		
		\tkzDefPoint(4*cos(pi/6 +2*pi/5),2*sin(pi/6 + 2* pi/5)){7};
		\tkzDefPoint(4*cos(-pi/7 +2*pi/5),2*sin(-pi/7 + 2* pi/5)){8};
		\tkzDefPoint(4*cos(-pi/14 +2*pi/5),2*sin(-pi/14 + 2* pi/5)){11};
		\tkzDefPoint(4*cos(-pi/8 +4*pi/3),2*sin(-pi/8 + 4*pi/3)){9}
		\tkzDefPoint(4*cos(-pi/12 +4*pi/3),2*sin(-pi/12 + 4*pi/3)){10}

		\draw[line width=0.7mm,bend right=10,mypurple,dash pattern={on 10pt off 2pt on 5pt off 2pt}](0) edge (2);
		
		\draw[line width=0.9mm,bend right=30,red](1) edge (2);
		\draw[line width=0.9mm,bend left=30,red](2) edge (3);
		\draw[line width=0.9mm,bend left=30,red](3) edge (4);
		\draw[line width=0.9mm,bend left=30,red](4) edge (5);
		\draw[line width=0.9mm,bend left=30,red](5) edge (6);
		\draw[line width=0.9mm,bend right=30,red](0) edge (6);
		\draw[line width=0.9mm,bend right=30,red](1) edge (7);
		\draw[line width=0.9mm,bend left=30,red](9) edge (2);
        
		\draw[line width=0.7mm, bend left=20, blue, loosely dashed] (9) edge (0);
		
		\draw [line width=0.7mm, mypurple,dash pattern={on 10pt off 2pt on 5pt off 2pt}, bend left=20] (3) edge (0);
		\draw [line width=0.7mm, mypurple,dash pattern={on 10pt off 2pt on 5pt off 2pt}, bend left=20] (4) edge (0);
		\draw [line width=0.7mm, mypurple,dash pattern={on 10pt off 2pt on 5pt off 2pt}, bend left=20] (5) edge (0);
		
		\filldraw [fill=red,opacity=0.1] (0) to [bend right=10] (7) to [bend left=30] (1) to [bend right=30] (2) to [bend left=30] (3) to [bend left=30] (4) to [bend left=30] (5) to [bend left=30] (6) to [bend left=30] cycle;

		\foreach \X in {0,...,7}
		{
		\tkzDrawPoints[fill =red,size=4,color=red](\X);
		};
		\tkzDrawPoints[fill =red,size=4,color=red](9);
		
		\begin{scope}[xshift=2ex,yshift=.6ex]
		\tkzDefPoint(-1.5,1.3){gammaM};
		\tkzLabelPoint[blue](gammaM){\Large $\eta$}
		\tkzDefPoint(-2.3,-0.8){omega1};
		\tkzLabelPoint[mypurple](omega1){\Large $\omega_1$}
		\tkzDefPoint(-.9,-1){omega2};
		\tkzLabelPoint[mypurple](omega2){\Large $\omega_2$}
		\tkzDefPoint(0.5,-1){omega3};
		\tkzLabelPoint[mypurple](omega3){\Large $\omega_3$}
	    \tkzDefPoint(1.8,-0.7){omega4};
		\tkzLabelPoint[mypurple](omega4){\Large $\omega_4$}
		
		\end{scope}
    \end{tikzpicture}
\caption{\label{fig:syzygyl}  The nonprojective summands $M(\omega_i)$ of the higher syzygy $\Omega^i(M(\eta))$ of $M(\eta)$.}
\end{figure}
\end{cor}
\begin{proof}
The first statement is a corollary of \cref{prop:geo_mor}. The second part is obtained by iteration of kernel computations of \cref{prop:geo_mor}. Accordions appearing in the projective cover of $M(\delta)$ are described in \cite{OPS18}. Summands of the projective replacement used in the computations of the higher syzygies appear in the border of large cells using \cref{prop:projdissec}.
\end{proof}
\begin{remark}
We recover the result \cite[Prop. 3.1]{BS21}, but here we express it in our geometric model.
\end{remark}

\begin{Notation}
In what follows, for any $\delta \in \Accord$, we denote by $\NP(\delta)_0$ the subset of $\mathcal{M}_{\rpoint}$ which elements are endpoints of accordions in $\NP(\delta)$. 
\end{Notation}

	\section{Monogeneous resolving subcategories for gentle trees}
	\label{sec:TreePart12}
	\pagestyle{plain}
In this section we compute all monogeneous resolving subcategories of a gentle tree algebra. In \cref{ss:Monoaccordsets}  we propose a geometric model of certain resolving subcategories and in \cref{ss:MonoResCat}, we show that they are all the resolving closures. 
\subsection{Simpler criterion for resolving subcategories of gentle tree algebras}
\label{ss:ressimplier}

This section aims to prove a valuable criterion for resolving subcategories of $\rep(Q,R)$ where $(Q,R)$ is a representation-finite gentle quiver.

\begin{theorem}
\label{thm:equivres}
Let $(Q,R)$ be a gentle tree. Let $\mathscr{D}$ be an additive subcategory of $\rep(Q,R)$. Then $\mathscr{D}$ is resolving if and only if $\mathscr{D}$ satisfies \ref{R'1}, \ref{R'2}, and:
\begin{enumerate}[label=$(\mathsf{Syg})$, itemsep=1mm]
    \item \label{Syg} $\mathscr{D}$ is closed under syzygies of indecomposable representations.
\end{enumerate}
\end{theorem}

To prove this theorem, we will use the geometric model we introduced for gentle quivers. Let us prove first the following lemma.

\begin{lemma}
\label{lem:kerbyextindec}
    Let $(Q,R)$ be a gentle tree. Let $\mathscr{D}$ be an additive subcategory of $\rep(Q,R)$ satisfying \ref{R'1} and \ref{R'2}. Then $\mathscr{D}$ satisfies \ref{Syg} if and only if $\mathscr{D}$ satisfies \ref{R3}.
\end{lemma}

\begin{proof}
By definition, if $\mathscr{D}$ satisfies \ref{R3}, then $\mathscr{D}$ satisfies \ref{Syg}. 
Assume that $\mathscr{D}$ is closed under syzygies of indecomposable representations. To show that $\mathscr{D}$ is closed under kernels of epimorphisms whose codomain is an indecomposable representation, it suffices to show that $\mathscr{D}$ is closed under kernels of minimal epimorphisms of the same type. 

Consider a minimal epimorphism  \[ \begin{tikzcd}
	{\displaystyle f: \bigoplus_{i=1}^p \MM(\delta_i)} & \MM(\eta)
	\arrow[two heads, from=1-1, to=1-2]
\end{tikzcd}.\] Set $\Ker(f) = \bigoplus_{i=0}^p \MM(\kappa_i)$ as described in \cref{prop:geo_mor}.
Assume that $\kappa_0$ is a nontrivial accordion. Then $\eta$ starts in a large cell $C$. Let $\mu \in \partial C$ such that $\mu$ crosses $\eta$. Let $\varsigma \in \NP(\eta) \cap \partial C$ be the unique projective accordion sharing with $\mu$ a common endpoint $e$. By \cref{prop:arcsaccordions}, $\delta_1$ must cross $\varsigma$. This implies that $\delta_1$ crosses the accordion $\nu$ whose endpoints are $s(\eta)$ and $e$, corresponding to an indecomposable summand of $\Omega(\MM(\eta))$ by \cref{cor:geom_syzygy}. The accordion $\kappa_0$ is associated to an indecomposable summand of the extension of $\MM(\delta_1)$ and $\MM(\nu)$. The accordion $\kappa_p$ is defined dually.

Let $i \in \{1,\ldots, p-1\}$. As $\kappa_i$ are nontrivial accordions, there exists $\rho_i \in \NP(\eta)$ such that either $\kappa_i$ crosses $\rho_i$, or shares with $\rho_i$ a common endpoint. Thus $\delta_i$ and $\delta_{i+1}$ either cross or share with $\rho_i$ a common endpoint. So:
\begin{enumerate}[label=$\bullet$, itemsep=1mm]
    \item if both $\delta_{i+1}$ and $\delta_i$ share an endpoint with $\rho_i$, as those endpoints are common with $\kappa_i$, then $\rho_i=\kappa_i$;
    \item if $\rho_i$ shares a common endpoint only with $\delta_i$, which is its common endpoint with $\kappa_i$, then $\rho_i$ crosses $\delta_{i+1}$: we have that $\OvExt(\delta_{i+1}, \rho_i) \neq \varnothing$, and $\kappa_i \in \OvExt(\delta_{i+1}, \rho_i)$;
    \item the case where $\rho_i$ shares a common endpoint only with $\delta_{i+1}$ can be treated as the previous case, up to symmetry; and,
    \item if both $\delta_i$ and $\delta_{i+1}$ cross $\rho_i$ (see \cref{fig:Co-Z}), then $\OvExt(\delta_i, \rho_i) \neq \varnothing \neq \OvExt(\delta_{i+1}, \rho_i)$: there exists exactly one accordion $\varsigma_i \in \OvExt(\delta_i, \rho_i)$ which crosses $\delta_{i+1}$, and $\kappa_i \in \OvExt(\varsigma_i, \delta_{i+1})$.
\end{enumerate}
\begin{figure}[ht!]
\centering 
    \begin{tikzpicture}[mydot/.style={
					circle,
					thick,
					fill=white,
					draw,
					outer sep=0.5pt,
					inner sep=1pt
				}, scale = 1.2]
		\tikzset{
		osq/.style={
        rectangle,
        thick,
        fill=white,
        append after command={
            node [
                fit=(\tikzlastnode),
                orange,
                line width=0.3mm,
                inner sep=-\pgflinewidth,
                cross out,
                draw
            ] {}}}}
        \draw [line width=0.7mm,domain=70:80] plot ({4*cos(\x)}, {1.5*sin(\x)});
        \draw [line width=0.7mm,domain=100:110] plot ({4*cos(\x)}, {1.5*sin(\x)});
        \draw [line width=0.7mm,domain=130:140] plot ({4*cos(\x)}, {1.5*sin(\x)});
        \draw [line width=0.7mm,domain=185:205] plot ({4*cos(\x)}, {1.5*sin(\x)});
        \draw [line width=0.7mm,domain=250:260] plot ({4*cos(\x)}, {1.5*sin(\x)});
        \draw [line width=0.7mm,domain=280:290] plot ({4*cos(\x)}, {1.5*sin(\x)});
		\foreach \X in {0,1,...,23}
		{
		\tkzDefPoint(4*cos(pi/12*\X),1.5*sin(pi/12*\X)){\X};
		};

		\draw[line width=0.9mm,bend left =30,red](9) edge (19);
		
		\draw[line width=0.7mm,bend right=10,blue, loosely dashed](5) edge (17);
		\draw[line width=0.7mm,bend left=10,blue, loosely dashed](7) edge (13);
        \draw[line width=0.7mm,bend left=0,darkgreen, dash pattern={on 5pt off 2pt on 1pt off 2pt}](13) edge (19);
		\draw[line width=0.7mm,bend left=10,orange, densely dashdotted](5) edge (13);
		

		\foreach \X in {5,7,9,13,17,19}
		{
		\tkzDrawPoints[fill =red,size=4,color=red](\X);
		};

		\tkzDefPoint(-3.3,0.5){gamma};
		\tkzLabelPoint[blue](gamma){\Large $\delta_i$}
		\tkzDefPoint(1,1){gamma};
		\tkzLabelPoint[blue](gamma){\Large $\delta_{i+1}$}
		\tkzDefPoint(-0.2,1.3){gamma};
		\tkzLabelPoint[orange](gamma){\Large $\kappa_i$}
		\tkzDefPoint(1.2,-0.7){s};
		\tkzLabelPoint[red](s){\Large $\rho_i$}
        \tkzDefPoint(-2,-0.3){s};
		\tkzLabelPoint[darkgreen](s){\Large $\varsigma_i$}
    \end{tikzpicture}
\caption{\label{fig:Co-Z} Illustration of the last case we enumerated above. The other cases can be obtained by identifying the endpoints of $\rho_i$ with those of $\kappa_i$.}
\end{figure}
In either case, we get that $\MM(\kappa_i) \in \mathscr{D}$ for all $i \in \{0,\ldots,p\}$, and so $\Ker(f) \in \mathscr{D}$. Therefore, $\mathscr{D}$ satisfies \ref{R3}.
\end{proof}

\begin{proof}[Proof of \cref{thm:equivres}]
Let $\mathscr{D}$ be an additive subcategory of $\rep(Q,R)$. If $\mathscr{D}$ is resolving, it is obvious that $\mathscr{D}$ satisfies \ref{R'1}, \ref{R'2} and \ref{Syg}. 

Now assume that $\mathscr{D}$ satisfies \ref{R'1}, \ref{R'2} and \ref{Syg}. By \cref{lem:kerbyextindec}, these properties are equivalent to $\mathscr{D}$ satisfying \ref{R'1}, \ref{R'2} and \ref{R'3}. By \cref{thm:rescondonindec}, we get to the conclusion that $\mathscr{D}$ is resolving.
\end{proof}

\subsection{Monogeneous geometric resolving sets}
\label{ss:Monoaccordsets}

Let $(Q,R)$ be a gentle tree. In this section, we introduce \new{monogeneous accordion sets} on $\Surf(Q,R)$, which are sets of accordions that are closed under geometric operations coming from \cref{prop:geo_mor,prop:geom_ext}. In the next section, we will show that they correspond bijectively to the sets of indecomposable representations appearing in monogeneous resolving subcategories of $\rep(Q,R)$. 

Let $\delta \in \Accord$ equipped with an orientation \(\underrightarrow{\delta}\). We call  \(\underrightarrow{\delta}(0)\) the \emph{source} of $\delta$, and \(\underrightarrow{\delta}(1)\) the \emph{target}.  As $\delta$ cuts the disc into two connected parts, we consider the \emph{above part} of $\delta$ to be the part of the disc where the boundary is oriented from the target to the source of $\delta$, and the \emph{below part} is the one oriented from the source to the target of $\delta$. 

\begin{definition}
\label{def:colourendpoints}
Let $(Q,R)$ be a gentle tree and $\Surf(Q,R) = (\pmb{\Sigma}, \mathcal{M}, \Delta^{\gpoint})$. Let $\delta \in \Accord'$ equipped with an orientation $\underrightarrow{\delta}$. We define a \new{coloration} of $\NP(\delta)_0$ to be a partition of $\NP(\delta)_0$ as  follows:
\begin{enumerate}[label=\arabic*),itemsep=1mm]
    \item  the points of $\NP(\delta)_0$ in the above part of $\delta$ are colored red $\rsquare$, and those in the below part of $\delta$ are in green $\gsquare$;
    \item The source is colored red $\rsquare$ if all the accordions in $\NP(\delta)$ having the source of $\delta$ as an endpoint admit a green endpoint $\gsquare$ at the other end. It is colored orange $\osquare$ otherwise;
    \item the target is colored dually;
    \item in larger cells, the coloration is changed: all the vertices except the endpoint of the edge that is in $\partial\pmb{\Sigma}$ and which is not the orange $\osquare$ endpoint of $\delta$ are colored orange $\osquare$; and,
    \item the endpoint of the border of the large cell which has a projective accordion which admits the orange endpoint of $\delta$ as an endpoint is finally colored in pink $\psquare$.
\end{enumerate}
We denote by:
\begin{enumerate}[label=$\bullet$,itemsep=1mm]
    \item $\NP(\delta)_{0}^{{\rsquare}}$ the set of points in $\NP(\delta)_0$ colored in red $\rsquare$;
    \item $\NP(\delta)_{0}^{{\osquare}}$ the set of points in $\NP(\delta)_0$ colored in orange $\osquare$;
    \item $\NP(\delta)_{0}^{{\gsquare}}$ the set of points in $\NP(\delta)_0$ colored in green $\gsquare$;
    \item $\NP(\delta)_{0}^{{\psquare}}$ the set of points in $\NP(\delta)_0$ colored in pink $\psquare$.
\end{enumerate}
\end{definition}

\begin{remark}
Note that the coloration of $\NP(\delta)_0$ is unique up to exchanging the ${\rsquare}$ and ${\gsquare}$ points. 
\end{remark}

Given $\delta \in \Accord'$ equipped with an orientation $\underrightarrow{\delta}$, we order $\NP(\delta)$ as follows:
\begin{enumerate}[label=$\bullet$,itemsep=1mm]
    \item We order the projective accordions that have the source of $\delta$ as an endpoint in the opposite of the counterclockwise order;
    \item If we have a large cell $C$ in the beginning of $\delta$, we order the projective accordions in $\partial C$ inductively as follows:
    \begin{enumerate}[label=$\bullet$, itemsep=1mm]
        \item Let $\rho_1$ be the one with the source of $\delta$ as an endpoint, and, whenever it is possible, let $\rho_{i+1}$ be the unique projective accordion in $\partial C$ such that it has a common endpoint with $\rho_i$ distinct from $\rho_{i-1}$. 
        \item We order those accordions by their index.
    \end{enumerate}
    \item For the projective accordions crossing $\delta$, we order them thanks to their crossing point and $\underrightarrow{\delta}$.
    \item We do the dual in the first two points for the remaining projective accordions.
\end{enumerate}
We denote this order $\overset{\text{\tiny{$\underrightarrow{\delta}$}}}{\rightarrowtriangle}$.

Now, whenever we take $\delta'\in\Accord'$ whose endpoints are in $\NP(\delta)_0$, we equipped $\delta'$ with an orientation $\underrightarrow{\delta'}$ such that $\overset{\text{\tiny{$\underrightarrow{\delta'}$}}}{\rightarrowtriangle}$ and $\overset{\text{\tiny{$\underrightarrow{\delta}$}}}{\rightarrowtriangle}$ coincides in $\NP(\delta) \cap \NP(\delta')$. This orientation is unique as we have at least two projective accordions in $\NP(\delta) \cap \NP(\delta')$.

\begin{conv}
From now on, given such a coloration of $\NP(\delta)_0$, for any non-projective accordion whose endpoints are in $\NP(\delta)_0$, we consider it oriented so that its orientation is compatible with the one of $\delta$.
\end{conv}

\begin{ex} We consider \cref{ex:Nproj1}.
We give in \cref{fig:Colorex} a coloration of $\NP(\delta)_0$. Note that the points $\rpoint$ are points in $\mathcal{M}_{\rpoint} \setminus \NP(\delta)_0$.
     \begin{figure}[!ht]
\centering 
    \begin{tikzpicture}[mydot/.style={
					circle,
					thick,
					fill=white,
					draw,
					outer sep=0.5pt,
					inner sep=1pt
				}, scale = 1]
		\tikzset{
		osq/.style={
        rectangle,
        thick,
        fill=white,
        append after command={
            node [
                fit=(\tikzlastnode),
                orange,
                line width=0.3mm,
                inner sep=-\pgflinewidth,
                cross out,
                draw
            ] {}}}}
        \tikzset{
		psq/.style={
        circle,
        thick,
        fill=white,
        append after command={
            node [
                fit=(\tikzlastnode),
                darkpink,
                line width=0.3mm,
                inner sep=-\pgflinewidth,
                cross out,
                draw
            ] {}}}}
		\draw[line width=0.7mm,black] (0,0) ellipse (4cm and 1.5cm);
			\foreach \X in {0,1,...,37}
		{
		\tkzDefPoint(4*cos(pi/19*\X),1.5*sin(pi/19*\X)){\X};
		};
        
		\draw[line width=0.9mm,bend left =60,red, densely dashdotted](1) edge (5);
		\draw[line width=0.9mm,bend left =60,red, densely dashdotted](3) edge (5);
		\draw[line width=0.9mm,bend right =30,red, densely dashdotted](5) edge (35);
		\draw[line width=0.9mm,bend left =30,red,densely dashdotted](33) edge (35);
		\draw[line width=0.9mm,bend left =30,red](35) edge (37);
		\draw[line width=0.9mm,bend left =30,red, densely dashdotted](31) edge (33);
		\draw[line width=0.9mm,bend left =30,red, densely dashdotted](7) edge (31);
		\draw[line width=0.9mm,bend left =30,red, densely dashdotted](7) edge (29);
		\draw[line width=0.9mm,bend left =30,red, densely dashdotted](9) edge (29);
		\draw[line width=0.9mm,bend left =30,red, densely dashdotted](11) edge (29);
		\draw[line width=0.9mm,bend left =30,red, densely dashdotted](11) edge (27);
		\draw[line width=0.9mm,bend right =30,red, densely dashdotted](13) edge (11);
		\draw[line width=0.9mm,bend left =30,red, densely dashdotted](13) edge (15);
		\draw[line width=0.9mm,bend left =30,red, densely dashdotted](17) edge (23);
		\draw[line width=0.9mm,bend left =30,red, densely dashdotted](15) edge (17);
		\draw[line width=0.9mm,bend left =30,red, densely dashdotted](19) edge (23);
		\draw[line width=0.9mm,bend left =30,red, densely dashdotted](21) edge (23);
		\draw[line width=0.9mm,bend left =30,red](25) edge (27);

		\draw[line width=0.7mm,bend right=10,blue, loosely dashed](5) edge (23);
        
		\foreach \X in {25,37}
		{
		\tkzDrawPoints[fill =red,size=4,color=red](\X);
		};
		\foreach \X in {7,9,19,21}
		{
		\tkzDrawPoints[rectangle,fill =red,size=6,color=red](\X);
		};
		
		\foreach \X in {1,3,27,29}
		{
		\tkzDrawPoints[rectangle,size=6,color=dark-green,thick,fill=white](\X);
		};
		\foreach \X in {5,11,13,15,23,31,33}
		{
		\tkzDrawPoints[size=6,orange,osq](\X);
		};
		\foreach \X in {17,35}
		{
		\tkzDrawPoints[size=6,darkpink,psq](\X);
		};
		\begin{scope}[xshift=2.5ex]
		\tkzDefPoint(-1.5,0.7){gamma};
		\tkzLabelPoint[blue](gamma){\Large $\delta$}
        \tkzDefPoint(-3.6,-1){s};
		\tkzLabelPoint[orange](s){ $s(\delta)$}
        \tkzDefPoint(2.5,1.8){t};
		\tkzLabelPoint[orange](t){ $t(\delta)$}
		\end{scope}
    \end{tikzpicture}
\caption{\label{fig:Colorex} An example of a coloration of $\NP(\delta)_0$.}
\end{figure}
\end{ex}

\begin{definition} \label{def:geometricres}
Let $(Q,R)$ be a gentle tree, and $\Surf(Q,R) = (\pmb{\Sigma}, \mathcal{M}, \Delta^{\gpoint})$. For any $\delta \in \Accord'$, we define the \new{monogeneous geometric resolving set of $\delta$} as the set $\opResAc(\delta) = \opResAc'(\delta) \cup \Prj(\Delta^{\gpoint})$ where: \[\opResAc'(\delta) = \{ \eta \in \Accord' \mid s(\eta) \in \NP(\delta)_0^{{\rsquare}} \cup \NP(\delta)_0^{{\osquare}}, t(\eta) \in \NP(\delta)_0^{{\gsquare}} \cup \NP(\delta)_0^{{\osquare}} \}.\]
The \new{monogeneous geometric resolving subcategory of $\delta$}, denoted by \new{$\mathscr{U}_\delta$}, is the additive subcategory of $\rep(Q,R)$ generated by $\{\MM(\varsigma) \mid \varsigma \in \opResAc(\delta)\}$.
\end{definition}

\begin{ex} \label{ex:ResM}  Let $(\pmb{\Sigma}, \mathcal{M}, \Delta^{\gpoint})$ be the $\gpoint$-dissection of the marked surface as seen in \cref{ex:Nproj1} and in \cref{fig:Colorex}.  Consider $\delta \in \Accord'$ as previously. In \cref{fig:exResM}, we represent the accordions in $\opResAc(\delta)$ with dotted and dashed lines. More precisely,
\begin{enumerate}[label=$\bullet$,itemsep=1mm]
    \item the blue loosely-dotted line is the accordion $\delta$; 
    \item the red densely dotted lines are the neighboring projectives of $\delta$ as seen in \cref{ex:Nproj1} and in \cref{fig:Colorex}; and,
    \item the purple dash-dotted lines are all the accordions in $\opResAc'(\delta)$. \qedhere
\end{enumerate}  
    \begin{figure}[!ht]
\centering 
    \begin{tikzpicture}[mydot/.style={
					circle,
					thick,
					fill=white,
					draw,
					outer sep=0.5pt,
					inner sep=1pt
				}, scale = 1.2]
		\tikzset{
		osq/.style={
        rectangle,
        thick,
        fill=white,
        append after command={
            node [
                fit=(\tikzlastnode),
                orange,
                line width=0.3mm,
                inner sep=-\pgflinewidth,
                cross out,
                draw
            ] {}}}}
        \tikzset{
		psq/.style={
        circle,
        thick,
        fill=white,
        append after command={
            node [
                fit=(\tikzlastnode),
                darkpink,
                line width=0.3mm,
                inner sep=-\pgflinewidth,
                cross out,
                draw
            ] {}}}}
		\draw[line width=0.7mm,black] (0,0) ellipse (4cm and 1.5cm);
			\foreach \X in {0,1,...,37}
		{
		\tkzDefPoint(4*cos(pi/19*\X),1.5*sin(pi/19*\X)){\X};
		};

		\draw[line width=0.9mm,bend left =60,red, densely dashdotted](1) edge (5);
		\draw[line width=0.9mm,bend left =60,red,densely dashdotted](3) edge (5);
		\draw[line width=0.9mm,bend right =30,red, densely dashdotted](5) edge (35);
		\draw[line width=0.9mm,bend left =30,red,densely dashdotted](33) edge (35);
		\draw[line width=0.9mm,bend left =30,red](35) edge (37);
		\draw[line width=0.9mm,bend left =30,red, densely dashdotted](31) edge (33);
		\draw[line width=0.9mm,bend left =30,red, densely dashdotted](7) edge (31);
		\draw[line width=0.9mm,bend left =30,red, densely dashdotted](7) edge (29);
		\draw[line width=0.9mm,bend left =30,red, densely dashdotted](9) edge (29);
		\draw[line width=0.9mm,bend left =30,red, densely dashdotted](11) edge (29);
		\draw[line width=0.9mm,bend left =30,red, densely dashdotted](11) edge (27);
		\draw[line width=0.9mm,bend right =30,red, densely dashdotted](13) edge (11);
		\draw[line width=0.9mm,bend left =30,red, densely dashdotted](13) edge (15);
		\draw[line width=0.9mm,bend left =30,red, densely dashdotted](17) edge (23);
		\draw[line width=0.9mm,bend left =30,red, densely dashdotted](15) edge (17);
		\draw[line width=0.9mm,bend left =30,red, densely dashdotted](19) edge (23);
		\draw[line width=0.9mm,bend left =30,red,densely dashdotted](21) edge (23);
		\draw[line width=0.9mm,bend left =30,red](25) edge (27);
		
		\draw[line width=0.7mm,blue, loosely dashed](5) edge (23);
		
		\draw [line width=0.7mm, mypurple,dash pattern={on 10pt off 2pt on 5pt off 2pt}, bend right=10] (5) edge (33);
		\draw [line width=0.7mm, mypurple,dash pattern={on 10pt off 2pt on 5pt off 2pt}, bend right=-10] (5) edge (31);
		\draw [line width=0.7mm, mypurple,dash pattern={on 10pt off 2pt on 5pt off 2pt}, bend right=10] (1) edge (33);
		\draw [line width=0.7mm, mypurple,dash pattern={on 10pt off 2pt on 5pt off 2pt}, bend right=10] (3) edge (33);
		\draw [line width=0.7mm, mypurple,dash pattern={on 10pt off 2pt on 5pt off 2pt}, bend left=10] (5) edge (7);
		\draw [line width=0.7mm, mypurple,dash pattern={on 10pt off 2pt on 5pt off 2pt}, bend left=10] (5) edge (9);
		\draw [line width=0.7mm, mypurple,dash pattern={on 10pt off 2pt on 5pt off 2pt}, bend left=10] (5) edge (11);
		\draw [line width=0.7mm, mypurple,dash pattern={on 10pt off 2pt on 5pt off 2pt}, bend right=10] (23) edge (11);
		\draw [line width=0.7mm, mypurple,dash pattern={on 10pt off 2pt on 5pt off 2pt}, bend right=10] (23) edge (13);
		\draw [line width=0.7mm, mypurple,dash pattern={on 10pt off 2pt on 5pt off 2pt}, bend right=10] (23) edge (15);
		
		\draw [line width=0.7mm, mypurple,dash pattern={on 10pt off 2pt on 5pt off 2pt}, bend right=20] (19) edge (15);
		\draw [line width=0.7mm, mypurple,dash pattern={on 10pt off 2pt on 5pt off 2pt}, bend right=20] (21) edge (15);
		
		\draw [line width=0.7mm, mypurple,dash pattern={on 10pt off 2pt on 5pt off 2pt}, bend left=20] (23) edge (27);
		\draw [line width=0.7mm, mypurple,dash pattern={on 10pt off 2pt on 5pt off 2pt}, bend left=15] (23) edge (29);
		\draw [line width=0.7mm, mypurple,dash pattern={on 10pt off 2pt on 5pt off 2pt}, bend left=15] (23) edge (31);
		\draw [line width=0.7mm, mypurple,dash pattern={on 10pt off 2pt on 5pt off 2pt}, bend left=10] (11) edge (31);
		\draw [line width=0.7mm, mypurple,dash pattern={on 10pt off 2pt on 5pt off 2pt}, bend left=10] (9) edge (31);
		
		\foreach \X in {25,37}
		{
		\tkzDrawPoints[fill =red,size=4,color=red](\X);
		};
		\foreach \X in {7,9,19,21}
		{
		\tkzDrawPoints[rectangle,fill =red,size=6,color=red](\X);
		};
		
		\foreach \X in {1,3,27,29}
		{
		\tkzDrawPoints[rectangle,size=6,color=dark-green,thick,fill=white](\X);
		};
		\foreach \X in {5,11,13,15,23,31,33}
		{
		\tkzDrawPoints[size=6,orange,osq](\X);
		};
		\foreach \X in {17,35}
		{
		\tkzDrawPoints[size=6,darkpink,psq](\X);
		};
		\begin{scope}[xshift=2.5ex]
		\tkzDefPoint(-1.4,-0.2){gamma};
		\tkzLabelPoint[blue](gamma){\Large $\delta$}
		\end{scope}
    \end{tikzpicture}
\caption{\label{fig:exResM} An explicit calculation of  $\opResAc(\delta)$, associated to $\NP(\delta)$ of \cref{ex:Nproj1}. }
\end{figure}
\end{ex}

\begin{lemma} \label{lem:Accordset_ext}
Let $(Q,R)$ be a gentle tree, and $\Surf(Q,R) = (\pmb{\Sigma}, \mathcal{M}, \Delta^{\gpoint})$. Consider $\delta \in \Accord'$. Then, for any pair $(\varsigma_1, \varsigma_2) \in \opResAc(\delta)^2$, we have \[\OvExt(\varsigma_1, \varsigma_2) \cup \ArExt(\varsigma_1,\varsigma_2) \subseteq \opResAc(\delta).\]
\end{lemma}

\begin{proof}

    Consider first the case of overlap extensions. Let $(\varsigma_1, \varsigma_2)\in \opResAc(\delta)^2$ such that $\OvExt(\varsigma_1, \varsigma_2) \neq \varnothing$. By \cref{prop:geom_ext}, the extension is involving at most two accordions $\upsilon_1$ and $\upsilon_2$ and as the four curves share a common part the curves are oriented following the orientation of $\delta$ thus their sources are either the source of $\varsigma_1$ or the source of $\varsigma_2$ and the targets either the target of $\varsigma_1$ or the target of $\varsigma_2$. We schematize the situation in \cref{fig:resclo1extclos}. \begin{figure}[!ht]
\centering 
    \begin{tikzpicture}[mydot/.style={
					circle,
					thick,
					fill=white,
					draw,
					outer sep=0.5pt,
					inner sep=1pt
				}, scale = 1]
		\tikzset{
		osq/.style={
        rectangle,
        thick,
        fill=white,
        append after command={
            node [
                fit=(\tikzlastnode),
                orange,
                line width=0.3mm,
                inner sep=-\pgflinewidth,
                cross out,
                draw
            ] {}}}}
		\draw [line width=0.7mm,domain=50:130] plot ({4*cos(\x)}, {1.5*sin(\x)});
        \draw [line width=0.7mm,domain=230:310] plot ({4*cos(\x)}, {1.5*sin(\x)});
		\foreach \X in {0,1}
		{
		\tkzDefPoint(4*cos(pi/6*\X +pi/2),1.5*sin(pi/6*\X + pi/2)){\X};
		};
		\foreach \X in {2,3}
		{
		\tkzDefPoint(4*cos(pi/6*(\X-2) +3*pi/2),1.5*sin(pi/6*(\X-2) + 3*pi/2)){\X};
		};
		
		\draw[line width=0.7mm,bend right=10,blue!70,dotted](0) edge (2);
		
		\draw[line width=0.7mm,bend right=20,blue, loosely dotted](1) edge (3);
		
		\draw [line width=0.7mm, mypurple,dash pattern={on 10pt off 2pt on 5pt off 2pt}, bend right=20] (1) edge (0);
		\draw [line width=0.7mm, mypurple,dash pattern={on 10pt off 2pt on 5pt off 2pt}, bend left=20] (2) edge (3);
		
		\foreach \X in {0,1,2,3}
		{
		\tkzDrawPoints[fill =red,size=4,color=red](\X);
		};
		
		\begin{scope}[xshift=2.5ex]
		\tkzDefPoint(-.1,1){gammaM};
		\tkzLabelPoint[blue](gammaM){\Large $\varsigma_1$}
		\tkzDefPoint(-1.4,0.2){gammaP};
		\tkzLabelPoint[blue!70](gammaP){\Large $\varsigma_2$}
		\tkzDefPoint(0.1,-0.8){deltav};
		\tkzLabelPoint[mypurple](deltav){\Large $\upsilon_2$}
		\tkzDefPoint(-1.2,1.2){deltaw};
		\tkzLabelPoint[mypurple](deltaw){\Large $\upsilon_1$}
		\end{scope}
    \end{tikzpicture}
\caption{\label{fig:resclo1extclos} Accordions $\upsilon_1$ and $\upsilon_2$ involving in the overlap extension of $\varsigma_1$ by $\varsigma_2$.}
\end{figure}
Up to exchanging the role of $\upsilon_1$ and $\upsilon_2$, assume that $s(\upsilon_1)=s(\varsigma_2)$, $t(\upsilon_1)=t(\varsigma_1)$, $s(\upsilon_2)=s(\varsigma_1)$ and $t(\upsilon_2)=t(\varsigma_2)$. Then $\upsilon_1, \upsilon_2 \in \Accord$, and, as $ s(\varsigma_i) \in \NP(\delta)_0^{{\rsquare}} \cup \NP(\delta)_0^{{\osquare}}$ and $t(\varsigma_i)\in \NP(\delta)_0^{{\osquare}} \cup \NP(\delta)_0^{{\gsquare}}$ for $i \in \{1,2\}$, it occurs that $\upsilon_1, \upsilon_2 \in \opResAc(\delta)$. 

Consider arrow extensions. Let $(\varsigma_1, \varsigma_2)\in \opResAc(\delta)^2$ such that $\ArExt(\varsigma_1, \varsigma_2) \neq \varnothing$. The curve corresponding to the indecomposable summand of the arrow extension $\theta$ is such that it crosses all the projective curves crossed by $\varsigma_1$ and $\varsigma_2$. Thus the curves $\varsigma_1, \varsigma_2 $ and $\theta$ have a common orientation induced by $\delta$. The curve $\theta$ is such that $s(\theta)=s(\varsigma_1)$ and $t(\theta)=t(\varsigma_2)$ or the dual configuration. Up to renaming, consider this configuration. As $\varsigma_1\in\opResAc(\delta)$, $ s(\varsigma_1) \in \NP(\delta)_0^{{\rsquare}} \cup \NP(\delta)_0^{{\osquare}}$. Dually as $\varsigma_2\in\opResAc(\delta)$, $t(\varsigma_2)\in \NP(\delta)_0^{{\osquare}} \cup \NP(\delta)_0^{{\gsquare}}$. Thus $\theta\in \opResAc(\delta)$.
\end{proof}

\begin{lemma}\label{lem:Syzygyclo}
Let $(Q,R)$ be a gentle tree, and $\Surf(Q,R) = (\pmb{\Sigma}, \mathcal{M}, \Delta^{\gpoint})$. Consider $\delta \in \Accord'$. For any $\varsigma\in\opResAc(\delta)$, a non-projective summand of the syzygy, $\Omega(\varsigma)$ is such that $\Omega(\varsigma)\in \opResAc(\delta)$.
\end{lemma}

\begin{proof}
Let $(Q,R)$ be a gentle tree, and $\Surf(Q,R) = (\pmb{\Sigma}, \mathcal{M}, \Delta^{\gpoint})$. Consider $\delta \in \Accord'$. Let $\varsigma\in\opResAc(\delta)$ and let $\Omega(\varsigma)$ be a non-projective summand of the syzygy. By \cref{cor:geom_syzygy}, one can see that $\Omega(\varsigma)\in\opResAc(\delta)$.
\end{proof}

\begin{prop} \label{prop:UMRes}
    Let $(Q,R)$ be a gentle tree, and set $(\pmb{\Sigma},\mathcal{M}, \Delta^{\gpoint}) = \Surf(Q,R)$. For any $\delta \in \Accord'$, the monogeneous geometric resolving subcategory $\mathscr{U}_\delta$ is resolving.
\end{prop}

\begin{proof}
By using \cref{prop:red_ker} and \ref{prop:red_ext}, and since $\mathscr{U}_\delta$ is additive, it suffices to prove that $\ind(\mathscr{U}_\delta)$ is closed under taking summands from non-trivial extensions between indecomposable representations, and from syzygies by \cref{thm:equivres}.
By \cref{prop:geom_ext,cor:geom_syzygy}, we have the result from \cref{lem:Accordset_ext} and \cref{lem:Syzygyclo}.

\end{proof}

\subsection{Monogeneous resolving subcategories}
\label{ss:MonoResCat}

Let $(Q,R)$ be a gentle quiver. We recall that a resolving subcategory $\mathscr{R} \subseteq \rep(Q,R)$ is monogeneous whenever there exists $X \in \ind(Q,R)$ such that $\Res(X) = \mathscr{R}$. In this section, we will prove that the monogeneous resolving subcategories of $\rep(Q,R)$ are exactly the monogeneous geometric resolving subcategories. Here is the more precise statement.

\begin{theorem}
\label{thm:resclo1}
Let $(Q,R)$ be a gentle tree, and $\Surf(Q,R) = (\pmb{\Sigma}, \mathcal{M}, \Delta^{\gpoint})$. For any $\delta \in \Accord'$, we have that $\mathscr{U}_\delta = \Res(\MM(\delta))$.
\end{theorem}

Let us introduce some useful notations.

For any $\delta \in \Accord$ which is not contained in a cell of the $\rpoint$-dissection $(\pmb{\Sigma}, \mathcal{M}, \Prj(\Delta^{\gpoint}))$, define $\tc (\delta)$ to be the pair $(\eta, C)$ where $\eta \in \Prj(\Delta^{\gpoint})$, and $C \in \pmb{\Gamma}(\Prj(\Delta^{\gpoint}))$ admitting $\eta$ as one of its edges, such that $\delta$ crosses $\eta$ and $t(\delta)$ is a vertex of $C$. If $\delta \subset C$, then we set $\tc(\delta) = (\varnothing,C)$. We define $\sc(\delta)$ dually.

Given a $C \in \pmb{\Gamma}(\Prj(\Delta^{\gpoint}))$, we denote by $C(\delta)_0^{{\osquare}}$ the set of vertices of $C$ in $\NP(\delta)_0^{{\osquare}}$. We define similarly $C(\delta)_0^{{\gsquare}}$, $C(\delta)_0^{{\rpoint}}$ and $C(\delta)_0^{{\psquare}}$.

Define $\vartheta_1^R, \vartheta_2^R \in \Prj(\Delta^{\gpoint})$ such that:
\begin{enumerate}[label=$\bullet$,itemsep=1mm]
    \item $t(\vartheta_1^R) = t(\delta)$ and $\vartheta_1^R$ covers $\delta$ with respect to the counterclockwise order (\cref{def:counterclockorder}) around $t(\delta)$ on $\Prj(\Delta^{\gpoint}) \cup \{\delta\}$; and,
    \item $t(\vartheta_2^R) = s(\vartheta_1^R)$ and $\vartheta_2^R$ covers $\vartheta_1^R$ with respect to the counterclockwise order around $s(\vartheta_1^R)$ on $\Prj(\Delta^{\gpoint}) \cup \{\delta\}$.
\end{enumerate}
We define $w_R(\delta) = s(\vartheta_2^R)$. We define dually $w_L(\delta)$ as we defined $w_R(\delta)$ by exchanging the above and below of $\delta$ (i.e., exchanging red and green in the coloration). See \cref{fig:exsctc} to visualize the different objects introduced.

\begin{figure}[!ht]
    \centering
    \begin{tikzpicture}[mydot/.style={
					circle,
					thick,
					fill=white,
					draw,
					outer sep=0.5pt,
					inner sep=1pt
				}, scale = 1]
		\tikzset{
		osq/.style={
        rectangle,
        thick,
        fill=white,
        append after command={
            node [
                fit=(\tikzlastnode),
                orange,
                line width=0.3mm,
                inner sep=-\pgflinewidth,
                cross out,
                draw
            ] {}}}}
        \tikzset{
		psq/.style={
        circle,
        thick,
        fill=white,
        append after command={
            node [
                fit=(\tikzlastnode),
                darkpink,
                line width=0.3mm,
                inner sep=-\pgflinewidth,
                cross out,
                draw
            ] {}}}}
        \draw [line width=0.7mm,domain=10:22] plot ({5*cos(\x)}, {2*sin(\x)});
		\draw [line width=0.7mm,domain=25:70] plot ({5*cos(\x)}, {2*sin(\x)});
		\draw [line width=0.7mm,domain=78:85] plot ({5*cos(\x)}, {2*sin(\x)});
		\draw [line width=0.7mm,domain=95:102] plot ({5*cos(\x)}, {2*sin(\x)});
		\draw [line width=0.7mm,domain=110:170] plot ({5*cos(\x)}, {2*sin(\x)});
		\draw [line width=0.7mm,domain=225:235] plot ({5*cos(\x)}, {2*sin(\x)});
		\draw [line width=0.7mm,domain=250:283] plot ({5*cos(\x)}, {2*sin(\x)});
		\draw [line width=0.7mm,domain=300:355] plot ({5*cos(\x)}, {2*sin(\x)});
        
		\foreach \X in {0,...,43}
		{
		\tkzDefPoint(5*cos(pi/22*\X),2*sin(pi/22*\X)){\X};
		};

		\draw[line width=0.7mm,blue, loosely dashed](4) to[bend left=30] (2,0) to[bend right=30] (-1,-1) to [bend right=30] (31);
		
		\draw[line width=0.9mm,bend right=30,red, densely dashdotted](4) edge (43);
		\draw[line width=0.9mm,bend right=30,red, densely dashdotted](4) edge (2);
		\draw[line width=0.9mm,bend right=30,red, densely dashdotted](43) edge (41);
		\draw[line width=0.9mm,bend right=30,red, densely dashdotted](41) edge (39);
		\draw[line width=0.9mm,bend right=30,red, densely dashdotted](39) edge (37);
		\draw[line width=1.5mm,bend left=30,red, densely dashdotted](37) edge (8);
		\draw[line width=0.9mm,bend left=40,red, densely dashdotted](37) edge (10);
		\draw[line width=0.9mm,bend left=40,red, densely dashdotted](34) edge (12);
		\draw[line width=1.5mm,bend left=30,red, densely dashdotted](34) edge (14);
		\draw[line width=0.9mm,bend right=30,red, densely dashdotted](16) edge (14);
		\draw[line width=0.9mm,bend right=30,red, densely dashdotted](18) edge (16);
		\draw[line width=0.9mm,bend right=30,red, densely dashdotted](20) edge (18);
		\draw[line width=0.9mm,bend right=30,red, densely dashdotted](31) edge (20);
		\draw[line width=0.9mm,bend right=30,red, densely dashdotted](31) edge (28);

		\filldraw [fill=red,opacity=0.1] (4) to [bend right=30] (43) to [bend right=30] (41) to [bend right=30] (39) to [bend right=30] (37) to [bend left=30] (8) to [bend left=10] cycle;
		
		\filldraw [fill=red,opacity=0.1] (34) to [bend left=30] (14) to [bend left=30] (16) to [bend left=30] (18) to [bend left=30] (20) to [bend left=30] (31) to [bend right=10] cycle;
		
		\foreach \X in {8,10,12,28}
		{
		\tkzDrawPoints[rectangle,fill =red,size=6,color=red](\X);
		};
		
		\foreach \X in {2,34}
		{
		\tkzDrawPoints[rectangle,size=6,color=dark-green,thick,fill=white](\X);
		};
		\foreach \X in {4,14,16,18,31,37,39,41}
		{
		\tkzDrawPoints[size=6,orange,osq](\X);
		};
		\foreach \X in {20,43}
		{
		\tkzDrawPoints[size=6,darkpink,psq](\X);
		};
		
		\begin{scope}[xshift=2.5ex]
		\tkzDefPoint(-4.1,1){gammaM};
		\tkzLabelPoint[red](gammaM){\Large $C_L$}
		\tkzDefPoint(-1.8,1){gammaP};
		\tkzLabelPoint[red](gammaP){\Large $\pmb{\eta_L}$}
        \tkzDefPoint(3.75,.4){gammaP};
		\tkzLabelPoint[red](gammaP){$\vartheta_1^R$}
        \tkzDefPoint(4,-.1){gammaP};
		\tkzLabelPoint[red](gammaP){ $\vartheta_2^R$}
		\tkzDefPoint(1.1,1){v};
		\tkzLabelPoint[red](v){\Large $\pmb{\eta_R}$}
		\tkzDefPoint(2.5,1.5){w};
		\tkzLabelPoint[red](w){\Large $C_R$}
		\tkzDefPoint(0,-0.2){deltav};
		\tkzLabelPoint[blue](deltav){\Large $\delta$}
		\tkzDefPoint(-5,1.8){deltaw};
		\tkzLabelPoint[orange](deltaw){ $w_L(\delta)$}
		\tkzDefPoint(4.3,-0.9){w};
		\tkzLabelPoint[orange](w){$w_R(\delta)$}
		\end{scope}
    \end{tikzpicture}
    \caption{\label{fig:exsctc} Illustration of  $\sc(\delta) = (\eta_L, C_L)$ and $\tc(\delta) = (\eta_R, C_R)$ given $\delta \in \Accord$.}
\end{figure}


To prove \cref{thm:resclo1}, as we already know that $\mathscr{U}_\delta$ is a resolving subcategory that contains $\MM(\delta)$, it is enough to show that $\mathscr{U}_\delta$ is the smallest resolving subcategory containing $\MM(\delta)$, it is sufficient to show that indecomposable representations in $\mathscr{U}_\delta$ are in $\Res(\MM(\delta))$.

First of all, given $\delta \in \Accord'$, we enumerate necessary conditions on $\rpoint$-arcs of $\Surf(Q,R)$ to be accordions in $\opResAc'(\delta)=\{ \eta \in \Accord' \mid s(\eta) \in \NP(\delta)_0^{{\rsquare}} \cup \NP(\delta)_0^{{\osquare}}, t(\eta) \in \NP(\delta)_0^{{\gsquare}} \cup \NP(\delta)_0^{{\osquare}} \}.$

\begin{lemma} \label{lem:conditionsResAc}
Let $(\pmb{\Sigma}, \mathcal{M}, \Delta^{\gpoint})$ be a $\gpoint$-dissected marked disc with $\mathcal{M} \subset \partial \pmb{\Sigma}$. Let $\delta \in \Accord$. If $\delta \nsubseteq C$ for some $C \in \pmb{\Gamma}(\Prj(\Delta^{\gpoint}))$, then we set $\sc(\delta) = (\eta_L,C_L)$ and $\tc(\delta) = (\eta_R,C_R)$; otherwise, we set $C=C_L=C_R$ the cell containing $\delta$. The curve $\varsigma \in \opResAc'(\delta)$ must satisfy all of the following assertions: 
\begin{enumerate}[label=$(\roman*)$, itemsep=1mm]
        \item \label{1Accord} if $\varsigma$ crosses $\eta_L$ and $C_L$ is a large cell, then $s(\varsigma) = s(\delta)$, and, dually, if $\varsigma$ crosses $\eta_R$ and $C_R$ is a large cell, then $t(\varsigma) = t(\delta)$;
        \item \label{2Accord} if $\varsigma \subset C_L$, then $s(\varsigma)=s(\delta)$, and, dually, if $\varsigma \subset C_R$, then $t(\varsigma) = t(\delta)$;
        \item \label{3Accord} if $s(\varsigma)\in \NP(\delta)_0^{{\rsquare}}$ and $t(\varsigma) \in C_L(\delta)_0$, then $t(\varsigma)=w_L(\delta)$, and, dually, if $t(\varsigma)\in \NP(\delta)_0^{{\gsquare}}$ and $s(\varsigma) \in C_R(\delta)_0$, then $s(\varsigma)=w_R(\delta)$;
        \item \label{4Accord} If $s(\varsigma) \in C_L(\delta)_0^{{\osquare}}$ then  $s(\varsigma)=s(\delta)$ or $s(\varsigma)$ is an orange endpoint of $\eta_L$, and if $t(\varsigma) \in C_R(\delta)_0^{{\osquare}}$ then $t(\varsigma)=t(\delta)$ or $t(\varsigma)$ is an orange endpoint of $\eta_R$.
    \end{enumerate}
\end{lemma}

\begin{proof}
   Recall that we obtain the green dissection from a clockwise rotation of the neighboring projectives. Thus \cref{fig:exsctc} shows that all the listed cases correspond to the accordions.
\end{proof}

Now we will treat the case where $\delta$ is contained in a cell of $\pmb{\Gamma}(\Delta^{\gpoint})$.

\begin{lemma} \label{lem:UMinResinacell}
Let $(\pmb{\Sigma}, \mathcal{M}, \Delta^{\gpoint})$ be a $\gpoint$-dissected marked disc with $\mathcal{M} \subset \partial \pmb{\Sigma}$. Let $\delta \in \Accord'$  such that $\delta \subset C$ for some $C \in \pmb{\Gamma}(\Prj(\Delta^{\gpoint}))$. Then \[\opResAc(\delta) \subseteq \{\gamma_{(E)} \mid E \in \ind \left(\Res(\MM(\delta)) \right) \}.\]
\end{lemma}

\begin{proof}
 Up to a change of orientation of $\delta$, we assume that the unique cell is $C_R$ and thus we can define $w_R(\delta)$. We illustrate the situation in \cref{fig:resclo1inacell}. 
\begin{figure}[!ht]
\centering 
\begin{tikzpicture}[mydot/.style={
					circle,
					thick,
					fill=white,
					draw,
					outer sep=0.5pt,
					inner sep=1pt
				}, scale = 1]
		\tikzset{
		osq/.style={
        rectangle,
        thick,
        fill=white,
        append after command={
            node [
                fit=(\tikzlastnode),
                orange,
                line width=0.3mm,
                inner sep=-\pgflinewidth,
                cross out,
                draw
            ] {}}}}
        \tikzset{
		psq/.style={
        circle,
        thick,
        fill=white,
        append after command={
            node [
                fit=(\tikzlastnode),
                darkpink,
                line width=0.3mm,
                inner sep=-\pgflinewidth,
                cross out,
                draw
            ] {}}}}
		\draw [line width=0.7mm,domain=40:140] plot ({4*cos(\x)}, {2*sin(\x)});
        \draw [line width=0.7mm,domain=210:355] plot ({4*cos(\x)}, {2*sin(\x)});
		\foreach \X in {0,1}
		{
		\tkzDefPoint(4*cos(pi/3*\X +2*pi/5),2*sin(pi/3*\X + 2* pi/5)){\X};
		};
		\foreach \X in {2,3,4,5,6}
		{
		\tkzDefPoint(4*cos(pi/7*(\X-2) +4*pi/3),2*sin(pi/7*(\X-2) + 4*pi/3)){\X};
		};
		
		\tkzDefPoint(4*cos(pi/6 +2*pi/5),2*sin(pi/6 + 2* pi/5)){7};
		\tkzDefPoint(4*cos(-pi/7 +2*pi/5),2*sin(-pi/7 + 2* pi/5)){8};
		\tkzDefPoint(4*cos(-pi/14 +2*pi/5),2*sin(-pi/14 + 2* pi/5)){11};
		\tkzDefPoint(4*cos(-pi/8 +4*pi/3),2*sin(-pi/8 + 4*pi/3)){9}
		\tkzDefPoint(4*cos(-pi/12 +4*pi/3),2*sin(-pi/12 + 4*pi/3)){10}

		\draw[line width=0.7mm,bend right=10,blue, loosely dashed](0) edge (2);
		
		\draw[line width=0.9mm,bend right=30,red](1) edge (2);
		\draw[line width=0.9mm,bend left=30,red, densely dashdotted](2) edge (3);
		\draw[line width=0.9mm,bend left=30,red, densely dashdotted](3) edge (4);
		\draw[line width=0.9mm,bend left=30,red, densely dashdotted](4) edge (5);
		\draw[line width=0.9mm,bend left=30,red, densely dashdotted](5) edge (6);
		\draw[line width=0.9mm,bend right=30,red, densely dashdotted](0) edge (6);
		\draw[line width=0.9mm,bend right=30,red, densely dashdotted](0) edge (8);
		\draw[line width=0.9mm,bend right=30,red](1) edge (7);
		\draw[line width=0.9mm,bend left=30,red](9) edge (2);
		\draw[line width=0.9mm,bend left=30,red](10) edge (2);
		\draw[line width=0.9mm,bend left=30,red, densely dashdotted](11) edge (0);

		\draw [line width=0.7mm, mypurple,dash pattern={on 10pt off 2pt on 5pt off 2pt}, bend left=20] (3) edge (0);
		\draw [line width=0.7mm, mypurple,dash pattern={on 10pt off 2pt on 5pt off 2pt}, bend left=20] (4) edge (0);
		\draw [line width=0.7mm, mypurple,dash pattern={on 10pt off 2pt on 5pt off 2pt}, bend left=20] (5) edge (0);
		\draw [line width=0.7mm, mypurple,dash pattern={on 10pt off 2pt on 5pt off 2pt}, bend left=20] (5) edge (8);
		\draw [line width=0.7mm, mypurple,dash pattern={on 10pt off 2pt on 5pt off 2pt}, bend left=20] (5) edge (11);
		
		\filldraw [fill=red,opacity=0.1] (0) to [bend right=10] (7) to [bend left=30] (1) to [bend right=30] (2) to [bend left=30] (3) to [bend left=30] (4) to [bend left=30] (5) to [bend left=30] (6) to [bend left=30] cycle;
		
		\foreach \X in {1,7,9,10}
		{
		\tkzDrawPoints[fill =red,size=4,color=red](\X);
		};
		
		\foreach \X in {2}
		{
		\tkzDrawPoints[rectangle,fill =red,size=6,color=red](\X);
		};
		
		\foreach \X in {8,11}
		{
		\tkzDrawPoints[rectangle,size=6,color=dark-green,thick,fill=white](\X);
		};
		\foreach \X in {0,3,4,5}
		{
		\tkzDrawPoints[size=6,orange,osq](\X);
		};
		\foreach \X in {6}
		{
		\tkzDrawPoints[size=6,darkpink,psq](\X);
		};
		
		\begin{scope}[xshift=2.5ex]
		\tkzDefPoint(-1.1,1.3){gammaM};
		\tkzLabelPoint[blue](gammaM){\Large $\delta$}
		\tkzDefPoint(-2.5,1.3){gammaP};
		\tkzLabelPoint[red](gammaP){\Large $C$}
		\tkzDefPoint(2.5,-1.4){v};
		\tkzLabelPoint[orange](v){ $w_R(\delta)$}
		\end{scope}
    \end{tikzpicture}
\caption{\label{fig:resclo1inacell}  The case where $\delta$ is contained in a cell $C$ of the $\rpoint$-dissected marked disc $(\pmb{\Sigma}, \mathcal{M}, \Prj(\Delta^{\gpoint}))$. All the dashed lines are $\rpoint$-arcs in $\opResAc(\delta)$. The dual case of a left end-cell is obtained by rotating the figure of $180^\circ$, and swap red points and green points, and re-labelling $w_R(\delta)$ by $w_L(\delta)$.}
\end{figure}
We can note that:
\begin{enumerate}[label=$\bullet$, itemsep=1mm]
    \item There is exactly one vertex in $\NP(\delta)_0^{{\psquare}}$;
    \item there is no $\varsigma \in \Accord$ such that $\varsigma$ crosses $\delta$;
\end{enumerate}
 Suppose that $\varsigma \in \opResAc'(\delta)$ such that $\delta \neq \varsigma$, then either:
\begin{enumerate}[label = $\arabic*)$, itemsep=1mm]
    \item \label{2caseresinacell} if $t(\varsigma) \in \NP(\delta)_0^{{\osquare}}$, then $t(\varsigma) = t(\delta)$ by \cref{lem:conditionsResAc} (iv) as we assumed the existence of $w_R(\delta)$ ; or,
    \item \label{3caseresinacell} if $t(\varsigma) \in \NP(\delta)_0^{{\gsquare}}$, then $s(\varsigma) = w_R(\delta)$ by \cref{lem:conditionsResAc} (iii).
\end{enumerate}

If $\varsigma$ satisfies \ref{2caseresinacell}, as $\varsigma \notin \Prj(\Delta^{\gpoint})$ and $\delta \neq \varsigma$, then $s(\delta) \neq w_R(\delta)$. Consider the vertex $v_1 \in C(\delta)_0^{{\osquare}}$ defined as the other endpoint of the projective accordion that is the smallest with respect to the counterclockwise order around $s(\delta)$ among the projective accordions,  such that it is greater than $\delta$ with respect to the previous order. Note that, as $\delta$ is in a cell and thus does not cross a projective, in the projective resolution $(\PP_{\MM(\delta)}^\ast)$ of $\MM(\delta)$, $\PP_{\MM(\delta)}^0 = \MM(\nu)$ where $\nu \in \Prj(\Delta^{\gpoint})$ is such that $s(\nu) = s(\delta)$ and $t(\nu) = v_1$. Then the accordion $\delta_1$ such that $s(\delta_1) = v_1$ (and so $t(\delta_1) = t(\delta)$) corresponds to the indecomposable representation $\MM(\delta_1)$ which is isomorphic to the kernel of the unique (by \cref{prop:homtree}) epimorphism $\PP_{\MM(\delta)}^0 \twoheadrightarrow \MM(\delta)$. As $\Res(\MM(\delta))$ is closed under taking kernels of epimorphisms, we have that $\MM(\delta_1) \in \Res(\MM(\delta))$. By induction, for any $\varsigma$ satisfying \ref{2caseresinacell}, we have $\MM(\varsigma) \in \Res(\MM(\delta))$.

If $\varsigma \in \opResAc'(\delta)$ satisfies \ref{3caseresinacell}, then $\varsigma \in \ArExt(\rho_\varsigma,\delta_R)$ where:
\begin{enumerate}[label=$\bullet$, itemsep=1mm]
    \item $\delta_R \in \Accord$ with endpoints given by $s(\delta_R) = w_R(\delta)$, and $t(\delta_R) = t(\delta)$; and,
    \item $\rho_\varsigma \in \Accord$ with endpoints given by $s(\rho_\varsigma) = t(\delta)$ and $t(\rho_\varsigma) = t(\delta) \in \NP(\gamma_{(M)})_0^{{\gsquare}}$.
\end{enumerate} 
As $\delta_R$ satisfies \ref{2caseresinacell}, then $\MM(\delta_R) \in \Res(\MM(\delta))$. Moreover, $\MM(\rho_\varsigma) \in \proj(Q,R)$. Therefore, by \cref{thm:ExplicitExtGentle}, and \cref{def:resolv}, $\MM(\varsigma) \in \Res(\MM(\delta))$.

Thus, we get the desired result.
\end{proof}

We now consider the case where $\delta$ is not contained in a cell of $\pmb{\Gamma}(\Prj(\Delta^{\gpoint}))$.

\begin{lemma} \label{lem:UMinRes}
Let $(\pmb{\Sigma}, \mathcal{M}, \Delta^{\gpoint})$ be a $\gpoint$-dissected marked disc with $\mathcal{M} \subset \partial \pmb{\Sigma}$. Let $\delta \in \Accord'$ such that $\delta \nsubseteq C$ for all $C \in \pmb{\Gamma}(\Prj(\Delta^{\gpoint}))$. Then $\opResAc(\delta) \subseteq \{\gamma_{(E)} \mid E \in \ind \left(\Res(\MM(\delta)) \right) \}$.
\end{lemma}

\begin{proof}
We illustrated a general situation, up to dualization, in \cref{fig:exsctc}. As $\delta$ is not contained in a cell of $(\pmb{\Sigma}, \mathcal{M}, \Prj(\Delta^{\gpoint}))$, we can set $\sc(\delta) = (\eta_L, C_L)$ and $\tc(\delta) = (\eta_R, C_R)$.

Let us show first that any $\varsigma \in \Accord$ such that $s(\varsigma) = s(\delta)$ and $t(\varsigma) \in \NP(\delta)_0^{{\gsquare}}$ is in $\{\gamma_{(E)} \mid E \in \ind(\Res(\MM(\delta))) \}$. Consider $\varsigma$ as stated. There exists $\rho \in \NP(\delta)$ and  $t(\varsigma)$ is an endpoint of $\rho$ select the orientation of $\rho$ such that it is $t(\rho)$.  
\begin{enumerate}[label=$\bullet$,itemsep=1mm]
    \item If $s(\rho) \in C_R(\delta)_0$, then $t(\rho) = t(\delta)$: in this case, $\varsigma = \delta$ and we are done.
    \item Else, if $s(\rho) = s(\delta)$, then $\rho = \varsigma$, and we are done.
    \item Otherwise, $\rho$ crosses $\delta$, we have that $\OvExt(\rho, \delta) \neq \varnothing$ it follows from \cref{prop:geom_ext} that we have $\varsigma \in  \OvExt(\rho, \delta)$ (see \cref{fig:resclo1include}). Thus, as  $\Res(\MM(\delta))$ is closed under extensions, we also have $\MM(\varsigma) \in \Res(\MM(\delta)))$.
\end{enumerate} 
By an analogous argument, we can prove that any $\varsigma'\in \Accord$ such that $t(\varsigma') = t(\delta)$ and $s(\varsigma') \in \NP(\delta)_0^{{\rsquare}}$ (denoted by $w$ in \cref{fig:resclo1include}) is in $\{\gamma_{(E)} \mid E \in \ind(\Res(\MM(\delta)))\}$.
\begin{figure}[!ht]
\centering 
    \begin{tikzpicture}[mydot/.style={
					circle,
					thick,
					fill=white,
					draw,
					outer sep=0.5pt,
					inner sep=1pt
				}, scale = 1]
		\tikzset{
		osq/.style={
        rectangle,
        thick,
        fill=white,
        append after command={
            node [
                fit=(\tikzlastnode),
                orange,
                line width=0.3mm,
                inner sep=-\pgflinewidth,
                cross out,
                draw
            ] {}}}}
		\draw [line width=0.7mm,domain=50:130] plot ({4*cos(\x)}, {1.5*sin(\x)});
        \draw [line width=0.7mm,domain=230:310] plot ({4*cos(\x)}, {1.5*sin(\x)});
		\foreach \X in {0,1}
		{
		\tkzDefPoint(4*cos(pi/6*\X +pi/2),1.5*sin(pi/6*\X + pi/2)){\X};
		};
		\foreach \X in {2,3}
		{
		\tkzDefPoint(4*cos(pi/6*(\X-2) +3*pi/2),1.5*sin(pi/6*(\X-2) + 3*pi/2)){\X};
		};
		
		\draw[line width=0.9mm,bend right=10,red, densely dashdotted](0) edge (2);
		
		\draw[line width=0.7mm,bend right=30,blue, loosely dotted](1) edge (3);
		
		\draw [line width=0.7mm, mypurple,dash pattern={on 10pt off 2pt on 5pt off 2pt}, bend right=20] (1) edge (2);
		\draw [line width=0.7mm, mypurple,dash pattern={on 10pt off 2pt on 5pt off 2pt}, bend right=20] (0) edge (3);
		\foreach \X in {0}
		{
		\tkzDrawPoints[rectangle,fill =red,size=6,color=red](\X);
		};
		
		\foreach \X in {2}
		{
		\tkzDrawPoints[rectangle,size=6,color=dark-green,thick,fill=white](\X);
		};
        
		\begin{scope}[xshift=2ex]
		\tkzDefPoint(-.7,1){gammaM};
		\tkzLabelPoint[red](gammaM){\Large $\rho$}
		\tkzDefPoint(-1.1,0.4){gammaP};
		\tkzLabelPoint[blue](gammaP){\Large $\delta$}
		\tkzDefPoint(-0.4,-1.6){v};
		\tkzLabelPoint[dark-green](v){\Large $v$}
		\tkzDefPoint(0,1.9){w};
		\tkzLabelPoint[red](w){\Large $w$}
		\tkzDefPoint(0.6,0.5){deltav};
		\tkzLabelPoint[mypurple](deltav){\Large $\varsigma'$}
		\tkzDefPoint(-1.5,-0.5){deltaw};
		\tkzLabelPoint[mypurple](deltaw){\Large $\varsigma$}
		\end{scope}
    \end{tikzpicture}
\caption{\label{fig:resclo1include} Accordions (in purple dashed lines) corresponding to accordions appearing in $\OvExt(\rho, \delta)$.}
\end{figure}  By repeating the procedure on the curves obtained, we have that $\MM(\varsigma'') \in \Res(\MM(\delta))$ for all $\varsigma'' \in \Accord$ such that $s(\varsigma'') \in \NP(\delta)_0^{{\rsquare}}$ and $t(\varsigma'') \in \NP(\delta)_0^{{\gsquare}}$.
   
Now we focus on $\varsigma \in \opResAc'(\delta)$ such that  $s(\delta) \neq s(\varsigma)$, $t(\delta) \neq t(\varsigma)$ and either:
\begin{enumerate}[label=$\bullet$,itemsep = 1mm]
    \item $s(\varsigma),t(\varsigma)\in\NP(\gamma_{(M)})_0^{{\osquare}}$; or,
    \item $s(\varsigma) \in \NP(\delta)_0^{{\rsquare}}$ and $t(\varsigma) \in C_L$; or,
    \item $s(\varsigma) \in C_R$ and $t(\varsigma) \in \NP(\delta)_0^{{\gsquare}}$.
\end{enumerate}
By defining $\delta_L^\circ\in \Accord$ such that $s(\delta_L^\circ) = s(\delta)$ and $t(\delta_L^\circ)$ is the endpoint of $\eta_L$ which is not an endpoint to the unique edge of $\partial C_L\cap\partial\pmb{\Sigma}$.  We define $\delta_R^\circ \in \Accord$ such that $t(\delta^\circ_R)=t(\delta)$ and $s(\delta_R^\circ)$ is the endpoint of $\eta_R$ which is not an endpoint to the unique edge of $\partial C_R\cap\partial\pmb{\Sigma}$.  
Then the previous non-projective accordions $\varsigma$ appear either in $\opResAc'(\delta_L^\circ)$ or in $\opResAc'(\delta_R^\circ)$. As $\delta_L^\circ \subset C_L$ and $\delta_R^\circ \subset C_R$, by proving that $\MM(\delta_L^\circ),\MM(\delta_R^\circ) \in \Res(\MM(\delta))$, the result follows by \cref{lem:UMinResinacell}.

Consider $\vartheta$ the smallest accordion in $\Prj(\Delta^{\gpoint})$ such that  $\vartheta$ crosses $\delta$, including the case where its endpoints are $s(\delta)$, and $ s(\delta^\circ_R) \in C(\gamma_{(M)})_0^{{\osquare}}$, with respect to the clockwise order around $s(\delta^\circ_R)$ on $\Prj(\Delta^{\gpoint})$.
We depict the situation in \cref{fig:resclo1reduce}. \begin{figure}[!ht]
\centering 
    \begin{tikzpicture}[mydot/.style={
					circle,
					thick,
					fill=white,
					draw,
					outer sep=0.5pt,
					inner sep=1pt
				}, scale = 1]
		\tikzset{
		osq/.style={
        rectangle,
        thick,
        fill=white,
        append after command={
            node [
                fit=(\tikzlastnode),
                orange,
                line width=0.3mm,
                inner sep=-\pgflinewidth,
                cross out,
                draw
            ] {}}}}
        \tikzset{
		psq/.style={
        circle,
        thick,
        fill=white,
        append after command={
            node [
                fit=(\tikzlastnode),
                darkpink,
                line width=0.3mm,
                inner sep=-\pgflinewidth,
                cross out,
                draw
            ] {}}}}
		\draw [line width=0.7mm,domain=60:150] plot ({4*cos(\x)}, {2*sin(\x)});
        \draw [line width=0.7mm,domain=270:360] plot ({4*cos(\x)}, {2*sin(\x)});
		\foreach \X in {0,1,2,3,4}
		{
		\tkzDefPoint(4*cos(pi/11*\X+2*pi/5),2*sin(pi/11*\X+2* pi/5)){\X};
		};
		\foreach \X in {5,6,7,8}
		{
		\tkzDefPoint(4*cos(pi/11*(\X-5)+5*pi/3),2*sin(pi/11*(\X-5)+ 5*pi/3)){\X};
		};
		
		\draw[line width=0.7mm,bend left=10,blue, loosely dashed](0) edge (-4,-1);
		
		\draw[line width=0.9mm,bend right=30,red, densely dashdotted](1) edge (5);
		\draw[line width=0.9mm,bend right=30,red, densely dashdotted](2) edge (5);
		\draw[line width=0.9mm,bend right=30,red, densely dashdotted](3) edge (5);
		\draw[line width=1.5mm,bend right=30,red, densely dashdotted](4) edge (5);
		\draw[line width=0.9mm,bend left=30,red, densely dashdotted](4) edge (-3,0);
		\draw[line width=0.9mm,bend left=40,red, densely dashdotted](5) edge (6);
		\draw[line width=0.9mm,bend left=40,red, densely dashdotted](6) edge (7);
		\draw[line width=0.9mm,bend left=40,red, densely dashdotted](7) edge (8);
		\draw[line width=0.9mm,bend right=30,red, densely dashdotted](0) edge (8);

		\draw [line width=0.7mm, mypurple,dash pattern={on 10pt off 2pt on 5pt off 2pt}, bend left=20] (5) edge (0);
		
		\filldraw [fill=red,opacity=0.1] (0) to [bend right=15] (1) to [bend right=30] (5) to [bend left=30] (6) to [bend left=30] (7) to [bend left=30] (8) to [bend left=30] cycle;
		\foreach \X in {1,2,3,4}
		{
		\tkzDrawPoints[rectangle,fill =red,size=6,color=red](\X);
		};
		
		\foreach \X in {0,5,6,7}
		{
		\tkzDrawPoints[size=6,orange,osq](\X);
		};
		\foreach \X in {8}
		{
		\tkzDrawPoints[size=6,darkpink,psq](\X);
		};
		
		\begin{scope}[xshift=2.5ex]
		\tkzDefPoint(-2.8,-0.5){gammaM};
		\tkzLabelPoint[blue](gammaM){\Large $\delta$}
		\tkzDefPoint(-1.5,-1){gammaP};
		\tkzLabelPoint[red](gammaP){\Large $\pmb{\vartheta}$}
		\tkzDefPoint(0.3,0.5){v};
		\tkzLabelPoint[red](v){\Large $\eta_R$}
		\tkzDefPoint(2.5,-0.3){w};
		\tkzLabelPoint[red](w){\Large $C_R$}
		\tkzDefPoint(1.3,0){deltav};
		\tkzLabelPoint[mypurple](deltav){\Large $\delta_R^\circ$}
		\end{scope}
    \end{tikzpicture}
\caption{\label{fig:resclo1reduce} Reduction to the case where $\delta$ is contained in a cell of $\pmb{\Gamma}(\Prj(\Delta^{\gpoint}))$.}
\end{figure} Then $\MM(\vartheta)$ corresponds to one summand of $\PP_{\MM(\delta)}^0$ in the minimal projective resolution $(\PP_M^\bullet)$. In the kernel of the epimorphism $\PP_{\MM(\delta)}^0 \twoheadrightarrow \MM(\delta)$, there is an indecomposable summand $K$ such that $\gamma_{(K)} = \delta_R^\circ$, by \cref{prop:geo_mor}. As $\Res(\MM(\delta))$ is closed under kernels of epimorphisms, then $\MM(\delta_R^\circ) \in \Res(\MM(\delta))$. By dual arguments, we also get that $\MM(\delta_L^\circ) \in \Res(\MM(\delta))$.
\end{proof}

Finally, we get the following result.

\begin{prop} \label{prop:UMinResM}
    Let $(Q,R)$ be a gentle tree, and $\Surf(Q,R) = (\pmb{\Sigma}, \mathcal{M}, \Delta^{\gpoint})$. For any $\delta \in \Accord'$, the category $\mathscr{U}_\delta$ is a subcategory of $\Res(\MM(\delta))$.
\end{prop}

\begin{proof}
To prove the result, as $\Res(\MM(\delta))$ and $\mathscr{U}_{\delta}$ are additive, we only have to prove that indecomposable representations in $\mathscr{U}_\delta$ are in $\Res(\MM(\delta))$. Therefore, we will use the geometric model to do so. 

Let $\Surf(Q,R) = (\pmb{\Sigma}, \mathcal{M}, \Delta^{\gpoint})$. As $(Q,R)$ is a gentle tree, by \cref{thm:surftree}, $\Sigma$ is homeomorphic to a disc, and $\mathcal{M} \subset \partial \pmb{\Sigma}$.
If $M \in \proj(Q,R)$, then we have that $\opResAc(\delta) = \Prj(\Delta^{\gpoint})$ and we are done. In the remainder of the proof, we assume that $\MM(\delta) \in \pmb{\ind \setminus \proj}(Q,R)$. 

Either $\delta$ is contained in a cell of $\pmb{\Gamma}(\Prj(\Delta^{\gpoint}))$ or not. By \cref{lem:UMinResinacell,lem:UMinRes}, we get that $\opResAc(\delta) \subset \{\gamma_{(E)} \mid E \in \ind(\Res(\MM(\delta)))\}$. Hence the result.
\end{proof}

\begin{proof}[Proof of \cref{thm:resclo1}]
This result follows from \cref{prop:UMinResM,prop:UMRes}.
\end{proof}

\begin{remark}
    \cref{thm:resclo1} offers a way to consider the poset of non-projective indecomposable objects with Res-relation geometrically. Indeed  accordions $\delta$ and $\eta$, such that $\MM(\eta)\Resleq\MM(\delta)$, are such that $\eta\in\opResAc(\delta)$.
\end{remark}
	
	\section*{Acknowledgements}
	
    B.D. thanks the \emph{Institut des Sciences Mathematiques (UQAM)} and the \emph{Engineering and Physical Sciences Research Council (EP/W007509/1)} for their partial funding support. The authors acknowledge the \emph{CHARMS program grant (ANR-19-CE40-0017-02)} for their funding support.
    
	The authors, together, thank Yann Palu and Baptiste Rognerud for the various discussions they had and their advice, all along this work.
    The authors also thank the anonymous referee for his extensive review and helpful comments that greatly improved this work.
	
	\bibliography{Article4_2}
	\bibliographystyle{alpha}
\end{document}